\documentclass[reqno]{amsart}
\usepackage[margin=1in]{geometry}
\usepackage{microtype}
\usepackage[T1]{fontenc}
\usepackage{lmodern}
\usepackage{amsmath,amssymb,amsthm,mathtools}
\usepackage{graphicx}
\usepackage[x11names,dvipsnames]{xcolor}
\usepackage[colorlinks]{hyperref}
\hypersetup{
  colorlinks=true,
  linkcolor=blue,
  citecolor=teal,
  filecolor=red,
  urlcolor=Plum
}
\newcommand{\R}{\mathbb{R}}
\newcommand{\Z}{\mathbb{Z}}
\newcommand{\G}{\mathcal{G}}
\newcommand{\NN}{\mathcal{N}}
\newcommand{\Hc}{\mathcal{H}}
\newcommand{\Lc}{\mathcal{L}}
\newcommand{\K}{\mathcal{K}}
\newcommand{\D}{\mathcal{D}}
\newcommand{\DiffI}{\operatorname{Diff}_1}
\newcommand{\AutI}{\operatorname{Aut}_1}
\newcommand{\CoderII}{\operatorname{Coder}_{\ge2}}
\newcommand{\hh}{\mathfrak{h}}
\newcommand{\pp}{\mathfrak{p}}
\newcommand{\X}{\mathfrak{X}}
\newcommand{\Of}{\widehat{\mathcal{O}}}
\newcommand{\Ic}{\mathcal{I}}
\newcommand{\Linf}{L_\infty}
\newcommand{\Pinf}{P_\infty}
\newcommand{\Aff}{\operatorname{Aff}^{+}(\R)}
\newcommand{\aff}{\mathfrak{aff}(\R)}
\newcommand{\sss}{\mathsf{s}}
\newcommand{\ttt}{\mathsf{t}}
\newcommand{\tto}{\rightrightarrows}
\newcommand{\too}{\longrightarrow}
\newcommand{\defn}[1]{\textbf{#1}}
\newcommand{\dd}[1]{\frac{\partial}{\partial #1}}
\DeclareMathOperator{\ad}{ad}
\DeclareMathOperator{\Coder}{Coder}
\DeclareMathOperator{\MultiDer}{MultiDer}

\DeclareMathOperator{\pr}{pr}
\DeclareMathOperator{\id}{id}

\DeclareMathOperator{\Lie}{Lie}
\DeclareMathOperator{\Iso}{Iso}
\DeclareMathOperator{\im}{im}
\DeclareMathOperator{\sgn}{sign}

\DeclareMathOperator{\res}{res}
\DeclareMathOperator{\BCH}{BCH}
\DeclareMathOperator{\divv}{div}
\theoremstyle{plain}
\newtheorem{introtheorem}{Theorem}

\newtheorem{proposition}{Proposition}[section]
\newtheorem{lemma}[proposition]{Lemma}
\newtheorem{corollary}[proposition]{Corollary}
\theoremstyle{definition}
\newtheorem{definition}[proposition]{Definition}
\newtheorem{example}[proposition]{Example}
\newtheorem{remark}[proposition]{Remark}
\numberwithin{equation}{section}
\begin{document}

\title{Shifted Lagrangian geometry of curves in Radko surfaces}

\author[Yassine Ait Mohamed]{Yassine Ait Mohamed}
\address{D\'epartement de math\'ematiques, Universit\'e de Sherbrooke, 2500 Bd de l'Universit\'e, Sherbrooke, QC, J1K 2R1, Canada}
\email{yassine.ait.mohamed@usherbrooke.ca}

\author[Maxence Mayrand]{Maxence Mayrand}
\address{D\'epartement de math\'ematiques, Universit\'e de Sherbrooke, 2500 Bd de l'Universit\'e, Sherbrooke, QC, J1K 2R1, Canada}
\email{maxence.mayrand@usherbrooke.ca}

\begin{abstract}
We describe the shifted Lagrangian geometry of curves in Radko surfaces, through their global integrations and their homotopy Poisson algebras. We construct the adjoint symplectic groupoid of the surface explicitly, show that admissible curves integrate to embedded Lagrangian subgroupoids, which are $1$-shifted Lagrangians, and characterize when their shifted intersections are clean. For a component of the degeneracy locus, a canonical Hamiltonian flow on its shifted self-intersection recovers the modular period. On the infinitesimal side, the homotopy Poisson algebra of a connected admissible curve not contained in the degeneracy locus is $\Linf$-isomorphic to its underlying cochain complex, and its $\Linf$-isomorphism class is determined by the total contact multiplicity. The brackets of arity at least two can be removed by a change of tubular neighbourhood exactly when all contacts are transverse. At tangential contacts, the binary bracket retains the modular vector field. We establish this distinction using a general characterization of the $\Linf$-isomorphisms induced by changes of tubular neighbourhood, valid for arbitrary coisotropic submanifolds.
\end{abstract}

\maketitle

\raggedbottom

\tableofcontents

\section{Introduction}\label{RADK854RZZ}
A coisotropic submanifold of a Poisson manifold admits two complementary descriptions. Globally, its conormal Lie algebroid can, under suitable integrability assumptions, be integrated to a Lagrangian subgroupoid of a symplectic groupoid integrating the ambient manifold \cite{crainicfernandes2004,cattaneo2004integration,cattaneofelder2004branes}.
Infinitesimally, its Chevalley--Eilenberg complex carries the homotopy Poisson structure of Cattaneo and Felder \cite{Cattaneo2007}. These descriptions lead to different questions: whether the integrating subgroupoid is embedded, how it intersects other Lagrangian subgroupoids, and which features of the original geometry are retained by the homotopy Poisson algebra.

Shifted geometry provides a common setting for these questions \cite{ptvv2013,cptvv2017,safronov2017,melanisafronov2018,melanisafronov2018II}. A symplectic groupoid presents a $1$-shifted symplectic differentiable stack, and a Lagrangian subgroupoid presents a $1$-shifted Lagrangian morphism \cite{xu2004momentum,bcwz2004}. In the differential-geometric theory developed by the second author, a clean homotopy fibre product of two such Lagrangians carries a $0$-shifted Poisson structure whenever the induced Dirac structure is smooth \cite{Maxence2025}. In previous work, we extended the integration and quantization picture to $1$-shifted coisotropics under suitable hypotheses, constructing a flat $\Pinf$-algebra and a curved $A_\infty$-quantization \cite{aitmohamedmayrand2026}. Explicit examples require an understanding of the integrating subgroupoids, the cleanness of their homotopy fibre products, and the dependence of the higher brackets on choices.

In this paper, we study these questions for curves in Radko surfaces. A Radko surface is a compact oriented Poisson surface whose bivector vanishes transversely along a nonempty union of circles $D$ \cite{radko2002classification}. We give an explicit description of its adjoint symplectic groupoid, relating the affine local model and the global gluing to the constructions of Gualtieri and Li \cite{gualtieri2014symplectic}. Every curve in the surface is coisotropic, and its Lagrangian integrations are $1$-shifted Lagrangians, but its contacts with $D$ produce different behaviour: a curve may cross $D$ transversely, have tangential contacts, or coincide with a component of $D$. This setting allows us to compute both the global integrations and the homotopy Poisson algebras explicitly.

Our results describe how these contacts govern the integrating Lagrangian subgroupoids, their shifted intersections, and the information retained by their homotopy Poisson algebras. In particular, we distinguish isomorphisms of the underlying $\Linf$-algebras from those induced by changes of tubular neighbourhood. This distinction is already visible on the Poisson sphere: a meridian and a circle tangent to the degeneracy equator have $\Linf$-isomorphic homotopy Poisson algebras, but only the meridian admits a tubular neighbourhood in which all brackets of arity at least two vanish (Example \ref{RADK207TNC}). We first state the general criterion that makes this distinction precise, before describing the global and infinitesimal results for curves.

Let $N\subset M$ be coisotropic, and write $\mathfrak{a}=C^\infty(\wedge^\bullet\nu_N)$. A tubular neighbourhood $\psi$ of $N$, that is, an open embedding into $M$ of a neighbourhood of the zero section of $\nu_N$ which is the identity on $N$ and induces the identity on $\nu_N$, determines a flat $\Pinf$-structure $\lambda^\psi$ on $\mathfrak{a}$, and we call $(\mathfrak{a},\lambda^\psi)$ the \emph{homotopy Poisson algebra} of $(N,\psi)$. Its unary bracket $\lambda_1$ is independent of $\psi$, and it is the Chevalley--Eilenberg differential of $TN^\circ$ up to sign. The higher brackets are obtained from the formal expansion $\sigma_\psi$ of $\sigma$ along $N$ in the tubular neighbourhood $\psi$, as higher derived brackets \cite{Cattaneo2007,voronov2005higher,cattaneo2008equivalences}.

For the statement of the criterion, recall that an $\Linf$-isomorphism is a coalgebra automorphism $F$ of $S(\mathfrak{a}[1])$ intertwining the coderivations encoding the brackets. Its linear part is the first Taylor coefficient $F_1$. A formal change of tubular neighbourhood is a formal diffeomorphism of $\nu_N$ along $N$ that induces the identity on $N$ and on $\nu_N$. The construction recalled in \S\ref{RADK088RYM} associates to it a coalgebra automorphism with identity linear part, and an $\Linf$-isomorphism whenever the change carries one formal Poisson bivector to the other. We call such an isomorphism \emph{induced by a change of tubular neighbourhood}.

When $F_1=\id$, its logarithm is the well-defined coderivation $\log F=\sum_{k\ge1}\frac{(-1)^{k+1}}{k}\,(F-\id)^k$. The series is finite on each element of $S(\mathfrak{a}[1])$. A \emph{multiderivation} is a graded symmetric multilinear map on $\mathfrak{a}[1]$ satisfying a Leibniz rule in each argument for the product of $\mathfrak{a}$. Our first result gives an intrinsic criterion for an isomorphism to be induced by a formal change of neighbourhood.

\begin{introtheorem}\label{RADK919BTK}
Let $\psi,\psi'$ be tubular neighbourhoods of a coisotropic submanifold of a Poisson manifold. An $\Linf$-isomorphism $F\colon(\mathfrak{a},\lambda^\psi)\to(\mathfrak{a},\lambda^{\psi'})$ is induced by a change of tubular neighbourhood if and only if $F_1=\id$ and the Taylor coefficients of positive arity of $\log F$ are multiderivations of $\mathfrak{a}$.
\end{introtheorem}

The inducing formal change is unique. The logarithm is essential, since the Taylor coefficients of $F$ itself need not be multiderivations. The theorem distinguishes abstract linearization of a homotopy Poisson algebra from linearization by a change of tubular neighbourhood. For coisotropic submanifolds of symplectic manifolds, Sch\"atz and Zambon compare geometric equivalences with gauge equivalences of the associated $\Linf$-algebras \cite{schatzzambon2017}. Its proof, together with the comparison with the isomorphisms of Cattaneo and Sch\"atz \cite{cattaneo2008equivalences}, occupies \S\ref{RADK088RYM}.

We now return to a Radko surface $(\Sigma,\sigma)$ with degeneracy locus $D$. These are precisely the compact oriented log symplectic surfaces with nonempty degeneracy locus \cite{gmp2014,gualtieri2014symplectic}. Near a component $D_j$ of $D$ there are coordinates $(t,r)\in S^1_{\ell_j}\times(-c_j,c_j)$ in which
\[
\sigma=r\,\dd r\wedge\dd t,\quad S^1_{\ell_j}=\R/\ell_j\Z.
\]
The positive number $\ell_j$ is the modular period of $D_j$, the period of the modular vector field $X_{\mathrm{vol}}$ \cite{weinstein1997modular} of $\sigma$ along $D_j$ \cite{radko2002classification}. Together with the oriented degeneracy locus and a regularized volume, the modular periods classify Radko structures \cite{radko2002classification}. 

In a symplectic surface, the questions above have simple answers. The adjoint symplectic groupoid of a connected symplectic surface is its pair groupoid, two curves always have a clean homotopy fibre product, and every curve has a tubular neighbourhood whose homotopy Poisson algebra has no brackets of arity at least two. Along the degeneracy locus of a Radko surface, the rank of $\sigma$ drops from $2$ to $0$, so that a curve can meet $D$ transversely, tangentially to any order, or contain a component of $D$.

Let $\G_\Sigma$ denote the adjoint symplectic groupoid integrating $(\Sigma,\sigma)$. In \S\ref{RADK261MGJ}, we give explicit formulas for this groupoid and its symplectic form, using the affine local model and the integration constructions of Gualtieri and Li \cite{gualtieri2014symplectic}. Over each component of $\Sigma\setminus D$ it is the pair groupoid, and its isotropy group at a point of $D$ is the orientation-preserving affine group $\Aff=\{(a,b):a>0,\ b\in\R\}$, with $(a_2,b_2)(a_1,b_1)=(a_2a_1,\,b_1+a_1b_2)$. Near $D_j$, the groupoid $\G_\Sigma$ is a quotient of the action groupoid of $\Aff$ acting on $S^1_{\ell_j}\times\R$ by $(a,b)\cdot(t,r)=(t+br,\,ar)$, by a discrete subgroupoid of loops which records the modular period (\S\S\ref{RADK265ACT} and \ref{RADK898CDG}). We call $\{(1,b):b\in\R\}$ the \emph{translation subgroup} of $\Aff$, and the conjugates of $\{(a,0):a>0\}$ its \emph{dilation subgroups}.

Throughout, a curve is a closed embedded one-dimensional submanifold. We call it \emph{admissible} if each component is contained in $\Sigma\setminus D$, is a component of $D$, or meets $D$ in finitely many points. We say that $N$ is \emph{tangent to $D$} at $x\in N\cap D$ if $T_xN=T_xD$, which holds in particular if $N$ contains the component of $D$ through $x$. At a point of $N\cap D$, a component not contained in $D$ has a tangential contact with $D$ if it is tangent to $D$, and a transverse contact otherwise. We do not consider curves which meet $D$ in infinitely many points without containing a component of it, for which $N\cap D$ can be an arbitrary closed subset of $D_j$. Our second result integrates admissible curves inside the adjoint groupoid itself, and not only in the source-simply-connected one.

\begin{introtheorem}\label{RADK953BQH}
Let $N$ be an admissible curve in a Radko surface. Then $TN^\circ$ integrates to an embedded Lagrangian subgroupoid $\Lc\tto N$ of $\G_\Sigma$, whose isotropy group at $p\in N\cap D$ is the translation subgroup if $N$ is tangent to $D$ at $p$, and a dilation subgroup otherwise.
\end{introtheorem}

The subgroupoid $\Lc$ is source-connected, two points of $N\setminus D$ are joined by an arrow of $\Lc$ exactly when they lie in the same component of $N\setminus D$, and the arrows with source in $N\cap D$ are loops. In general, $TN^\circ$ only integrates to an immersed subgroupoid, the image of a morphism which need not be injective \cite{moerdijkmrcun2006}, and the Lagrangian subgroupoid of \cite{cattaneo2004integration} lives in the source-simply-connected integration. Here it is embedded in the adjoint integration, which is Hausdorff, so that $\Lc$ is a $1$-shifted Lagrangian in an explicit Lie groupoid and can be intersected with others.

Given two admissible curves $N_1$ and $N_2$, their subgroupoids $\Lc_1$ and $\Lc_2$ define a homotopy fibre product $\Lc_1\times^h_{\G_\Sigma}\Lc_2$, which is a groupoid over the space of triples $(x,g,y)$ with $x\in N_1$, $y\in N_2$ and $g$ an arrow of $\G_\Sigma$ from $x$ to $y$. It is clean if this space is a submanifold with the expected tangent spaces, and this is the condition under which the general theory produces a $0$-shifted Poisson structure, since the relevant Dirac structure is then automatically smooth for $1$-shifted Lagrangians. Cleanness is weaker than transversality, so this structure is in general Poisson and not symplectic. The cleanness of $\Lc_1\times^h_{\G_\Sigma}\Lc_2$ is controlled by the contacts of $N_1$ and $N_2$ with $D$ at their common points.

\begin{introtheorem}\label{RADK367ESC}
Let $N_1$ and $N_2$ be admissible curves in a Radko surface. Then $\Lc_1\times^h_{\G_\Sigma}\Lc_2$ is clean if and only if, at every point of $N_1\cap N_2\cap D$, either both curves contain the component of $D$ through this point, or exactly one of them is tangent to $D$.
\end{introtheorem}

This condition is not ordinary transversality of the curves in $\Sigma$: two transverse great circles on the Poisson sphere of Example \ref{RADK748QSA} can have a non-clean homotopy fibre product. For a single curve, the shifted self-intersection $\Lc\times^h_{\G_\Sigma}\Lc$ is clean if and only if every component of $N$ which meets $D$ is a component of $D$ (Corollary \ref{RADK784SLF}).

When it is clean, the homotopy fibre product is Morita equivalent to a disjoint union of finitely many points and, for each common component $D_j$ of $D$, of the trivial bundle of groups $(S^1_{\ell_j}\times\R)\times\R$ over a cylinder (Proposition \ref{RADK741QHM}). The $0$-shifted Poisson structure of \cite[Corollary 4.12]{Maxence2025} restricts to a symplectic form on this cylinder, for which the Hamiltonian vector field of a canonical function projects to the modular vector field of $\sigma$ along $D_j$. The modular period is thus recovered as the period of this flow, although the symplectic cylinder alone does not determine it. The regularized volume does not enter, in agreement with the result of Bursztyn and Radko that, for a fixed oriented degeneracy locus, the modular periods determine a Radko structure up to gauge equivalence and Poisson diffeomorphism, and hence up to Morita equivalence \cite[Theorems 6.2 and 6.3]{bursztynradko2003}.

We next compute the homotopy Poisson algebra of a component $N_c$ of an admissible curve. The algebra $\mathfrak{a}_c=C^\infty(N_c)\oplus C^\infty(\nu_{N_c})$ is concentrated in degrees $0$ and $1$, and we write $(\mathfrak{a}_c,\lambda_1)$ for this algebra with $\lambda_1$ as its only bracket. In an adapted chart $(x,y)$ for a tubular neighbourhood, write $N_c=\{y=0\}$, let $e$ be the frame of $\nu_{N_c}$ dual to $y$, and write
\[
\sigma_\psi=\left(\sum_{k\ge0}\upsilon_k(x)y^k\right)\dd x\wedge\dd y.
\]
The higher brackets are expressed directly in terms of the coefficients $\upsilon_k$. Up to graded antisymmetry, the only possibly nonzero brackets are
\[
\lambda^\psi_n(f,g_1e,\dots,g_{n-1}e)
=(n-1)!\,\upsilon_{n-1}\,\frac{\partial f}{\partial x}\,g_1\cdots g_{n-1}\,e,
\]
where $f,g_1,\dots,g_{n-1}\in C^\infty(N_c)$ and $n\ge1$.
In particular, $\lambda_1$ is determined by $\upsilon_0$ and the binary bracket by $\upsilon_1$.

Suppose first that $N_c$ is not contained in $D$. At each point $p\in N_c\cap D$, the \emph{contact multiplicity} $m_p\in\{1,2,\dots\}\cup\{\infty\}$ is the order of vanishing on $N_c$ of a local defining function of $D$, so that transverse contacts have multiplicity one. These multiplicities govern the cohomology of $\lambda_1$ and the $\Linf$-isomorphism class of the homotopy Poisson algebra. Here $\Linf$-isomorphisms are understood algebraically: no continuity, locality or compatibility with the commutative products is required. The homotopy Poisson algebras carry their products throughout, but the abstract classification concerns their underlying $\Linf$-algebras. An $\Linf$-isomorphism $(\mathfrak{a}_c,\lambda^\psi)\to(\mathfrak{a}_c,\lambda_1)$ is induced by a change of tubular neighbourhood if $\lambda_1$ is the homotopy Poisson structure of a tubular neighbourhood $\psi''$ and the isomorphism is induced by a change from $\psi$ to $\psi''$ as above.

\begin{introtheorem}\label{RADK420ZBG}
Let $N_c$ be a component of an admissible curve in a Radko surface, with $N_c\not\subset D$, and let $\psi$ be a tubular neighbourhood of $N_c$.
\begin{enumerate}
\item The algebra $(\mathfrak{a}_c,\lambda^\psi)$ is $\Linf$-isomorphic to $(\mathfrak{a}_c,\lambda_1)$. Two such components have $\Linf$-isomorphic algebras if and only if they have the same total multiplicity $\sum_pm_p\in\{0,1,2,\dots\}\cup\{\infty\}$.
\item An $\Linf$-isomorphism $(\mathfrak{a}_c,\lambda^\psi)\to(\mathfrak{a}_c,\lambda_1)$ can be chosen to be induced by a change of tubular neighbourhood if and only if all contacts of $N_c$ with $D$ are transverse. For every tubular neighbourhood $\psi$, the binary bracket is $\lambda^\psi_2(f,v)=X(f)\,v$ for a vector field $X$ on $N_c$ with $X(p)=X_{\mathrm{vol}}(p)$ at every tangential contact $p$.
\end{enumerate}
\end{introtheorem}

Here $\dim H^1(\mathfrak{a}_c,\lambda_1)=1+\sum_pm_p$ when all multiplicities are finite, and $H^1$ has the dimension of the continuum otherwise. Every component has a tubular neighbourhood in which $\lambda_n=0$ for $n\ge3$, so that its homotopy Poisson algebra is a differential graded Poisson algebra. When all contacts are transverse, $|N_c\cap D|=2k$, and in a suitable tubular neighbourhood, coordinate $\theta\in\R/2\pi\Z$ and frame $e$, the only nonzero bracket is $\lambda_1(f)=\sin(k\theta)\,\frac{\partial f}{\partial\theta}\,e$, or $\lambda_1(f)=\frac{\partial f}{\partial\theta}\,e$ if $k=0$. This agrees with the Mazzeo--Melrose formula for $b$-manifolds \cite[Theorem 27]{gmp2014} and with the $b$-symplectic normal form of Geudens and Zambon \cite{GeudensZambon2021}. At a tangential contact of finite multiplicity $m_p\ge2$, the jet of order $m_p-2$ of $X$ at $p$ is intrinsic, whereas its jet of order $m_p-1$ depends on the tubular neighbourhood.

A component $N_c=D_j$ has different behaviour. Deformations of such Lagrangian submanifolds contained in the degeneracy locus are studied by Geudens and Zambon \cite{geudenszambon2020lagrangian}. Its differential $\lambda_1$ vanishes, its binary bracket is $\lambda_2(f,v)=X_{\mathrm{vol}}(f)\,v$ for every tubular neighbourhood, and every flat $\Linf$-algebra $\Linf$-isomorphic to $(\mathfrak{a}_c,\lambda^\psi)$ has a nonzero binary bracket (Propositions \ref{RADK890FAH} and \ref{RADK667KEG}). All these algebras are $\Linf$-isomorphic to each other, whatever the modular period. The period is nevertheless recovered by every $\Linf$-isomorphism whose linear part is compatible with the products, as the period of $X_{\mathrm{vol}}|_{D_j}$ (Proposition \ref{RADK669PER}).

For a connected admissible curve not contained in $D$, the $\Linf$-isomorphism class of the homotopy Poisson algebra retains only a numerical invariant, the total contact multiplicity, which also governs the cohomology of $TN^\circ$ and, through the van Est map, the differentiable cohomology of the integrating subgroupoid when the curve meets $D$ (Remark \ref{RADK745VEP}). The geometric isomorphisms of Theorem \ref{RADK919BTK} retain more, namely the modular vector field at the tangential contacts, which is invisible to the abstract classification but obstructs geometric linearization. The global structures retain the modular periods, through the canonical function on the shifted self-intersections of the degeneracy curves, whose Hamiltonian flow has period $\ell_j$.

\subsection*{Organization of the paper}
In \S\ref{RADK088RYM}, we recall the homotopy Poisson algebra of a coisotropic submanifold, study changes of tubular neighbourhood and prove Theorem \ref{RADK919BTK}. In \S\ref{RADK261MGJ}, we describe the adjoint symplectic groupoid of a Radko surface. In \S\ref{RADK793NZP}, we construct the Lagrangian subgroupoids of admissible curves and prove Theorem \ref{RADK953BQH}. In \S\ref{RADK389HFP}, we determine the clean homotopy fibre products of these subgroupoids and their induced structures, and prove Theorem \ref{RADK367ESC}. Finally, in \S\ref{RADK035TQA}, we compute and classify the homotopy Poisson algebras of curves and prove Theorem \ref{RADK420ZBG}.

\subsection*{Acknowledgments}
Y.A.M.\ thanks the Institut des sciences math\'ematiques (ISM) for financial support. M.M.\ acknowledges support from a Discovery Grant (RGPIN-2023-04587) from the Natural Sciences and Engineering Research Council of Canada (NSERC).

\section{Homotopy Poisson algebras and changes of tubular neighbourhood}\label{RADK088RYM}

In this section, we characterize the $\Linf$-isomorphisms between homotopy Poisson algebras of a coisotropic submanifold that are induced by changes of tubular neighbourhood, and prove Theorem \ref{RADK919BTK}. Let $(M,\sigma)$ be a Poisson manifold and let $N\subset M$ be a coisotropic submanifold. We write
\[
\nu=TM|_N/TN,
\quad
\mathfrak{a}=C^\infty\bigl(\wedge^\bullet\nu\bigr),
\quad
\mathfrak{a}[1]^k=\mathfrak{a}^{k+1}.
\]
A tubular neighbourhood of $N$ determines a homotopy Poisson structure on $\mathfrak{a}$. We first recall the coalgebra operations needed to take logarithms of isomorphisms, then construct the homotopy Poisson brackets by higher derived brackets. Finally, we compare formal changes of neighbourhood with coalgebra automorphisms.

\subsection{Coalgebras and \texorpdfstring{$\infty$}{infinity}-isotopies}\label{RADK418QWN}

All vector spaces are real and $\Z$-graded, and all tensor operations use the Koszul sign rule. For homogeneous operators, we use the graded commutator $[T_1,T_2] =T_1\circ T_2-(-1)^{|T_1||T_2|}\,T_2\circ T_1$.
For a graded vector space $V$, let
\[
S(V)=\bigoplus_{n\ge0}S^nV,
\quad
S^{\le n}V=\bigoplus_{k=0}^{n}S^kV.
\]
We equip $S(V)$ with its usual coproduct $\Delta$, counit and coaugmentation. It is the cofree conilpotent cocommutative coalgebra on $V$ \cite{ladamarkl1995,lodayvallette2012}. An element of $S^nV$ has \defn{length} $n$.

Coalgebra morphisms will have degree zero and preserve the counit and coaugmentation. Such a morphism $F\colon S(V)\to S(W)$ and a homogeneous coderivation $Q$ of $S(V)$ are determined by their \defn{Taylor coefficients} $F_n=\pr_W\circ F|_{S^nV}$, $Q_n=\pr_V\circ Q|_{S^nV}$.
Here $F_0=0$, whereas a coderivation may have a nonzero coefficient $Q_0\in V$. Conversely, every family of coefficients of the appropriate degree, with $F_0=0$ in the morphism case, determines such a map \cite[\S2 and Lemma 2.4]{ladamarkl1995}, \cite[\S4.1]{kontsevich2003}. Explicitly,
\begin{equation}\label{RADK170LCB}
F(v_1\cdots v_n)
=\sum_{\mathcal B}\varepsilon(\mathcal B)\,
F_{\#B_1}(v_{B_1})\cdots F_{\#B_k}(v_{B_k}),
\end{equation}
where $\mathcal B=\{B_1,\dots,B_k\}$ runs over the partitions of $\{1,\dots,n\}$ into nonempty blocks. The factors in $v_B$ are ordered by their indices, and $\varepsilon(\mathcal B)$ is the Koszul sign of the corresponding permutation. Similarly,
\begin{equation}\label{RADK450GWZ}
Q(v_1\cdots v_n)
=\sum_{I\subset\{1,\dots,n\}}\varepsilon(I)\,
Q_{\#I}(v_I)\,v_{I^c},
\end{equation}
where $\varepsilon(I)$ is the Koszul sign of $v_1\cdots v_n\mapsto v_Iv_{I^c}$, and the term indexed by $I=\varnothing$ uses $Q_0=Q_0(1)$.

We use the symmetric convention: an \defn{$\Linf$-structure} on $V$ is a coderivation $Q$ of degree $1$ with $Q\circ Q=0$, and it is \defn{flat} if $Q_0=0$. An \defn{$\Linf$-morphism} $F\colon(V,Q)\to(V',Q')$ of flat $\Linf$-structures is a coalgebra morphism with $F_0=0$ and $F\circ Q=Q'\circ F$. It is an \defn{$\Linf$-isomorphism} if it is bijective, and its \defn{linear part} is $F_1$ \cite[\S2.1, Definition 2.2 and \S3.1]{cattaneo2008equivalences}.

When applying this convention to a graded algebra $\mathfrak{a}$, the coderivation acts on $S(\mathfrak{a}[1])$. In the rest of this subsection, $\mathfrak{a}$ may be any graded commutative algebra. We write $\deg b$ for the degree of a homogeneous element $b\in\mathfrak{a}$ and $|b|=\deg b-1$ for its degree in $\mathfrak{a}[1]$.

\begin{definition}\label{RADK013YGA}
Let $n\ge1$. A graded symmetric map $l\colon S^n(\mathfrak{a}[1])\to \mathfrak{a}[1]$ of degree $d$ is a \defn{multiderivation} of $\mathfrak{a}$ if, for all homogeneous $b,c\in \mathfrak{a}$ and $a_2,\dots,a_n\in \mathfrak{a}[1]$,
\begin{equation}\label{RADK692LBY}
\begin{aligned}
l(bc,a_2,\dots,a_n)
&=(-1)^{\deg c\,(|a_2|+\cdots+|a_n|)}\,
l(b,a_2,\dots,a_n)\,c\\
&\quad+(-1)^{d\,\deg b}\,b\,l(c,a_2,\dots,a_n).
\end{aligned}
\end{equation}
A flat $\Linf$-structure $Q$ on $\mathfrak{a}[1]$ is a \defn{flat $\Pinf$-structure} on $\mathfrak{a}$ if its Taylor coefficients $Q_n$, $n\ge1$, are multiderivations.
\end{definition}

By graded symmetry, a multiderivation satisfies a Leibniz rule in each argument,
\begin{equation}\label{RADK049HQH}
l(\dots,bc,\dots)=\pm\,l(\dots,b,\dots)\,c\pm b\,l(\dots,c,\dots),
\end{equation}
with signs depending only on the degrees. For $\mathfrak{a}=C^\infty(\wedge^\bullet\nu)$, a multiderivation is therefore determined, over a chart with coordinates $x^\mu$ on $N$ and a frame $(e_i)$ of $\nu$, by its values on the $x^\mu$ and the $e_i$.

Let $\AutI(\mathfrak{a})$ be the set of degree-zero coalgebra morphisms $F\colon S(\mathfrak{a}[1])\to S(\mathfrak{a}[1])$ with $F_0=0$ and $F_1=\id$, called \defn{$\infty$-isotopies} \cite[Chapter 3]{dsv2023}, and $\CoderII(\mathfrak{a})$ the Lie algebra of coderivations $\delta$ of $S(\mathfrak{a}[1])$ of degree $0$ with $\delta_0=\delta_1=0$. By \eqref{RADK170LCB} and \eqref{RADK450GWZ}, $F-\id$ and $\delta$ map $S^n(\mathfrak{a}[1])$ into $S^{\le n-1}(\mathfrak{a}[1])$ and vanish on $S^{\le1}(\mathfrak{a}[1])$, so
\begin{equation}\label{RADK818ULM}
(F-\id)^n|_{S^{\le n}(\mathfrak{a}[1])}=0,\quad \delta^n|_{S^{\le n}(\mathfrak{a}[1])}=0 .
\end{equation}
Hence the series $\log F=\sum_{k\ge1}\frac{(-1)^{k+1}}{k}(F-\id)^k$ and $e^\delta=\sum_{k\ge0}\frac1{k!}\delta^k$ are finite on each element.

The set $\AutI(\mathfrak{a})$ is a group \cite[Chapter 3]{dsv2023}. To see that $\log F$ is a coderivation, write
\begin{align*}
\Delta\circ(F-\id)^k
&=(F\otimes F-\id\otimes\id)^k\circ\Delta,\\
\log(F\otimes F)
&=\log(F\otimes\id)+\log(\id\otimes F).
\end{align*}
The second identity follows because $F\otimes\id$ and $\id\otimes F$ commute, and all series are finite on each length filtration step. It follows that
\[
\Delta\circ\log F
=\bigl((\log F)\otimes\id+\id\otimes(\log F)\bigr)\circ\Delta.
\]
Thus $\log F\in\CoderII(\mathfrak{a})$. Conversely, for $\delta\in\CoderII(\mathfrak{a})$, the coderivation identity gives $\Delta\circ e^\delta=(e^\delta\otimes e^\delta)\circ\Delta$.
Since $e^\delta$ is the identity on $S^{\le1}(\mathfrak{a}[1])$, it belongs to $\AutI(\mathfrak{a})$. The identities $\exp\circ\log=\id$ and $\log\circ\exp=\id$ therefore give a bijection
\begin{equation}\label{RADK042DFG}
\log\colon\AutI(\mathfrak{a})\longrightarrow\CoderII(\mathfrak{a}),
\quad
\log^{-1}=\exp.
\end{equation}

Taking the components of length one in $\log F=(F-\id)-\frac12(F-\id)^2+\cdots$, \eqref{RADK170LCB} gives, for $F\in\AutI(\mathfrak{a})$,
\begin{align}
(\log F)_2(a_1,a_2)&=F_2(a_1,a_2),\label{RADK675CTY}\\
(\log F)_3(a_1,a_2,a_3)&=F_3(a_1,a_2,a_3)-\tfrac12\sum\varepsilon\,F_2\big(F_2(a_i,a_j),a_k\big),\label{RADK504RBS}
\end{align}
where the sum runs over the decompositions $\{1,2,3\}=\{i<j\}\sqcup\{k\}$, and $\varepsilon$ is the Koszul sign of $a_1a_2a_3\mapsto a_ia_ja_k$.

Since all series are finite on each $S^{\le n}(\mathfrak{a}[1])$, identities of formal power series give, for $\delta,\delta'\in\CoderII(\mathfrak{a})$ and a coderivation $Q$ with $Q_0=0$,
\begin{equation}\label{RADK496NHH}
\begin{aligned}
e^\delta\circ Q\circ e^{-\delta}
&=\sum_{k\ge0}\frac1{k!}\,\ad_\delta^k(Q),\\
e^\delta\circ e^{\delta'}
&=e^{\BCH(\delta,\delta')},
\end{aligned}
\end{equation}
where $\BCH$ is the Baker--Campbell--Hausdorff series, and the series converge coefficientwise, since $\ad_\delta$ raises the lowest non-vanishing arity by at least one.

\subsection{The homotopy Poisson algebra of a coisotropic submanifold}\label{RADK059DGU}

We now recall the higher-derived-bracket construction in the formal neighbourhood of $N$. The formal setting is needed because the brackets depend on all normal derivatives of the ambient Poisson bivector.

Let $\Ic$ be the ideal of smooth functions on $\nu$ vanishing on its zero section. The algebra of \defn{formal functions} along $N$ is
\[
\Of=\varprojlim_{k\ge1}C^\infty\bigl(\operatorname{tot}(\nu)\bigr)/\Ic^k.
\]
We use the same symbol $\Ic$ for the kernel of the restriction $\Of\to C^\infty(N)$. Let $\X^\bullet(\nu)$ be the corresponding completion of the multivector fields on $\nu$. In local coordinates $x^\mu$ on $N$, with a frame $(e_i)$ of $\nu$ and dual fibre coordinates $y^i$, its coefficients are formal power series in $y$ with coefficients smooth in $x$. Functions on $N$ are identified with their pullbacks to $\nu$, and formal vector fields $\X(\nu)$ act as continuous derivations of $\Of$.

We write $\hh=\X^\bullet(\nu)[1]$, a graded Lie algebra for the Schouten bracket. Thus a $k$-vector field has degree $|u|=k-1$ in $\hh$. Our bracket convention is fixed by
\begin{equation}\label{RADK808KXY}
[X,f]=X(f),\quad [X,Y]=X\circ Y-Y\circ X,
\end{equation}
for formal vector fields $X,Y$ and formal functions $f$. So $[X,u]$ is the Lie derivative of $u$ along $X$, and, by the Leibniz rule,
\begin{equation}\label{RADK739THN}
[u,f]=(-1)^{|u|}\,\iota_{df}u ,
\end{equation}
where $\iota_{df}(X\wedge Y)=X(f)\,Y-Y(f)\,X$.

\begin{definition}\label{RADK189AVW}
Let $\mathcal E$ be the Euler vector field of $\nu$. A formal multivector field $u$ has \defn{weight} $w$ if $[\mathcal E,u]=w\,u$. In coordinates, $x^\mu$ and $\dd{x^\mu}$ have weight $0$, $y^i$ has weight $1$, and $\dd{y^i}$ has weight $-1$. We write $u_{(w)}$ for the component of weight $w$ of $u$, and $\X_{\ge k}(\nu)$ for the set of formal vector fields $Z$ with $Z_{(w)}=0$ for $w<k$.
\end{definition}

A $k$-vector field has weights $\ge-k$, and weights add under the Schouten bracket, so $\X_{\ge1}(\nu)$ is a Lie subalgebra and $\ad_Z$ raises the weight filtration for $Z\in\X_{\ge1}(\nu)$. In coordinates,
\begin{equation}\label{RADK254EXY}
Z=Z^\mu\dd{x^\mu}+Z^i\dd{y^i}\in\X_{\ge1}(\nu)\iff Z^\mu\in\Ic\ \text{ and }\ Z^i\in\Ic^2 .
\end{equation}

\begin{definition}\label{RADK223NJS}
A \defn{formal diffeomorphism} $\varphi$ of $\nu$ along $N$ is specified by a continuous unital $\R$-algebra automorphism $\varphi^*$ of $\Of$, with $(\varphi\circ\varphi')^*=\varphi'^*\circ\varphi^*$. It is a \defn{formal change of tubular neighbourhood} if $\varphi^*$ induces the identity on $\Of/\Ic=C^\infty(N)$ and on $\Ic/\Ic^2=C^\infty(\nu^*)$. These form a group $\DiffI(\nu)$.
\end{definition}

In coordinates, the condition reads $\varphi^*(x^\mu)-x^\mu\in\Ic$ and $\varphi^*(y^i)-y^i\in\Ic^2$. If $\gamma$ is a diffeomorphism between neighbourhoods of the zero section which is the identity on $N$ and induces the identity on $\nu$, then $f\mapsto f\circ\gamma$ defines an element of $\DiffI(\nu)$, the \defn{expansion} of $\gamma$. An element $\varphi\in\DiffI(\nu)$ acts on $\X^\bullet(\nu)$ by the automorphism $\varphi^*$ of Gerstenhaber algebras which is $\varphi^*X=\varphi^*\circ X\circ(\varphi^*)^{-1}$ on vector fields. For the expansion of $\gamma$, it is the pullback $\gamma^*=(\gamma^{-1})_*$ of multivector fields, where $\gamma_*$ is the pushforward.

For $Z\in\X_{\ge1}(\nu)$, we write $e^Z$ for the formal diffeomorphism whose pullback on functions is the operator $\exp(Z)$. Accordingly, $\log\varphi$ denotes the formal vector field $\log(\varphi^*)$. The next lemma justifies this notation.

\begin{lemma}\label{RADK768RKF}
\begin{enumerate}
\item For $Z\in\X_{\ge1}(\nu)$, the series $e^Z=\sum_k\frac1{k!}Z^k$ converges on $\Of$ and defines an element of $\DiffI(\nu)$.
\item The map $\exp\colon\X_{\ge1}(\nu)\to\DiffI(\nu)$ is bijective, with inverse $\log$.
\item For $\varphi=e^Z$ and $u\in\X^\bullet(\nu)$, one has
\begin{equation}\label{RADK781EGD}
\varphi^*u=e^{\ad_Z}u .
\end{equation}
\end{enumerate}
\end{lemma}

\begin{proof}
By \eqref{RADK254EXY}, $Z$ maps $\Ic^w$ into $\Ic^{w+1}$, so $e^Z$ converges, is an automorphism with inverse $e^{-Z}$, and lies in $\DiffI(\nu)$. Conversely, for $\varphi\in\DiffI(\nu)$, writing $f\in\Ic^w$ as $\sum_{|\beta|=w}y^\beta f_\beta$ shows that $\varphi^*-\id$ maps $\Ic^w$ into $\Ic^{w+1}$. Thus $\log\varphi^*$ converges. To check the derivation property, work modulo $\Ic^k$ and write $A=\varphi^*-\id$. The finite binomial series $(\id+A)^t$ is polynomial in $t$ and is multiplicative for every nonnegative integer $t$, hence for every real $t$. Differentiating at $t=0$ shows that $\log\varphi^*$ is a continuous derivation of $\Of$, hence a formal vector field determined by its values on the $x^\mu$ and $y^i$. It lies in $\X_{\ge1}(\nu)$ by \eqref{RADK254EXY}. For (3), both sides are automorphisms of $\X^\bullet(\nu)$ agreeing on functions and vector fields, by \eqref{RADK808KXY}.
\end{proof}

By Definition \ref{RADK223NJS}, $(e^Z\circ e^{Z'})^*=e^{Z'}e^Z$ on $\Of$, so
\begin{equation}\label{RADK212LKY}
e^Z\circ e^{Z'}=e^{\BCH(Z',Z)},
\end{equation}
and the Lie algebra of $\DiffI(\nu)$ is $\X_{\ge1}(\nu)$ with the opposite bracket.

A \defn{tubular neighbourhood} of $N$ is an open embedding $\psi$ into $M$ of a neighbourhood of the zero section of $\nu$ which is the identity on $N$ and induces the identity on $\nu$. We write $\sigma_\psi\in\X^2(\nu)$ for the expansion along $N$ of the pullback $\psi^*\sigma=(\psi^{-1})_*\sigma$. For two tubular neighbourhoods, the \defn{transition} $\varphi_{\psi\psi'}\in\DiffI(\nu)$ is the expansion of $\psi^{-1}\circ\psi'$, and
\begin{equation}\label{RADK722CVJ}
\varphi_{\psi\psi'}^*\,\sigma_\psi=\sigma_{\psi'} .
\end{equation}
We write
\[
\DiffI(\psi,\psi')=\big\{\varphi\in\DiffI(\nu):\varphi^*\sigma_\psi=\sigma_{\psi'}\big\}.
\]

We now embed $\mathfrak{a}[1]$ in $\hh$, as in \cite[\S2.6]{Cattaneo2007} and \cite[\S4.1]{cattaneo2008equivalences}. A section $a=a^{i_1\cdots i_k}(x)\,e_{i_1}\wedge\cdots\wedge e_{i_k}$ of $\wedge^k\nu$ is identified with the vertical $k$-vector field $a^{i_1\cdots i_k}(x)\,\dd{y^{i_1}}\wedge\cdots\wedge\dd{y^{i_k}}$. Let $P\colon\hh\to\mathfrak{a}[1]$ be the restriction to $N$ followed by the projection onto $\wedge^\bullet\nu$, so that in coordinates $P$ sets $y=0$ and removes the terms containing some $\dd{x^\mu}$, and let $\pp=\ker P$. Then $(\hh,\mathfrak{a}[1],P)$ is a V-algebra \cite[Definition 2.3]{cattaneo2008equivalences}: $\mathfrak{a}[1]$ is an abelian subalgebra, $\pp$ is a subalgebra, and $\hh=\mathfrak{a}[1]\oplus\pp$.

\begin{definition}
For $u\in\hh$, let $\D_u$ be the coderivation of $S(\mathfrak{a}[1])$ with Taylor coefficients
\[
(\D_u)_0=P(u),\quad (\D_u)_n(a_1,\dots,a_n) =P\bigl[\cdots[[u,a_1],a_2],\dots,a_n\bigr], \quad n\ge1.
\]
\end{definition}

These coefficients are graded symmetric, since $\mathfrak{a}[1]$ is abelian, and $\D_u$ has degree $|u|$. For the frame element $e_i$ and the coordinate $x^\mu$, the brackets with $e_i$ and $x^\mu$ are, by \eqref{RADK739THN},
\begin{equation}\label{RADK133NVD}
[u,e_i]=-\frac{\partial u}{\partial y^i},\quad [u,x^\mu]=(-1)^{|u|}\,\iota_{dx^\mu}u ,
\end{equation}
where $\frac{\partial u}{\partial y^i}$ is the derivative of the coefficients of $u$. We use two results of Voronov. For every odd $u\in\hh$, $[\D_u,\D_u]=\D_{[u,u]}$ \cite[Theorem 1]{voronov2005higher}. For $u,v\in\pp$,
\begin{equation}\label{RADK256NMA}
\D_{[u,v]}=[\D_u,\D_v]
\end{equation}
\cite[\S5, Theorem 3]{voronov2005arbitrary}, since $\D_u$ is the coderivation attached by Voronov to the derivation $\ad_u$ of $\hh$ preserving $\pp$.

We write $\Coder\big(S(\mathfrak{a}[1])\big)$ for the space of coderivations of $S(\mathfrak{a}[1])$, and $\MultiDer(\mathfrak{a})$ for the space of coderivations of degree $0$ whose coefficients of positive arity are multiderivations. Since $\mathfrak{a}[1]$ is abelian and $P$ is multiplicative, the Leibniz rule of the Schouten bracket shows that $(\D_u)_n$ is a multiderivation of $\mathfrak{a}$ of degree $|u|$, for every $u\in\hh$ and $n\ge1$.
The Fourier transform of Cattaneo and Felder \cite[Theorem 4.1 and Proposition 4.4]{Cattaneo2007} identifies formal multivectors with their higher derived brackets. In the form needed here, \eqref{RADK133NVD} shows that the coefficients of $\D_u$ determine all normal derivatives of the coefficients of $u$. Consequently, $\D\colon\hh\to\Coder\bigl(S(\mathfrak{a}[1])\bigr)$ is injective and restricts to a bijection
\begin{equation}\label{RADK572NAF}
\D\colon\X(\nu)\longrightarrow\MultiDer(\mathfrak{a}).
\end{equation}

With our conventions, \eqref{RADK133NVD} and $[Z,x^\mu]=Z^\mu$ give, for $Z=Z^\mu\dd{x^\mu}+Z^i\dd{y^i}$,
\begin{align}
(\D_Z)_n(e_{i_1},\dots,e_{i_n})&=(-1)^n\,\frac{\partial^nZ^i}{\partial y^{i_1}\cdots\partial y^{i_n}}(x,0)\,e_i ,\label{RADK914PAN}\\
(\D_Z)_n(e_{i_1},\dots,e_{i_{n-1}},x^\mu)&=(-1)^{n-1}\,\frac{\partial^{n-1}Z^\mu}{\partial y^{i_1}\cdots\partial y^{i_{n-1}}}(x,0).\label{RADK862PHS}
\end{align}
Since a multiderivation is determined by its values on the $x^\mu$ and the $e_i$, which are symmetric in the frame arguments, these formulas also show that every element of $\MultiDer(\mathfrak{a})$ is of the form $\D_Z$, the local vector fields $Z$ gluing by injectivity of $\D$.

The weights of $u$ control the arities of $\D_u$.

\begin{lemma}\label{RADK332ENG}
Let $u\in\hh$ be homogeneous of degree $d$. Then $(\D_u)_n=(\D_{u_{(n-1-d)}})_n$ for every $n\ge0$. In particular, the following hold.
\begin{enumerate}
\item For a vector field $Z$, $(\D_Z)_n$ depends only on $Z_{(n-1)}$, and $\X_{\ge1}(\nu)=\{Z\in\X(\nu):(\D_Z)_0=(\D_Z)_1=0\}$. Hence $\D$ is an injective morphism of Lie algebras from $\X_{\ge1}(\nu)$ to $\CoderII(\mathfrak{a})$, with image $\MultiDer(\mathfrak{a})\cap\CoderII(\mathfrak{a})$.
\item For a bivector $\pi$, $(\D_\pi)_n$ depends only on $\pi_{(n-2)}$, and $P(\pi)=0$ if and only if $\pi_{(-2)}=0$. In particular, $\sigma_\psi$ has weight $\ge-1$, since $P(\sigma_\psi)=\frac12\sigma_\psi(dy^j,dy^k)|_N\,e_j\wedge e_k=0$ for $N$ coisotropic \cite[Lemma 5.17]{aitmohamedmayrand2026}.
\end{enumerate}
\end{lemma}

\begin{proof}
A section of $\wedge^k\nu$ has degree $k-1$ and weight $-k$, and $P$ only retains terms of degree $m-1$ and weight $-m$. Counting degrees and weights, $P[\cdots[u_{(w)},a_1],\dots,a_n]=0$ unless $w=n-1-d$. For (1), $d=0$, so $(\D_Z)_0$ and $(\D_Z)_1$ depend only on $Z_{(-1)}$ and $Z_{(0)}$, and vanish exactly when $Z\in\X_{\ge1}(\nu)$ by \eqref{RADK572NAF}. Since $\X_{\ge1}(\nu)\subset\pp$ by \eqref{RADK254EXY}, the rest follows from \eqref{RADK256NMA} and \eqref{RADK572NAF}. For (2), $d=1$ and $P(\pi)=\pi_{(-2)}$.
\end{proof}

In particular, $\D$ is continuous from the weight filtration to the arity filtration, so it commutes with convergent Lie series in $\X_{\ge1}(\nu)$.

For a formal Poisson bivector $\pi$ with $P(\pi)=0$, we set $Q_\pi=\D_\pi$, and $Q_\psi=Q_{\sigma_\psi}$. By Voronov's identity, $Q_\pi\circ Q_\pi=\frac12\D_{[\pi,\pi]}=0$, so $Q_\pi$ is a flat $\Linf$-structure on $\mathfrak{a}[1]$, and a flat $\Pinf$-structure on $\mathfrak{a}$, since the $(\D_u)_n$ are multiderivations. On the unshifted algebra, we set
\begin{equation}\label{RADK840EAG}
\lambda^\psi_n(a_1,\dots,a_n)=(-1)^{\,1+\sum_{i=1}^n(n-i)|a_i|}\,(Q_\psi)_n(a_1,\dots,a_n).
\end{equation}
These maps have degree $2-n$, are graded antisymmetric, and are derivations in the last argument. We call $(\mathfrak{a},\lambda^\psi)$ the \defn{homotopy Poisson algebra} of $(N,\psi)$, and write $\lambda_n$ for $\lambda^\psi_n$ when $\psi$ is fixed.

The brackets are explicit in coordinates. Let $\{f,g\}=\sigma_\psi(df,dg)$. For $f,g\in C^\infty(N)$, \eqref{RADK133NVD} and \eqref{RADK739THN} give
\begin{equation}\label{RADK377FMB}
\begin{aligned}
\lambda_n(e_{i_1},\dots,e_{i_{n-2}},f,g)&=(-1)^{n-1}\,\frac{\partial^{n-2}\{f,g\}}{\partial y^{i_1}\cdots\partial y^{i_{n-2}}}\Big|_{y=0},\quad n\ge2,\\
\lambda_n(e_{i_1},\dots,e_{i_{n-1}},f)&=(-1)^{n-1}\,\frac{\partial^{n-1}\{f,y^j\}}{\partial y^{i_1}\cdots\partial y^{i_{n-1}}}\Big|_{y=0}\,e_j,\quad n\ge1,\\
\lambda_n(e_{i_1},\dots,e_{i_n})&=\frac{(-1)^{n-1}}{2}\,\frac{\partial^{n}\{y^j,y^k\}}{\partial y^{i_1}\cdots\partial y^{i_n}}\Big|_{y=0}\,e_j\wedge e_k,\quad n\ge0 .
\end{aligned}
\end{equation}
These are the formulas of \cite[Proposition 2.1]{Cattaneo2007} up to the sign $(-1)^{m+1}$, where $m$ is the number of arguments in $C^\infty(N)$.

\begin{remark}\label{RADK135FUR}
The differential $\lambda_1^\psi$ does not depend on $\psi$, and $-\lambda_1^\psi$ is the Chevalley--Eilenberg differential of the conormal Lie algebroid $TN^\circ$ with the conventions of \S\ref{RADK219CNV} \cite[Remark 4.2]{cattaneo2008equivalences}. Indeed, both are derivations, and for $f\in C^\infty(N)$ and $\xi\in TN^\circ$ the second formula of \eqref{RADK377FMB} gives $\langle\lambda_1^\psi(f),\xi\rangle=\sigma(df,\xi)=-\sigma(\xi,df)$, and similarly on the frame $e_i$. So the algebra automorphism $a\mapsto(-1)^{\deg a}a$ of $\mathfrak{a}$ identifies $(\mathfrak{a},\lambda_1^\psi)$ with the Chevalley--Eilenberg complex of $TN^\circ$.
\end{remark}

\subsection{Changes of tubular neighbourhood}\label{RADK601TBN}

We now integrate $\D$ to a map $\Phi\colon\DiffI(\nu)\to\AutI(\mathfrak{a})$ and prove Theorem \ref{RADK919BTK}. By Lemma \ref{RADK332ENG}(1), \eqref{RADK212LKY}, \eqref{RADK042DFG} and Lemma \ref{RADK768RKF}, we may set
\begin{equation}\label{RADK999XZC}
\Phi(\varphi)=\exp\big(\D_{\log\varphi}\big)\in\AutI(\mathfrak{a}),\quad\text{that is,}\quad \Phi(e^Z)=e^{\D_Z}.
\end{equation}

\begin{proposition}\label{RADK286QVV}
Let $\varphi,\varphi'\in\DiffI(\nu)$, and let $\pi$ be a formal Poisson bivector with $P(\pi)=0$.
\begin{enumerate}
\item $\Phi(\varphi)\circ Q_\pi\circ\Phi(\varphi)^{-1}=Q_{\varphi^*\pi}$.
\item $\Phi(\varphi\circ\varphi')=\Phi(\varphi')\circ\Phi(\varphi)$.
\item $\Phi$ is injective.
\item If $\varphi\in\DiffI(\psi,\psi')$, then $\Phi(\varphi)$ is an $\Linf$-isomorphism $(\mathfrak{a},\lambda^\psi)\to(\mathfrak{a},\lambda^{\psi'})$.
\end{enumerate}
\end{proposition}

\begin{proof}
Let $Z=\log\varphi$ and $Z'=\log\varphi'$. By Lemma \ref{RADK332ENG}(2), $\pi$ has weight $\ge-1$, so $\ad_Z^k\pi$ has weight $\ge k-1$ and $\D_{\ad_Z^k\pi}$ vanishes in arities $\le k$. Since $Z,\pi\in\pp$, \eqref{RADK496NHH}, \eqref{RADK256NMA} and \eqref{RADK781EGD} give
\begin{align*}
\Phi(\varphi)\circ\D_\pi\circ\Phi(\varphi)^{-1}
&=\sum_{k\ge0}\frac1{k!}\,\ad_{\D_Z}^k(\D_\pi)\\
&=\sum_{k\ge0}\frac1{k!}\,\D_{\ad_Z^k\pi}\\
&=\D_{e^{\ad_Z}\pi}
\\&=\D_{\varphi^*\pi},
\end{align*}
and $P(\varphi^*\pi)=0$ since $\varphi^*\pi$ has weight $\ge-1$. For (2), \eqref{RADK212LKY}, \eqref{RADK496NHH} and \eqref{RADK256NMA} give
\[
\Phi(\varphi')\circ\Phi(\varphi) =e^{\BCH(\D_{Z'},\D_Z)} =e^{\D_{\BCH(Z',Z)}} =\Phi(\varphi\circ\varphi').
\] Part (3) follows from \eqref{RADK042DFG}, \eqref{RADK572NAF} and Lemma \ref{RADK768RKF}, and (4) from (1) with $\pi=\sigma_\psi$.
\end{proof}

\begin{corollary}\label{RADK250MSL}
For tubular neighbourhoods $\psi,\psi',\psi''$, the map $\Phi(\varphi_{\psi\psi'})$ is an $\Linf$-isomorphism $(\mathfrak{a},\lambda^\psi)\to(\mathfrak{a},\lambda^{\psi'})$ with identity linear part. These isomorphisms satisfy
\[
\Phi(\varphi_{\psi'\psi''})\circ\Phi(\varphi_{\psi\psi'})
=\Phi(\varphi_{\psi\psi''}).
\]
\end{corollary}

\begin{proof}
This follows from \eqref{RADK722CVJ}, Proposition \ref{RADK286QVV}(2) and (4), and $\varphi_{\psi\psi''}=\varphi_{\psi\psi'}\circ\varphi_{\psi'\psi''}$.
\end{proof}

We next compare $\Phi$ with the $\Linf$-isomorphisms attached by Cattaneo and Sch\"atz to isotopies of tubular neighbourhoods \cite{cattaneo2008equivalences}. A family $(\varphi_\tau)_{\tau\in[0,1]}$ in $\DiffI(\nu)$ is smooth if, in coordinates, the coefficients of $\varphi_\tau^*(x^\mu)$ and $\varphi_\tau^*(y^i)$ depend smoothly on $(\tau,x)$, and similarly for families of multidifferential operators. We use the formula for the derivative of the exponential map \cite{rossmann2002}
\begin{equation}\label{RADK122EXZ}
\Big(\frac{d}{d\tau}e^{B_\tau}\Big)\circ e^{-B_\tau}=\sum_{k\ge0}\frac{1}{(k+1)!}\,\ad_{B_\tau}^k\Big(\frac{d}{d\tau}B_\tau\Big),
\end{equation}
for a smooth family of operators $B_\tau$ which either lowers an exhaustive increasing filtration or raises a complete decreasing filtration by at least one step. Indeed, all series are then finite on each filtration step, respectively modulo each filtration step, and the usual proof applies.

\begin{lemma}\label{RADK831LKE}
Let $(\varphi_\tau)$ be a smooth family in $\DiffI(\nu)$ with $\varphi_0=\id$, and $Z_\tau=\log\varphi_\tau$. Then $Y_\tau=\big(\frac{d}{d\tau}\varphi_\tau^*\big)\circ(\varphi_\tau^*)^{-1}$ is the vector field
\begin{align*}
Y_\tau
&=\sum_{k\ge0}\frac{1}{(k+1)!}\,
\ad_{Z_\tau}^k\Bigl(\frac{dZ_\tau}{d\tau}\Bigr)
\in\X_{\ge1}(\nu),\\
\frac{d}{d\tau}\varphi_\tau^*
&=\ad_{Y_\tau}\circ\varphi_\tau^*
\quad\text{on }\hh.
\end{align*}
\end{lemma}

\begin{proof}
Each weight component of $\log\varphi_\tau^*$ is a polynomial in finitely many jets of $\varphi_\tau^*(x^\mu)$ and $\varphi_\tau^*(y^i)$, so $Z_\tau$ is smooth in $\tau$. On $\Of$, $\varphi_\tau^*=e^{Z_\tau}$ and $Z_\tau$ raises the weight filtration, so \eqref{RADK122EXZ} gives $Y_\tau\in\X_{\ge1}(\nu)$, a closed Lie subalgebra. On $\hh$, $\varphi_\tau^*=e^{\ad_{Z_\tau}}$ by \eqref{RADK781EGD}, and \eqref{RADK122EXZ} with $B_\tau=\ad_{Z_\tau}$ gives $\ad_{Y_\tau}$, since $\ad$ is a morphism of Lie algebras.
\end{proof}

Let $\psi_\tau=\psi\circ\gamma_\tau$ be an isotopy of tubular neighbourhoods, with $\gamma_0=\id$, and let $\varphi_\tau$ be the expansion of $\gamma_\tau$. By \eqref{RADK722CVJ} and Lemma \ref{RADK831LKE}, the automorphisms $\varphi_\tau^*$ carry $\sigma_\psi$ to $\sigma_{\psi_\tau}$, with generator $\ad_{Y_\tau}$ in the sense of \cite[eq.(7)]{cattaneo2008equivalences}. Since $Y_\tau\in\pp$, the derivation $\ad_{Y_\tau}$ preserves $\pp$, which is condition (9) of \cite{cattaneo2008equivalences}. The Cauchy problem (10) of \cite{cattaneo2008equivalences} has only the zero solution, since $P[Y_\tau,a]=(\D_{Y_\tau})_1(a)=0$. So the $\Linf$-isomorphism of \cite[Theorems 3.2 and 4.3]{cattaneo2008equivalences} attached to the isotopy is $R_1$, where $R_\tau$ solves the Cauchy problem \cite[eq.(14)]{cattaneo2008equivalences}
\begin{equation}\label{RADK453RFE}
\frac{d}{d\tau}R_\tau=\D_{Y_\tau}\circ R_\tau,\quad R_0=\id .
\end{equation}

\begin{proposition}\label{RADK296FDU}
Let $(\varphi_\tau)$ be a smooth family in $\DiffI(\nu)$ with $\varphi_0=\id$. The solution of \eqref{RADK453RFE} is $R_\tau=\Phi(\varphi_\tau)$. In particular, the Cattaneo--Sch\"atz isomorphism of an isotopy is $\Phi(\varphi_1)$, and it depends only on $\varphi_1$.
\end{proposition}

\begin{proof}
By Lemma \ref{RADK831LKE}, \eqref{RADK256NMA} and Lemma \ref{RADK332ENG},
\[ \D_{Y_\tau}=\sum_{k\ge0}\frac1{(k+1)!}\,\ad_{\D_{Z_\tau}}^k\Big(\frac{d}{d\tau}\D_{Z_\tau}\Big).
\]
By \eqref{RADK122EXZ} on $S(\mathfrak{a}[1])$ with $B_\tau=\D_{Z_\tau}$, which lowers the length filtration, this is $\big(\frac{d}{d\tau}e^{\D_{Z_\tau}}\big)\circ e^{-\D_{Z_\tau}}$. So $\Phi(\varphi_\tau)$ solves \eqref{RADK453RFE}, and the solution is unique by \cite[Proposition 3.4]{cattaneo2008equivalences}.
\end{proof}

An $\Linf$-isomorphism $F\colon(\mathfrak{a},\lambda^\psi)\to(\mathfrak{a},\lambda^{\psi'})$ is \defn{induced by a change of tubular neighbourhood} if $F=\Phi(\varphi)$ for some $\varphi\in\DiffI(\psi,\psi')$.

\begin{proposition}\label{RADK054KVK}
Let $F\colon(\mathfrak{a},\lambda^\psi)\to(\mathfrak{a},\lambda^{\psi'})$ be an $\Linf$-isomorphism. The following are equivalent.
\begin{enumerate}
\item $F$ is induced by a change of tubular neighbourhood.
\item $F=e^{\D_Z}$ for some $Z\in\X_{\ge1}(\nu)$.
\item $F_1=\id$, and every Taylor coefficient of positive arity of $\log F$ is a multiderivation of $\mathfrak{a}$.
\end{enumerate}
In that case, the element $\varphi\in\DiffI(\psi,\psi')$ with $F=\Phi(\varphi)$ is unique.
\end{proposition}

\begin{proof}
By definition, (1) implies (2), and (2) implies (3) by Lemma \ref{RADK332ENG}(1) and \eqref{RADK042DFG}. Assume (3). Then $F\in\AutI(\mathfrak{a})$ and $\log F\in\MultiDer(\mathfrak{a})\cap\CoderII(\mathfrak{a})$, so $\log F=\D_Z$ for a unique $Z\in\X_{\ge1}(\nu)$ by Lemma \ref{RADK332ENG}(1). Set $\varphi=e^Z$, so that $F=\Phi(\varphi)$. Since $F$ is an $\Linf$-isomorphism, Proposition \ref{RADK286QVV}(1) gives
\[
\D_{\varphi^*\sigma_\psi} =F\circ Q_\psi\circ F^{-1} =Q_{\psi'} =\D_{\sigma_{\psi'}}.
\]
Injectivity in \eqref{RADK572NAF} therefore gives $\varphi^*\sigma_\psi=\sigma_{\psi'}$, which proves (1). Uniqueness is Proposition \ref{RADK286QVV}(3).
\end{proof}

\begin{proof}[\textnormal{\bfseries Proof of Theorem \ref{RADK919BTK}}]
This is the equivalence of (1) and (3) in Proposition \ref{RADK054KVK}.
\end{proof}

\begin{remark}\label{RADK827PMG}
The element $\varphi$ of Proposition \ref{RADK054KVK} need not be the transition $\varphi_{\psi\psi'}$. Indeed, $\DiffI(\psi,\psi')=\varphi_{\psi\psi'}\circ\DiffI(\psi',\psi')$, and by \eqref{RADK781EGD} the stabilizer $\DiffI(\psi',\psi')$ consists of the $e^X$ with $X\in\X_{\ge1}(\nu)$ and $[X,\sigma_{\psi'}]=0$, since the component of lowest weight of $e^{\ad_X}\sigma_{\psi'}-\sigma_{\psi'}$ is that of $[X,\sigma_{\psi'}]$. For instance, for a formal function $g$ of weight $\ge2$, the Hamiltonian vector field $X_g=\sigma_{\psi'}^\sharp(dg)$ lies in $\X_{\ge1}(\nu)$ and commutes with $\sigma_{\psi'}$.
\end{remark}

\begin{remark}\label{RADK881TQT}
An induced isomorphism need not have multiderivations as Taylor coefficients. Let $M=\R^4$ with $\sigma=\dd{x^1}\wedge\dd{y^1}+\dd{x^2}\wedge\dd{y^2}$, let $N=\{y=0\}$, with frame $e_i=\dd{y^i}$, and let $\psi$ be the identity. Let $Z=\tfrac12\big((y^1)^2\dd{y^1}+(y^2)^2\dd{y^2}\big)$, whose time-one flow is $\gamma(x,y)=\big(x,\frac{y^1}{1-y^1/2},\frac{y^2}{1-y^2/2}\big)$, and $\psi'=\psi\circ\gamma$. The expansion of $\gamma$ is $e^Z=\varphi_{\psi\psi'}$, so $F=\Phi(e^Z)$ is induced, by \eqref{RADK722CVJ}. Let $\delta=\D_Z$, which is concentrated in arity $2$ by Lemma \ref{RADK332ENG}(1). The Lie derivative of $\dd{y^1}\wedge\dd{y^2}$ along $Z$ is $-(y^1+y^2)\,\dd{y^1}\wedge\dd{y^2}$, so \eqref{RADK914PAN} and \eqref{RADK133NVD} give
\[
\delta_2(e_1,e_1)=e_1,\quad \delta_2(e_2,e_2)=e_2,\quad \delta_2(e_1,e_2)=0,\quad \delta_2(e_1\wedge e_2,e_j)=e_1\wedge e_2.
\]
By \eqref{RADK675CTY} and \eqref{RADK504RBS} with $(\log F)_3=0$, we get $F_3=\tfrac12\sum\varepsilon\,\delta_2(\delta_2(a_i,a_j),a_k)$, with all Koszul signs equal to $+1$ for the arguments below. Hence $F_3(e_1\wedge e_2,e_1,e_2)=e_1\wedge e_2$, while $F_3(e_1,e_1,e_2)=\tfrac12\delta_2(e_1,e_2)=0$ and $F_3(e_2,e_1,e_2)=0$. If $F_3$ were a multiderivation, \eqref{RADK049HQH} would give $F_3(e_1\wedge e_2,e_1,e_2)=0$.
\end{remark}

\section{The adjoint symplectic groupoid of a Radko surface}\label{RADK261MGJ}

In this section, we give an explicit coordinate description of the adjoint symplectic groupoid of a Radko surface $(\Sigma,\sigma)$, using the integration theory of Gualtieri and Li \cite{gualtieri2014symplectic}. This fixes the ambient symplectic groupoid $\G_\Sigma$ in which we will construct the Lagrangian integrations of admissible curves in \S\ref{RADK793NZP} and compute their homotopy fibre products in \S\ref{RADK389HFP}.

\subsection{Conventions and Radko surfaces}\label{RADK219CNV}

We first fix our conventions and recall the normal form of a Radko surface near its degeneracy locus.

For a bivector $\pi$ on a manifold $M$, we use the same symbol for the bundle map
\[
\pi\colon T^*M\longrightarrow TM,\quad
\pi(\xi)=\iota_\xi\pi=\pi(\xi,\cdot).
\]
If $\pi$ is Poisson, our conventions for the Poisson bracket and Hamiltonian vector fields are $\{f,g\}=\pi(df,dg)$, $X_f=\pi(df)$.
The cotangent Lie algebroid $T^*_\pi M$ has anchor $\pi$ and bracket $[df,dg]=d\{f,g\}$. A $2$-form $\omega$ and a bivector $\pi$ are \defn{inverse} if
\[
\iota_{\pi(\xi)}\omega=\xi\quad\text{for every }\xi\in T^*M.
\]
Equivalently, their graphs agree:
\[
\Gamma_\omega=\big\{(v,\iota_v\omega):v\in TM\big\}
=\big\{(\pi(\xi),\xi):\xi\in T^*M\big\}.
\]
We then write $\pi=\omega^{-1}$. In particular, $\iota_{X_f}\omega=df$.

Following \cite{gualtieri2014symplectic}, we allow the arrow manifold of a Lie groupoid $\G\tto M$ to be non-Hausdorff. We compose $g$ and $h$ in the order $gh$ when $\sss(g)=\ttt(h)$, and write $\G_x$ for the isotropy group at $x\in M$. For an open set $U\subset M$, the restriction is $\G|_U=\sss^{-1}(U)\cap\ttt^{-1}(U)$.
The source-connected part $\G^c$ is the union, over $x\in M$, of the connected components of $\sss^{-1}(x)$ containing $1_x$. The Lie algebroid $\Lie(\G)$ is $\ker d\sss$ along the units, with anchor $d\ttt$ and bracket induced by right-invariant vector fields.

A \defn{symplectic groupoid} is a Lie groupoid equipped with a symplectic form $\Omega$ which is \defn{multiplicative}: $m^*\Omega=\pr_1^*\Omega+\pr_2^*\Omega$,
where $m$ is multiplication and $\pr_1,\pr_2$ are the projections from the space of composable pairs. It \defn{integrates} $(M,\pi)$ if the target map $\ttt\colon(\G,\Omega^{-1})\to(M,\pi)$ is Poisson. In this case,
\[
\Lie(\G)\longrightarrow T^*_\pi M,\quad
v\longmapsto\iota_v\Omega|_{TM}
\]
is an isomorphism of Lie algebroids, and the induced Dirac structure
\[
\big\{\big(d\ttt(v),\iota_v\Omega|_{TM}\big):v\in\Lie(\G)\big\}
\]
is the graph of $\pi$. These are the conventions of \cite[\S\S2.1 and 3.1, Lemma 3.1, Example 3.8 and Remark 4.9]{Maxence2025}. Gualtieri and Li use the opposite convention, in which $\sss$ is Poisson \cite[\S1]{gualtieri2014symplectic}.

A \defn{Radko surface} is a compact oriented surface $\Sigma$ with a Poisson bivector $\sigma$ which vanishes transversely along a nonempty curve $D$ and is nondegenerate on $\Sigma\setminus D$ \cite{radko2002classification}. Writing $\sigma=f\sigma_0$ with $\sigma_0$ a nowhere vanishing bivector, this means that $df\neq0$ along $D=f^{-1}(0)$. We do not assume $\Sigma$ connected.

\begin{figure}[ht]
\centering
\includegraphics[width=0.6\linewidth]{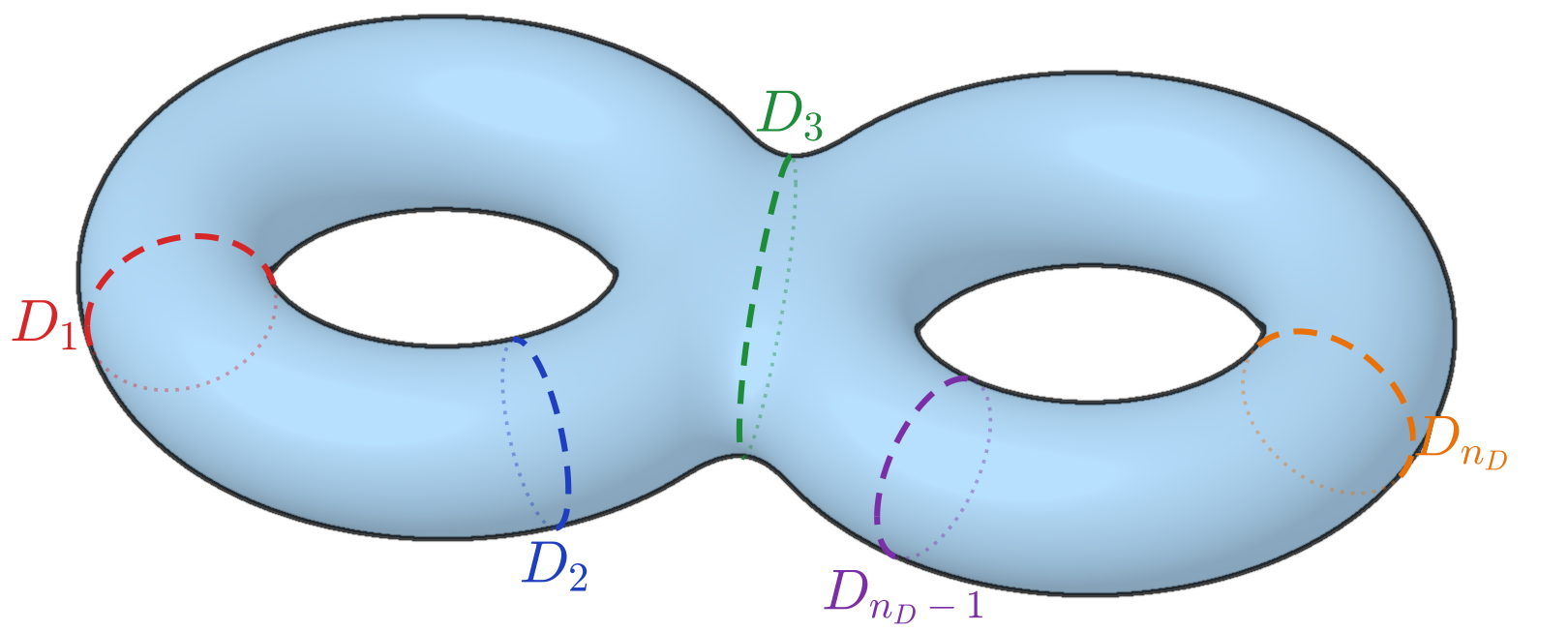}
\caption{A Radko surface and its degeneracy locus.}\label{RADK701FGA}
\end{figure}

The degeneracy locus $D=f^{-1}(0)$ is a closed one-dimensional submanifold of the compact surface $\Sigma$, hence a finite union of circles $D_1,\dots,D_{n_D}$ (Figure \ref{RADK701FGA}). Its normal bundles $\nu_{D_j}$ are trivial, since both $\Sigma$ and $D_j$ are orientable. Write $V_i$ for the connected components of $\Sigma\setminus D$ and $\omega_i$ for the inverse of $\sigma$ on $V_i$. The symplectic leaves are the $V_i$ and the individual points of $D$. Thus the components $D_j$ and all leaves contained in them are compact, so $(\Sigma,\sigma)$ is a \defn{proper} log symplectic surface in the sense of \cite[Definition 1.10]{gualtieri2014symplectic}.

A log symplectic manifold is linearizable along its degeneracy locus \cite[Theorem 1.12]{gualtieri2014symplectic}, and the linear models of a proper one are classified in \cite[Proposition 1.11]{gualtieri2014symplectic}. In dimension two, the leaves in $D_j$ are points and $\nu_{D_j}$ is trivial. The local model is therefore determined by one number $\ell_j>0$, the \defn{modular period} of $D_j$ (see \eqref{RADK432CKK} below).

More explicitly, Radko's normal form \cite[\S2.4, eq. (4)]{radko2002classification}, see also \cite[Lemma 2]{radko2002classification}, gives coordinates
\[
(\theta,r)\in(\R/2\pi\Z)\times(-c_j,c_j),\quad
D_j=\{r=0\},\quad
\sigma=\frac{2\pi}{\ell_j}\,r\,\dd r\wedge\dd\theta.
\]
Set $t=\ell_j\theta/(2\pi)$ and $S^1_{\ell_j}=\R/\ell_j\Z$. We identify a neighbourhood $U_j$ of $D_j$ with the Poisson cylinder
\[
\big(S^1_{\ell_j}\times(-c_j,c_j),\ r\,\dd r\wedge\dd t\big),
\]
with $D_j=S^1_{\ell_j}\times\{0\}$. It will be convenient to extend this local model to the whole cylinder

\begin{equation}\label{RADK871UDA}
\widehat U_j=S^1_{\ell_j}\times\R,\quad \sigma=r\,\dd r\wedge\dd t,
\end{equation}
so that $U_j=S^1_{\ell_j}\times(-c_j,c_j)$ is an open Poisson submanifold of $\widehat U_j$. Since the $D_j$ are pairwise disjoint compact sets, we choose the $c_j$ so that the $U_j$ are pairwise disjoint. Then $U_j\cap D=D_j=\{r=0\}$.

On $\{r\neq0\}$, the inverse symplectic form is $\omega=\frac1r\,dt\wedge dr$.
Indeed,
\[
\sigma(dr)=r\,\dd t,\quad \iota_{\sigma(dr)}\omega=dr,\quad \sigma(dt)=-r\,\dd r,\quad \iota_{\sigma(dt)}\omega=dt.
\]
The normal-form identification is Poisson, so $\omega_i=\omega$ on $U_j\cap V_i$.

Let $\mathrm{vol}$ be a volume form on a neighbourhood of $D_j$. With the convention $X_f=\sigma(df)$, so that $X_f(g)=\sigma(df,dg)$, the \defn{modular vector field} $X_{\mathrm{vol}}$ of $\sigma$ \cite{weinstein1997modular} is defined by $X_{\mathrm{vol}}(f)=\divv(X_f)$,
where divergence is defined by $d(\iota_X\mathrm{vol})=\divv(X)\,\mathrm{vol}$.

The restriction of $X_{\mathrm{vol}}$ to $D_j$ does not depend on the volume form \cite[\S2.4]{radko2002classification}, since changing $\mathrm{vol}$ changes $X_{\mathrm{vol}}$ by a Hamiltonian vector field, which vanishes along $D_j$. On $U_j$, take $\mathrm{vol}=dr\wedge dt$. Then
$X_f=r\,\frac{\partial f}{\partial r}\,\dd t-r\,\frac{\partial f}{\partial t}\,\dd r$ and  
$\divv(X_f)=-\frac{\partial f}{\partial t},$ so that
\begin{equation}\label{RADK432CKK}
X_{\mathrm{vol}}\big|_{D_j}=-\dd t .
\end{equation}

In particular, the modular period $\ell_j$ is the period of the modular vector field along $D_j$.

We next show that each leaf $V_i$ meets each cylinder $U_j$ in at most one of the two components of $U_j\setminus D_j$. Write $\sigma=f\sigma_0$ as above. The function $f$ has constant sign on each $V_i$. On $U_j$, $\sigma_0=u\,\dd r\wedge\dd t$ with $u$ nowhere vanishing, so $f=r/u$ has opposite signs on the components $\{r>0\}$ and $\{r<0\}$ of $U_j\setminus D_j$. Each of them is connected and contained in $\Sigma\setminus D$, hence in a single leaf, and the two leaves are distinct. Therefore
\begin{equation}\label{RADK247VWC}
V_i\cap U_j\in\big\{\emptyset,\ S^1_{\ell_j}\times(0,c_j),\ S^1_{\ell_j}\times(-c_j,0)\big\}.
\end{equation}
If $V_i$ meets no $U_j$, then its closure misses $D$, since $\bigcup_jU_j$ is a neighbourhood of $D$, so $V_i$ is open and closed in $\Sigma$. Hence each $V_i$ is a component of the compact surface $\Sigma$ or contains one of the $2n_D$ components of the sets $U_j\setminus D_j$, and there are finitely many leaves $V_i$.

\subsection{Adjoint integrations}\label{RADK237YVT}

We recall the description of source-connected integrations used to identify $\G_\Sigma$ \cite{gualtieri2014symplectic}. Let $E$ be an integrable Lie algebroid over $M$. Its source-connected integrations, with their identifications of Lie algebroids with $E$, form a category $\mathbf{Gpd}(E)$. Morphisms induce the identity on $E$. There is at most one morphism between two objects, and every such morphism is a surjective local diffeomorphism \cite[Definition 1.19 and the discussion following it]{gualtieri2014symplectic}. The initial object is the source-simply-connected integration $\G^{\mathrm{ssc}}$. The terminal object, when it exists, is the \defn{adjoint integration} $\G^{\mathrm{adj}}$.

Let $\Lambda(\G^{\mathrm{ssc}})$ be the poset, ordered by inclusion, of wide Lie subgroupoids $\NN\subset\G^{\mathrm{ssc}}$ which are \defn{discrete}, meaning that $\sss|_\NN$ is a local diffeomorphism, \defn{totally disconnected}, meaning that $\sss=\ttt$ on $\NN$, and \defn{normal}. By \cite[Theorems 1.14 and 1.20]{gualtieri2014symplectic}, the assignment
\[
\Lambda(\G^{\mathrm{ssc}})\longrightarrow\mathbf{Gpd}(E),\quad
\NN\longmapsto\G^{\mathrm{ssc}}/\NN
\]
is an equivalence. The quotient map is a surjective local diffeomorphism, and the quotient is Hausdorff if $\G^{\mathrm{ssc}}$ is Hausdorff and $\NN$ is closed. In particular, a greatest element $\NN\in\Lambda(\G^{\mathrm{ssc}})$ determines the adjoint integration: $\G^{\mathrm{adj}}=\G^{\mathrm{ssc}}/\NN$.

To pass from local to global integrations, an open cover $\mathcal O$ of $M$ is an \defn{orbit cover} if every orbit of $E$ is contained in one of its members. For $O,O'\in\mathcal O$, let $\widetilde\G_O$ and $\widetilde\G_{OO'}$ be the source-simply-connected integrations of $E|_O$ and $E|_{O\cap O'}$. The restriction map
\[
\res_{O,O'}\colon\Lambda(\widetilde\G_O)\longrightarrow\Lambda(\widetilde\G_{OO'})
\]
sends $\NN$ to the preimage of $\NN\cap(\widetilde\G_O|_{O\cap O'})^c$ under the canonical morphism
\[
\widetilde\G_{OO'}\longrightarrow(\widetilde\G_O|_{O\cap O'})^c.
\]
A family $(\NN_O)_{O\in\mathcal O}$ is \defn{compatible} if
\[
\res_{O,O'}(\NN_O)=\res_{O',O}(\NN_{O'})\quad\text{for all }O,O'\in\mathcal O.
\]
For a locally finite orbit cover, compatible families glue to source-connected integrations $\G$ satisfying $(\G|_O)^c=\widetilde\G_O/\NN_O$.
This gives an equivalence between compatible families, ordered by inclusion, and $\mathbf{Gpd}(E)$ \cite[Theorem 3.4 and Corollary 3.7]{gualtieri2014symplectic}. We will use the following \emph{gluing criterion}: if each $\Lambda(\widetilde\G_O)$ has a greatest element $\NN_O$ and these elements form a compatible family, then the glued groupoid is the adjoint integration of $E$. Indeed, every compatible family is contained in $(\NN_O)$.

\subsection{The action groupoid of \texorpdfstring{$\Aff$}{Aff}}\label{RADK265ACT}

We first describe the source-simply-connected symplectic integration of the cylinder $(\widehat U_j,\sigma)$ in \eqref{RADK871UDA}, and then take the quotient giving its adjoint integration. This is the surface case of the action-groupoid model in \cite[\S A.2]{gualtieri2014symplectic}.

The cotangent Lie algebroid $T^*_\sigma\widehat U_j$ has the global frame $(dr,dt)$, with
\[
\sigma(dr)=r\,\dd t,\quad \sigma(dt)=-r\,\dd r,\quad [dr,dt]=d\big(\sigma(dr,dt)\big)=dr.
\]
So its anchor is onto where $r\neq0$ and vanishes along $D_j=\{r=0\}$, and the isotropy Lie algebra $T^*_x\widehat U_j$ at every $x\in D_j$ is isomorphic to the Lie algebra $\aff$ of the group $\Aff$ introduced below.

The group $\Aff=\{(a,b):a>0,\ b\in\R\}$ has the product
\begin{equation}\label{RADK795GLQ}
(a_2,b_2)(a_1,b_1)=(a_2a_1,\ b_1+a_1b_2),
\end{equation}
with unit $(1,0)$ and inverse $(a,b)^{-1}=(a^{-1},-b/a)$. It is isomorphic to the group of orientation-preserving affine transformations of $\R$ by $(a,b)\mapsto\big(x\mapsto(x-b)/a\big)$. We call $\{(1,b):b\in\R\}$ its \defn{translation subgroup}, and the conjugates of $\{(a,0):a>0\}$ its \defn{dilation subgroups}. The group acts on $\widehat U_j$ on the left by $(a,b)\cdot(t,r)=(t+br,\ ar)$.
Let $\G_j=\Aff\ltimes\widehat U_j\tto\widehat U_j$ be the action groupoid. We write its arrows as $(a,b,t,r)$, so that
\begin{align*}
\sss(a,b,t,r)&=(t,r),\\
\ttt(a,b,t,r)&=(t+br,\ ar),\\
1_{(t,r)}&=(1,0,t,r),\\
(a_2,b_2,t+b_1r,a_1r)\,(a_1,b_1,t,r)&=(a_2a_1,\ b_1+a_1b_2,\ t,\ r).
\end{align*}
As a manifold, $\G_j=\Aff\times\widehat U_j$, so $\G_j$ is Hausdorff, and every source fibre is diffeomorphic to $\R^2$. On $\Aff$ we use the global coordinate $\alpha=\log a$. Since $ar\neq0$ if and only if $r\neq0$, we have $\sss^{-1}(\{r\neq0\})=\ttt^{-1}(\{r\neq0\})$, a dense open subgroupoid of $\G_j$.

Let $\zeta_1,\zeta_2\in\aff$ be the velocities at $0$ of the curves $\alpha\mapsto(e^\alpha,0)$ and $b\mapsto(1,b)$. The Lie algebroid $\Lie(\G_j)$ is the trivial vector bundle $\aff\times\widehat U_j$, and the anchors of $\zeta_1$ and $\zeta_2$ are $r\dd r$ and $r\dd t$. The anchor preserves brackets and is injective on the dense set $\{r\neq0\}$. Since $\Big[r\,\dd r,r\,\dd t\Big]=r\,\dd t$,
we obtain $[\zeta_1,\zeta_2]=\zeta_2$ in $\Lie(\G_j)$ by continuity. Define
\begin{equation}\label{RADK554PBN}
\chi\colon T^*_\sigma\widehat U_j\too\Lie(\G_j),\quad \chi(dr)=\zeta_2,\quad \chi(dt)=-\zeta_1 .
\end{equation}
The map $\chi$ is an isomorphism of vector bundles and intertwines the anchors. It preserves brackets on the global frame, since
\[
[\chi(dr),\chi(dt)]=[\zeta_2,-\zeta_1] =\zeta_2=\chi([dr,dt]).
\]
The Leibniz rule then shows that $\chi$ is an isomorphism of Lie algebroids. As the source fibres of $\G_j$ are diffeomorphic to $\R^2$, this identifies $\G_j$ with the source-simply-connected integration of $T^*_\sigma\widehat U_j$.

\begin{lemma}\label{RADK286UAX}
The $2$-form
\begin{equation}\label{RADK958NSM}
\Omega=dt\wedge d\alpha+db\wedge dr+b\,dr\wedge d\alpha+r\,db\wedge d\alpha
\end{equation}
on $\G_j$ is closed, multiplicative and symplectic, and the target map is Poisson. Thus $(\G_j,\Omega)$ is the source-simply-connected symplectic groupoid integrating $(\widehat U_j,\sigma)$, and $v\mapsto\iota_v\Omega|_{T\widehat U_j}$ is the inverse of $\chi$. Moreover, $\Omega=\ttt^*\omega-\sss^*\omega$ on $\sss^{-1}(\{r\neq0\})$.
\end{lemma}

\begin{proof}
Up to sign, $\Omega$ is the multiplicative symplectic form of Radko and Shlyakhtenko \cite[\S5.2, (5--1)]{radkoshlyakhtenko2006}, written there on $(\R^2,x\,\dd x\wedge\dd y)$, and the sign comes from the opposite convention for the inverse of a bivector. Since $\Omega\wedge\Omega=2\,dt\wedge d\alpha\wedge db\wedge dr$, the form $\Omega$ is nondegenerate. On $\sss^{-1}(\{r\neq0\})$, $d(t+br)=dt+b\,dr+r\,db$ and $d(ar)/ar=d\alpha+dr/r$, so
\begin{align*}
\ttt^*\omega&=\big(dt+b\,dr+r\,db\big)\wedge\Big(d\alpha+\frac{dr}r\Big)\\
&=\Omega+\frac1r\,dt\wedge dr\\
&=\Omega+\sss^*\omega.
\end{align*}
So $\Omega=(\ttt,\sss)^*(\omega\oplus(-\omega))$ on the dense open subgroupoid $\sss^{-1}(\{r\neq0\})$, hence $\Omega$ is closed and multiplicative by continuity. Moreover, $(\ttt,\sss)$ is a local diffeomorphism into $\{r\neq0\}\times\{r\neq0\}$, and $\Omega^{-1}$ is $(\ttt,\sss)$-related to $\sigma\oplus(-\sigma)$. So $\Omega^{-1}$ and $\sigma$ are $\ttt$-related on this dense set, and everywhere by continuity. The groupoid $\G_j$ is source-simply-connected by the paragraph after \eqref{RADK554PBN}. At a unit, $b=0$, $\zeta_1=\dd\alpha$ and $\zeta_2=\dd b$. Formula \eqref{RADK958NSM} therefore gives
\begin{align*}
\iota_{\zeta_1}\Omega|_{T\widehat U_j}&=-dt=\chi^{-1}(\zeta_1),\\
\iota_{\zeta_2}\Omega|_{T\widehat U_j}&=dr=\chi^{-1}(\zeta_2).
\end{align*}
\end{proof}

The adjoint integration is obtained by quotienting $\G_j$ by the following subgroupoid of loops. Its discrete isotropy away from $D_j$ records the modular period $\ell_j$.

\begin{lemma}\label{RADK997RKH}
The subset \[
\NN_j=\big\{(1,0,t,0):t\in S^1_{\ell_j}\big\}\ \sqcup\ \big\{(1,k\ell_j/r,t,r):k\in\Z,\ t\in S^1_{\ell_j},\ r\neq0\big\}
\]
is a closed wide Lie subgroupoid of $\G_j$, with isotropy groups $\NN_{j,(t,0)}=\{(1,0)\}$ and $\NN_{j,(t,r)}=\{(1,k\ell_j/r):k\in\Z\}=\G_{j,(t,r)}$ for $r\neq0$. Moreover, $\sss|_{\NN_j}\colon\NN_j\to\widehat U_j$ is a local diffeomorphism, so $\NN_j$ is discrete.
\end{lemma}

\begin{proof}
For $r\neq0$ and $k\in\Z$, let $g_k(t,r)=(1,k\ell_j/r,t,r)$. Then
\[
\ttt(g_k(t,r))=(t+k\ell_j,r)=(t,r),\quad g_k(t,r)\,g_{k'}(t,r)=g_{k+k'}(t,r),\quad g_k(t,r)^{-1}=g_{-k}(t,r).
\]
Since $\NN_j$ only contains units over $r=0$, it is a wide subgroupoid. For $r\neq0$, $(a,b,t,r)$ is a loop if and only if $ar=r$ and $t+br=t$ in $S^1_{\ell_j}$, that is, $a=1$ and $b\in(\ell_j/r)\Z$.

Fix $n_0=(1,b_0,t_0,r_0)\in\NN_j$. If $r_0\neq0$, let $k_0\in\Z$ with $b_0=k_0\ell_j/r_0$, let $\varepsilon_0=\ell_j/(3|r_0|)$, and set
\[
B=\{(t,r):|r-r_0|<|r_0|/2\},\quad
\widetilde B=\Big\{(a,b,t,r):|a-1|<\tfrac12,\ (t,r)\in B,\ \Big|b-\frac{k_0\ell_j}{r}\Big|<\varepsilon_0\Big\}.
\]
Then $(a,b,t,r)\mapsto(a-1,\,b-k_0\ell_j/r)$ is a submersion from the open neighbourhood $\widetilde B$ of $n_0$ to $\R^2$. For $(t,r)\in B$, consecutive values of $k\ell_j/r$ differ by $\ell_j/|r|>2\ell_j/(3|r_0|)=2\varepsilon_0$, so $\NN_j\cap \widetilde B=\{(1,k_0\ell_j/r,t,r):(t,r)\in B\}$ is the zero set of this submersion. If $r_0=0$, let $B=\{(t,r):|r|<\ell_j/2\}$, $\widetilde B=\{(a,b,t,r):|a-1|<\frac12,\ (t,r)\in B,\ |b|<1\}$, and use the submersion $(a,b,t,r)\mapsto(a-1,b)$. Since $|k\ell_j/r|>2$ for $0<|r|<\ell_j/2$ and $k\neq0$, the set $\NN_j\cap \widetilde B$ is the set of units over $B$, the zero set of this submersion. So $\NN_j$ is an embedded submanifold, hence a Lie subgroupoid. In both cases, $\NN_j\cap \widetilde B$ is the image of a smooth section of $\sss$ over $B$, so $\sss|_{\NN_j}$ is a local diffeomorphism.

For closedness, let $(1,b_i,t_i,r_i)\in\NN_j$ converge to $(a,b,t,r)$. Then $a=1$. If $r_i=0$ for infinitely many $i$, then $b_i=0$ for these $i$, so $(a,b,t,r)=(1,0,t,0)\in\NN_j$. Otherwise, we may assume $r_i\neq0$ and $b_i=k_i\ell_j/r_i$ with $k_i\in\Z$. If $r\neq0$, then $k_i=b_ir_i/\ell_j\to br/\ell_j$, so $k_i=k$ eventually and $b=k\ell_j/r$. If $r=0$, then $k_i\ell_j=b_ir_i\to0$, so $k_i=0$ eventually and $b=0$.
\end{proof}

\begin{lemma}\label{RADK987EJB}
The subgroupoid $\NN_j$ is the greatest element of $\Lambda(\G_j)$.
\end{lemma}

\begin{proof}
By Lemma \ref{RADK997RKH}, $\NN_j$ is discrete and consists of loops. Let $g=(a,b,t,r)$ be an arrow from $x=(t,r)$ to $y=(t+br,ar)$. If $r\neq0$, then $ar\neq0$, so $\NN_{j,x}=\G_{j,x}$ and $\NN_{j,y}=\G_{j,y}$ by Lemma \ref{RADK997RKH}, and $g\,\G_{j,x}\,g^{-1}=\G_{j,y}$. If $r=0$, then $y=x$ and $\NN_{j,x}$ is trivial. In either case $g\,\NN_{j,x}\,g^{-1}=\NN_{j,y}$, so $\NN_j$ is normal.

Let $\NN'\in\Lambda(\G_j)$. Since $\sss=\ttt$ on $\NN'$, $\NN'=\coprod_x\NN'_x$ with $\NN'_x\subset\G_{j,x}$. Since $\sss|_{\NN'}$ is a local diffeomorphism, each $\NN'_x=(\sss|_{\NN'})^{-1}(x)$ is a discrete subgroup of $\G_{j,x}$, normal since $\NN'$ is normal. For $x=(t,r)$ with $r\neq0$, Lemma \ref{RADK997RKH} gives $\NN'_x\subset\G_{j,x}=\NN_{j,x}$. For $x=(t,0)$, $\NN'_x$ is central in the connected group $\G_{j,x}=\Aff$, whose centre is trivial, so $\NN'_x=\{(1,0)\}=\NN_{j,x}$. Hence $\NN'\subset\NN_j$.
\end{proof}

\begin{proposition}\label{RADK640SSJ}
The quotient $\overline\G_j=\G_j/\NN_j\tto\widehat U_j$ is a Hausdorff Lie groupoid, and it is the adjoint integration of $T^*_\sigma\widehat U_j$. The quotient map $q_j\colon\G_j\to\overline\G_j$ is a surjective local diffeomorphism, and there is a unique $2$-form $\overline\Omega_j$ on $\overline\G_j$ with
\begin{equation}\label{RADK150JPC}
q_j^*\,\overline\Omega_j=\Omega .
\end{equation}
The pair $(\overline\G_j,\overline\Omega_j)$ is a symplectic groupoid integrating $(\widehat U_j,\sigma)$.
\end{proposition}

\begin{proof}
By Lemmas \ref{RADK997RKH} and \ref{RADK987EJB}, $\NN_j$ is a closed greatest element of $\Lambda(\G_j)$, and $\G_j$ is the source-simply-connected integration of $T^*_\sigma\widehat U_j$ by the paragraph after \eqref{RADK554PBN}. Since $\G_j$ is Hausdorff, \cite[Theorems 1.14 and 1.20]{gualtieri2014symplectic} show that $\overline\G_j$ is a Hausdorff Lie groupoid and the adjoint integration, with $q_j$ a surjective local diffeomorphism.

For $k\in\Z$, let $R_k(g)=g\cdot(1,k\ell_j/r,t,r)$ on $\sss^{-1}(\{r\neq0\})$, where $\sss(g)=(t,r)$. Then $R_k(a,b,t,r)=(a,\,b+k\ell_j/r,\,t,\,r)$, so the terms $dt\wedge d\alpha$ and $db\wedge dr$ of \eqref{RADK958NSM} are unchanged, while
\[
R_k^*\Omega-\Omega =\frac{k\ell_j}{r}\,dr\wedge d\alpha +r\Big(-\frac{k\ell_j}{r^2}\,dr\Big)\wedge d\alpha =0.
\]
To descend the form, let $\phi$ be a local diffeomorphism between open subsets of $\G_j$ with $q_j\circ\phi=q_j$. Since the fibres of $q_j$ are the $\NN_j$-orbits, $\phi(g)=g\,n(g)$ for the smooth map $n(g)=g^{-1}\phi(g)$ into $\NN_j$, and $\sss\circ n=\sss$. Near $g_0$ with $\sss(g_0)\in\{r\neq0\}$, the integer $k$ with $n(g)=(1,k\ell_j/r,\sss(g))$ is locally constant, so $\phi=R_k$. Near $g_0$ with $\sss(g_0)\in\{r=0\}$, $n(g_0)$ is a unit, so, for $g$ near $g_0$, $n(g)$ lies in the set $\widetilde B$ of the proof of Lemma \ref{RADK997RKH} with $r_0=0$, hence is a unit, and $\phi=\id$. So $\phi^*\Omega=\Omega$.

For open sets $O\subset\G_j$ on which $q_j$ is injective, the forms $((q_j|_O)^{-1})^*\Omega$ on $q_j(O)$ differ on overlaps by a map $\phi$ as above, so they glue to a $2$-form $\overline\Omega_j$ satisfying \eqref{RADK150JPC}, which is unique.

Since $q_j$ is a local diffeomorphism compatible with the structure maps, $\overline\Omega_j$ is closed and nondegenerate and the target map of $\overline\G_j$ is Poisson by Lemma \ref{RADK286UAX}. To check multiplicativity, let $\overline m$ denote multiplication on $\overline\G_j$. The induced map $q_j\times q_j$ between the spaces of composable pairs is a surjective local diffeomorphism, and
\begin{align*}
(q_j\times q_j)^*\overline m^*\overline\Omega_j
&=m^*\Omega\\
&=\pr_1^*\Omega+\pr_2^*\Omega\\
&=(q_j\times q_j)^*\big(\pr_1^*\overline\Omega_j+\pr_2^*\overline\Omega_j\big).
\end{align*}
Hence $\overline\Omega_j$ is multiplicative.
\end{proof}

The gluing in \S\ref{RADK898CDG} uses the restrictions of these groupoids to the open sets $U_I=S^1_{\ell_j}\times I\subset\widehat U_j$, where $I\subset\R$ is an open interval. Since the target of $(a,b,t,r)$ is $(t+br,ar)$, we have $\G_j|_{U_I}=\{(a,b,t,r):r\in I,\ ar\in I\}$.

\begin{lemma}\label{RADK716LDR}
Let $I\subset\R$ be an open interval. Then $\G_j|_{U_I}$ is the source-simply-connected integration of $T^*_\sigma U_I$, and $\NN_j\cap\G_j|_{U_I}$ is a closed greatest element of $\Lambda(\G_j|_{U_I})$. The restriction $\overline\G_j|_{U_I}$ is the quotient $\G_j|_{U_I}/(\NN_j\cap\G_j|_{U_I})$. It is source-connected, it is the adjoint integration of $T^*_\sigma U_I$, and the restriction of $\overline\Omega_j$ is a multiplicative symplectic form on it.
\end{lemma}

\begin{proof}
The restriction $\G_j|_{U_I}$ is an open subgroupoid with Lie algebroid $T^*_\sigma U_I$. Its source fibre over $(t,r)\in U_I$ is $\big\{a>0:ar\in I\big\}\times\R$.
The first factor is an open interval containing $1$, so every source fibre is diffeomorphic to $\R^2$. The subgroupoid $\NN_j\cap\G_j|_{U_I}$ is closed, discrete, totally disconnected and normal. Restriction to $U_I$ does not change the isotropy groups, and the proof of Lemma \ref{RADK987EJB} only uses them.

Since $q_j$ commutes with $\sss$ and $\ttt$, $q_j^{-1}(\overline\G_j|_{U_I})=\G_j|_{U_I}$, and the fibres of the surjective local diffeomorphism $q_j\colon\G_j|_{U_I}\to\overline\G_j|_{U_I}$ are the orbits of $\NN_j\cap\G_j|_{U_I}$. So $\overline\G_j|_{U_I}$ is source-connected, as the image of $\G_j|_{U_I}$, and it is the adjoint integration by \cite[Theorems 1.14 and 1.20]{gualtieri2014symplectic}. The last claim holds since $\overline\G_j|_{U_I}$ is an open subgroupoid of $\overline\G_j$.
\end{proof}

\subsection{The adjoint symplectic groupoid \texorpdfstring{$\G_\Sigma$}{of Sigma}}\label{RADK898CDG}

We now recover the adjoint integration of the surface by gluing these coordinate models. We apply the criterion of \S\ref{RADK237YVT} to $E=T^*_\sigma\Sigma$. This Lie algebroid is integrable since $(\Sigma,\sigma)$ is a proper log symplectic manifold \cite[Theorem 2.15]{gualtieri2014symplectic}. Since $\sigma$ is injective on the dense set $\Sigma\setminus D$, an isomorphism of Lie algebroids $\Lie(\G)\to T^*_\sigma\Sigma$ covering the identity is determined by the anchors. So it is unique if it exists, and we omit it from the notation.

The cover $\{V_i\}_i\cup\{U_j\}_j$ is a finite orbit cover, since the leaves are the finitely many $V_i$ and the points of $D$. Over $V_i$, $\sigma$ is nondegenerate, so the anchor identifies $T^*_\sigma V_i$ with $TV_i$, whose source-simply-connected integration is the fundamental groupoid $\Pi_1(V_i)$. The greatest element of $\Lambda(\Pi_1(V_i))$ is its full isotropy $\NN_{V_i}$, which is closed, and $\Pi_1(V_i)/\NN_{V_i}=V_i\times V_i$ \cite[Example 1.23]{gualtieri2014symplectic}. Over $U_j$, by Lemma \ref{RADK716LDR} with $I=(-c_j,c_j)$, the source-simply-connected integration is $\G_j|_{U_j}$, and $\NN_j\cap\G_j|_{U_j}$ is a closed greatest element of $\Lambda(\G_j|_{U_j})$, with quotient $\overline\G_j|_{U_j}$.

The family $(\NN_{V_i})_i\cup(\NN_j\cap\G_j|_{U_j})_j$ is compatible. Indeed, the $V_i$ are pairwise disjoint, and so are the $U_j$, so, by \eqref{RADK247VWC}, the remaining overlaps are the components $U$ of the sets $U_j\setminus D_j$, with $U\subset V_i$ for a unique $i$. Then $U=U_I$ with $I=(0,c_j)$ or $I=(-c_j,0)$, and $T^*_\sigma U\cong TU$, so the source-simply-connected integration over $U$ is $\Pi_1(U)$, with full isotropy $\Iso(\Pi_1(U))$. Since $\NN_{V_i}$ is the full isotropy and the canonical morphism $\Pi_1(U)\to(\Pi_1(V_i)|_U)^c$ covers the identity, $\res_{V_i,U_j}(\NN_{V_i})=\Iso(\Pi_1(U))$. By Lemma \ref{RADK716LDR}, $\G_j|_U\cong\Pi_1(U)$, and, by Lemma \ref{RADK997RKH}, $\NN_j\cap\G_j|_U$ is the full isotropy of $\G_j|_U$, so $\res_{U_j,V_i}(\NN_j\cap\G_j|_{U_j})=\Iso(\Pi_1(U))$.

\begin{proposition}\label{RADK181QST}
The family $(\NN_{V_i})_i\cup(\NN_j\cap\G_j|_{U_j})_j$ glues to the adjoint integration $\G_\Sigma$ of $T^*_\sigma\Sigma$, and $\G_\Sigma$ is Hausdorff. There is a unique $2$-form $\Omega_\Sigma$ on $\G_\Sigma$ with
\begin{align*}
\G_\Sigma|_{V_i}&=V_i\times V_i, & \Omega_\Sigma|_{V_i\times V_i}&=\ttt^*\omega_i-\sss^*\omega_i,\\
\G_\Sigma|_{U_j}&=\overline\G_j|_{U_j}, & \Omega_\Sigma|_{\overline\G_j|_{U_j}}&=\overline\Omega_j,
\end{align*}
and $(\G_\Sigma,\Omega_\Sigma)$ is a symplectic groupoid integrating $(\Sigma,\sigma)$. Over $U_j$, the isomorphism $\Lie(\G_\Sigma)\to T^*_\sigma\Sigma$, $v\mapsto\iota_v\Omega_\Sigma|_{T\Sigma}$, corresponds under $q_j$ to the inverse of $\chi$.
\end{proposition}

\begin{proof}
By compatibility and the gluing criterion of \S\ref{RADK237YVT}, the family glues to the adjoint integration $\G_\Sigma$, with $(\G_\Sigma|_{V_i})^c=V_i\times V_i$ and $(\G_\Sigma|_{U_j})^c=\overline\G_j|_{U_j}$. Since $\G_\Sigma$ is source-connected, its orbits are the leaves. For $x\in V_i$, the source fibre of $\G_\Sigma|_{V_i}$ over $x$ is the connected fibre $\sss^{-1}(x)$, so $\G_\Sigma|_{V_i}=V_i\times V_i$ and $\ttt$ maps $\sss^{-1}(x)$ bijectively onto $V_i$. If moreover $x\in U_j$, the source fibre of $\G_\Sigma|_{U_j}$ over $x$ is $\sss^{-1}(x)\cap\ttt^{-1}(U_j)\cong V_i\cap U_j$, which is connected by \eqref{RADK247VWC}. If $x\in D_j$, it is $\sss^{-1}(x)$, since the orbit of $x$ is $\{x\}$. Hence $\G_\Sigma|_{U_j}=\overline\G_j|_{U_j}$.

The open subgroupoids $\G_\Sigma|_{V_i}$ and $\G_\Sigma|_{U_j}$ cover $\G_\Sigma$. They also cover composable pairs: both arrows of a composable pair $(g,h)$ lie in $V_i\times V_i$ if $\sss(g)=\ttt(h)\in V_i$, and both are loops in $\overline\G_j|_{U_j}$ if $\sss(g)\in D_j$. The only nonempty overlaps are the $\G_\Sigma|_U$, with $U$ as above, and $(\ttt,\sss)$ is injective on them. There, $\ttt^*\omega_i-\sss^*\omega_i$ and $\overline\Omega_j$ both equal $\ttt^*\omega-\sss^*\omega$, by Lemma \ref{RADK286UAX}, \eqref{RADK150JPC} and $\omega_i=\omega$ on $U$, so they patch to a unique $2$-form $\Omega_\Sigma$.

Closedness, nondegeneracy, multiplicativity and the Poisson property of $\ttt$ can be checked on the pieces, and hold over $U_j$ by Proposition \ref{RADK640SSJ} and Lemma \ref{RADK716LDR}, and over $V_i$ since $\ttt^*\omega_i-\sss^*\omega_i=(\ttt,\sss)^*(\omega_i\oplus(-\omega_i))$ and $\ttt=\pr_1\circ(\ttt,\sss)$, as in the proof of Lemma \ref{RADK286UAX}. The last assertion follows from Lemma \ref{RADK286UAX} and \eqref{RADK150JPC}.

For the Hausdorff property, let $g\neq g'$ in $\G_\Sigma$. If $(\ttt,\sss)(g)\neq(\ttt,\sss)(g')$, then $(\ttt,\sss)$ separates $g$ and $g'$. Otherwise, $\sss(g)=\sss(g')=x$ and $\ttt(g)=\ttt(g')=y$. If $x\in V_i$, then $g=g'=(y,x)$ in $V_i\times V_i$, a contradiction. So $x\in D_j$ for some $j$, and $g,g'$ lie in the open subgroupoid $\overline\G_j|_{U_j}$, which is Hausdorff by Proposition \ref{RADK640SSJ}.
\end{proof}

Gualtieri and Li also obtain the adjoint integration of a proper log symplectic manifold as a symplectic pair groupoid \cite[Theorem 2.15 and Corollary 2.19]{gualtieri2014symplectic}, which is isomorphic to $\G_\Sigma$ by uniqueness of the adjoint integration. At $x\in D_j$, the subgroupoid $\NN_j$ is trivial over $\{r=0\}$, so $q_j$ maps $\G_{j,x}=\Aff$ isomorphically onto $\G_{\Sigma,x}$. We use it to identify $\G_{\Sigma,x}$ with $\Aff$.

Since $\G_\Sigma$ is source-connected, $\ttt(g)$ lies in the leaf through $\sss(g)$, which is $\{\sss(g)\}$ if $\sss(g)\in D$, and $V_i$ if $\sss(g)\in V_i$, in which case $g\in\G_\Sigma|_{V_i}=V_i\times V_i$ by Proposition \ref{RADK181QST}. Hence, for every arrow $g$ of $\G_\Sigma$,
\begin{alignat}{2}
\ttt(g)&=\sss(g) &\quad&\text{if }\sss(g)\in D,\label{RADK472KHB}\\
g&=(\ttt(g),\sss(g))\in V_i\times V_i &\quad&\text{if }\sss(g)\in V_i.\label{RADK618ENP}
\end{alignat}

\section{Lagrangian subgroupoids of admissible curves}\label{RADK793NZP}

In this section, we construct the Lagrangian integrations of curves inside the adjoint symplectic groupoid $\G_\Sigma$. Let $N\subset\Sigma$ be a curve, meaning a closed embedded one-dimensional submanifold. Every such curve is coisotropic, so its conormal bundle $TN^\circ$ is a Lie subalgebroid of $T^*_\sigma\Sigma$. For an admissible curve, we will integrate this subalgebroid to an embedded Lagrangian subgroupoid of $\G_\Sigma$ and prove Theorem \ref{RADK953BQH}. The construction supplies the subgroupoids used in \S\S\ref{RADK389HFP} and \ref{RADK035TQA}.

\subsection{Admissible curves}\label{RADK341ACV}

We first fix the types of components of $N$ and local coordinates at the points of $N\cap D$. The curve $N$ has finitely many components $N_c$, each diffeomorphic to a circle.

\begin{definition}\label{RADK076ZXX}
A curve $N$ is \defn{admissible} if every component $N_c$ is of one of the following types:
\begin{enumerate}
\item $N_c\subset\Sigma\setminus D$;
\item $N_c=D_j$ for some $j$;
\item $N_c\cap D$ is a nonempty finite set.
\end{enumerate}
The curve $N$ is \defn{tangent to $D$} at $x\in N\cap D$ if $T_xN=T_xD$. For a component $N_c$ of type (3), a \defn{contact} at $x\in N_c\cap D$ is \defn{tangential} if $T_xN=T_xD$, and \defn{transverse} otherwise.
\end{definition}

\begin{figure}[ht]
\centering
\includegraphics[width=0.68\linewidth]{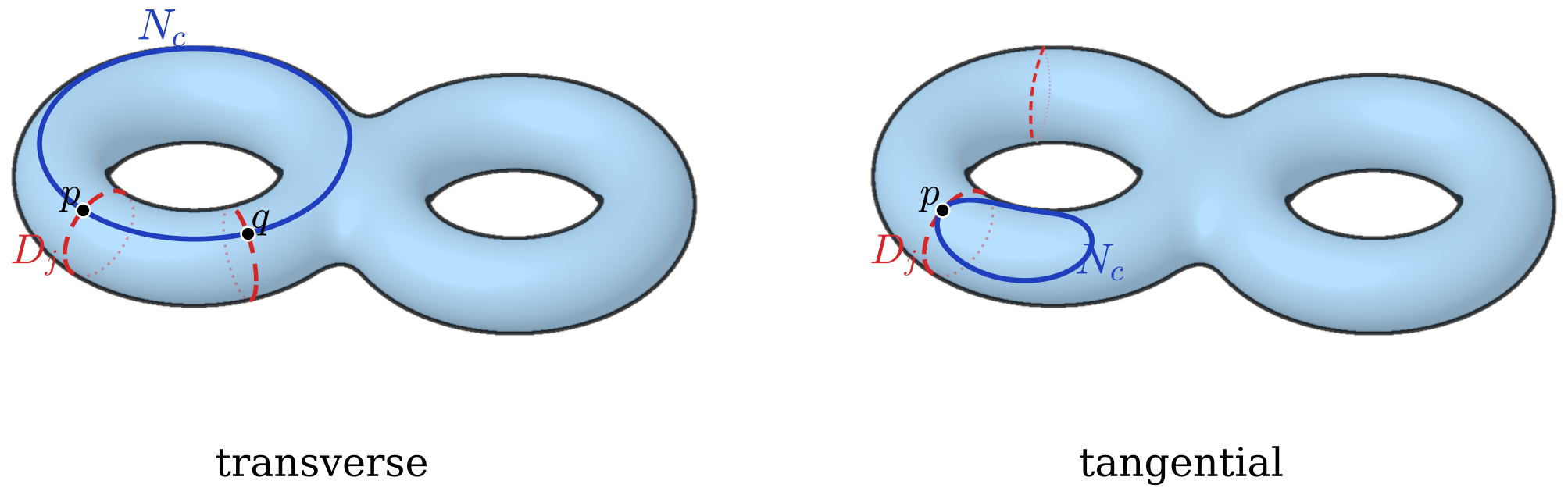}
\caption{Transverse and tangential contacts.}\label{RADK702FGB}
\end{figure}

Figure \ref{RADK702FGB} shows a component with two transverse contacts $p,q$, and a component with a single tangential contact $p$. From now on, $N$ is admissible. For $N_c$ of type (3), we write $\mathcal C_c=N_c\cap D$. Let $N_0$ and $N_D$ be the unions of the components of type (1) and (2), and let $\mathcal C$ be the union of the sets $\mathcal C_c$, which is finite. For $N_c$ of type (3), removing the $|\mathcal C_c|\ge1$ points of $\mathcal C_c$ from the circle $N_c$ leaves $|\mathcal C_c|$ open arcs, contained in $\Sigma\setminus D$. Writing $\K$ for the set of all these arcs, we obtain
\begin{equation}\label{RADK376JAG}
N=N_0\ \sqcup\ N_D\ \sqcup\ \Big(\coprod_{K\in\K}K\Big)\ \sqcup\ \mathcal C .
\end{equation}

A \defn{box} around $p\in D_j$ is an open set $W=J\times I\subset U_j$ containing $p$, where $J\subset S^1_{\ell_j}$ is an open arc of length less than $\ell_j$ and $I\subset(-c_j,c_j)$ is an open interval. Choose a lift $\widetilde J\subset\R$ of $J$. We identify $J$ with this lift whenever we use real coordinates, take derivatives, or write expressions such as $t+br$ and $t-\gamma(r)$. In particular, a condition such as $t+br\in\widetilde J$ refers to the real-valued expression, before reduction modulo $\ell_j$.

\begin{lemma}\label{RADK227HCV}
Let $N_c$ be of type (3), and let $p\in N_c\cap D_j$ be a transverse contact. Translating $t$, assume $p=(0,0)$. There are a box $W=J\times I$ around $p$ and a smooth function $\gamma\colon I\to J$ with $\gamma(0)=0$ such that, in the coordinates $(s,r)$ on $W$ with $s=t-\gamma(r)$,
\[
\sigma=r\,\dd r\wedge\dd s,\quad D_j\cap W=\{r=0\},\quad N_c\cap W=\{s=0\}=\big\{(\gamma(r),r):r\in I\big\}.
\]
\end{lemma}

\begin{lemma}\label{RADK562YPM}
Let $N_c$ be of type (3), and let $p\in N_c\cap D_j$ be a tangential contact. Translating $t$, assume $p=(0,0)$. There are a box $W=J\times I$ around $p$, with $I=(-\varepsilon,\varepsilon)$, and a smooth function $h\colon J\to I$ with $h(0)=h'(0)=0$ and $h(t)\neq0$ for $t\neq0$, such that in the coordinates \eqref{RADK871UDA}
\[
D_j\cap W=\{r=0\},\quad N_c\cap W=\big\{(t,h(t)):t\in J\big\}.
\]
\end{lemma}

\begin{proof}[Proof of Lemmas \ref{RADK227HCV} and \ref{RADK562YPM}]
Since $T_pD_j=\langle\dd t\rangle$, the restriction of $dr$ to $T_pN_c$ is nonzero at a transverse contact, while that of $dt$ is nonzero at a tangential contact. The inverse function theorem therefore gives an open arc $N_p\subset N_c$ containing $p$ which is a graph over $r$ in the first case and over $t$ in the second.

At a transverse contact, write this graph as $t=\gamma(r)$ for $r\in I_0$, with $\gamma(0)=0$. Since $N_c\setminus N_p$ is closed and misses $p$, there is a box $J\times I_1$ around $p$ disjoint from this closed set, with $I_1\subset I_0$. Shrink $I_1$ to an interval $I\ni0$ such that $\gamma(I)\subset J$. Then $W=J\times I$ satisfies $N_c\cap W=\big\{(\gamma(r),r):r\in I\big\}$.
The change of coordinates $s=t-\gamma(r)$ preserves the normal-form bivector: $\dd r\wedge\dd s=\dd r\wedge\dd t$,
where the first $\dd r$ is taken in the system $(s,r)$.

At a tangential contact, write the graph as $r=h(t)$ for $t\in J_0$. Tangency gives
\[
\big\langle\dd t+h'(0)\dd r\big\rangle=T_pN_c=\big\langle\dd t\big\rangle,
\]
so $h(0)=h'(0)=0$. Choose a box $J_1\times(-\varepsilon,\varepsilon)$ disjoint from $N_c\setminus N_p$, with $J_1\subset J_0$. Because $N_c\cap D$ is finite, we may shrink $J_1$ to an arc $J\ni0$ such that
\[
h(J)\subset I=(-\varepsilon,\varepsilon),\quad h^{-1}(0)=\{0\}.
\]
For $W=J\times I$, this gives $N_c\cap W=\{(t,h(t)):t\in J\}$, as required.
\end{proof}

From now on, $\dd r$ denotes the coordinate field of $r$ in whichever of the two systems is in use.

\subsection{Local Lagrangian subgroupoids}\label{RADK352LLS}

We construct the integration separately over the pieces of \eqref{RADK376JAG}. Away from the finite contact set $\mathcal C$, the models are pair groupoids and bundles of translation groups. At transverse and tangential contacts, we will construct smooth local subgroupoids joining the adjacent pair groupoids to the appropriate isotropy groups. We use \eqref{RADK618ENP} throughout.

Let $K\subset\Sigma\setminus D$ be a connected embedded one-dimensional submanifold, not necessarily closed. It lies in a single $V_i$, over which $\G_\Sigma=V_i\times V_i$ with the form $\ttt^*\omega_i-\sss^*\omega_i$ (Proposition \ref{RADK181QST}). Since $K$ is Lagrangian in $V_i$, the anchor $\sigma\colon T^*_\sigma V_i\to TV_i$ identifies $TK^\circ$ with $TK$. The pair groupoid $K\times K$ therefore integrates $TK^\circ$, is embedded, and has connected source fibres. It is Lagrangian because
\[
\big(\ttt^*\omega_i-\sss^*\omega_i\big)|_{K\times K}
=\ttt^*(\omega_i|_K)-\sss^*(\omega_i|_K)=0.
\]

Near $D_j$, we first construct the local models in $\G_j$ and then pass through $q_j$ to $\G_\Sigma$. Because a local diffeomorphism can identify distinct points, we will use the following \emph{embedding criterion}. Let $O\subset\G_j$ be open, and let $L\subset O$ be an embedded submanifold with
\begin{equation}\label{RADK284EQW}
q_j^{-1}(q_j(g))\cap O=\{g\}\quad\text{for every }g\in L.
\end{equation}
Then $q_j|_L$ is an injective immersion, and $q_j(L\cap O')=q_j(L)\cap q_j(O')$ for every open $O'\subset O$. So $q_j|_L$ is a homeomorphism onto its image, and $q_j(L)$ is an embedded submanifold of $\overline\G_j$, with $q_j|_L$ a diffeomorphism onto it \cite[Proposition 5.2]{lee2003smooth}.

The fibres of $q_j$ are the orbits of the right action of $\NN_j$ (proof of Proposition \ref{RADK640SSJ}). So, by Lemma \ref{RADK997RKH} and \eqref{RADK795GLQ}, if $g\neq g'$ in $\G_j$ and $q_j(g)=q_j(g')$, then $g$ and $g'$ have the same coordinates $a$, $t$ and $r\neq0$, and $b$-coordinates differing by $k\ell_j/r$ with $k\in\Z\setminus\{0\}$. In particular, $q_j$ is not injective on the restriction $\G_j|_W$ to a box $W$: for $r\in I\setminus\{0\}$, the arrows $(a,b,t,r)$ and $(a,b+\ell_j/r,t,r)$ have the same source and target, so one lies in $\G_j|_W$ if and only if the other does, and they have the same image under $q_j$. Requiring both endpoint coordinates to lie in the chosen lift $\widetilde J$ selects a unique preimage, as the next lemma shows.

\begin{lemma}\label{RADK788LQC}
Let $W=J\times I$ be a box around a point of $D_j$. Then
\[
\G^{\mathrm{loc}}_{j,W}=\big\{(a,b,t,r)\in\G_j:t\in\widetilde J,\ t+br\in\widetilde J,\ r\in I,\ ar\in I\big\}
\]
is an open subset of $\G_j|_W$ which contains the isotropy groups $\G_{j,x}$ for $x\in D_j\cap W$, and $q_j$ restricts to a diffeomorphism from $\G^{\mathrm{loc}}_{j,W}$ onto $\overline\G_j|_W$.
\end{lemma}

\begin{proof}
The conditions are open, and they hold on $\G_{j,x}$ since $r=0$ there. If $g\neq g'$ in $\G^{\mathrm{loc}}_{j,W}$ had the same image, then $r\neq0$ and the values of $t+br$ for $g$ and $g'$ would differ by $k\ell_j$ with $k\neq0$, which is impossible since both lie in $\widetilde J$. So $q_j$ is injective on $\G^{\mathrm{loc}}_{j,W}$. Conversely, an arrow of $\overline\G_j|_W$ is $q_j(a,b,t,r)$ with $t\in\widetilde J$, $r\in I$, $ar\in I$ and $t+br\in\widetilde J+\ell_j\Z$. If $r=0$, then $t+br=t\in\widetilde J$. If $r\neq0$, then $t+br-k\ell_j\in\widetilde J$ for some $k\in\Z$, and $(a,b-k\ell_j/r,t,r)=(a,b,t,r)\,(1,-k\ell_j/r,t,r)$ lies in $\G^{\mathrm{loc}}_{j,W}$ and has the same image. So $q_j$ maps $\G^{\mathrm{loc}}_{j,W}$ bijectively onto $\overline\G_j|_W$, hence diffeomorphically, since $q_j$ is a local diffeomorphism.
\end{proof}

\begin{proposition}\label{RADK001UFG}
Let $N_c=D_j$, and let $\widehat\Lc_c=\{(1,b,t,0):b\in\R,\ t\in S^1_{\ell_j}\}\subset\G_j$. Then $\Lc_c=q_j(\widehat\Lc_c)$ is an embedded Lagrangian subgroupoid of $\G_\Sigma$ over $D_j$ integrating $TD_j^\circ$, with source fibres $\R$.
\end{proposition}

\begin{proof}
In the coordinates \eqref{RADK871UDA}, we have $TD_j^\circ=\langle dr\rangle$ and $\chi(dr)=\zeta_2$ by \eqref{RADK554PBN}. All arrows of $\widehat\Lc_c$ are loops, with multiplication $(1,b_2,t,0)\,(1,b_1,t,0)=(1,b_1+b_2,t,0)$.
Thus $\widehat\Lc_c$ is a Lie subgroupoid with source fibres $\R$ and Lie algebroid $\langle\zeta_2\rangle=\chi(TD_j^\circ)$.

Since $r=0$ on $\widehat\Lc_c$, every $g\in\widehat\Lc_c$ is the only point in its $q_j$-fibre. The embedding criterion \eqref{RADK284EQW} therefore applies with $O=\G_j$. Moreover, every term of \eqref{RADK958NSM} contains $d\alpha$ or $dr$, so $\Omega|_{\widehat\Lc_c}=0$. As $\widehat\Lc_c$ has dimension $2$, it is Lagrangian. Finally, $\widehat\Lc_c\subset\G_j|_{U_j}$, so Proposition \ref{RADK181QST} and \eqref{RADK150JPC} give $\Lc_c\subset\G_\Sigma$ and $\Omega_\Sigma|_{\Lc_c}=0$.
\end{proof}

At a transverse contact, consider arrows $(a,b,\gamma(r),r)$ with $r\neq0$ and $r,ar\in I$. For representatives whose real target coordinate $\gamma(r)+br$ lies in the chosen lift of $J$, the target belongs to $N_c\cap W$ precisely when $\gamma(r)+br=\gamma(ar)$.
Thus $b$ is the divided difference $\beta(a,r)$ defined below. Its smooth extension across $r=0$ gives coordinates that straighten the curve inside $\G_j$.

\begin{lemma}\label{RADK450RVY}
Let $W=J\times I$ and $\gamma$ be as in Lemma \ref{RADK227HCV}, and set
\[
\beta(a,r)=
\begin{cases}
\dfrac{\gamma(ar)-\gamma(r)}{r}, & r\neq0,\\
(a-1)\,\gamma'(0), & r=0.
\end{cases}
\]
Then $\beta$ is smooth on $\{(a,r):a>0,\ r\in I,\ ar\in I\}$, and the functions $a$, $b'=b-\beta(a,r)$, $s=t-\gamma(r)$ and $r$ are coordinates on $O_W=\{(a,b,t,r)\in\G_j:t\in J,\ r\in I,\ ar\in I\}$. In these coordinates,
\begin{align*}
\sss(a,b',s,r)&=(s,r),\quad \ttt(a,b',s,r)=(s+b'r,\ ar),\\
(a_2,b'_2,s+b'_1r,a_1r)\,(a_1,b'_1,s,r)&=(a_2a_1,\ b'_1+a_1b'_2,\ s,\ r),\\
\Omega&=ds\wedge d\alpha+db'\wedge dr+b'\,dr\wedge d\alpha+r\,db'\wedge d\alpha .
\end{align*}
Here $(s+b'r,ar)$ stands for the point $\big(s+b'r+\gamma(ar),ar\big)$ of $\widehat U_j$, with first coordinate read modulo $\ell_j$.
\end{lemma}

\begin{proof}
Hadamard's lemma gives $\gamma(u)=\gamma'(0)u+u^2\gamma_2(u)$ for a smooth function $\gamma_2$ on $I$. Hence
\[
\beta(a,r)=\gamma'(0)(a-1)+r\big(a^2\gamma_2(ar)-\gamma_2(r)\big)
\]
for every $r$, proving smoothness. The coordinate change has the smooth inverse
\[
(a,b',s,r)\longmapsto\big(a,b'+\beta(a,r),s+\gamma(r),r\big).
\]

The identity $r\,\beta(a,r)=\gamma(ar)-\gamma(r)$ gives $t+br=(s+b'r)+\gamma(ar)$,
which proves the formulas for source and target. Substituting $b_i=b'_i+\beta(a_i,r_i)$, with $r_1=r$ and $r_2=a_1r$, into \eqref{RADK795GLQ} reduces the multiplication formula to $\beta(a_1,r)+a_1\,\beta(a_2,a_1r)=\beta(a_2a_1,r)$.
For $r\neq0$, the divided differences telescope. At $r=0$, the identity becomes
\[
(a_1-1)\gamma'(0)+a_1(a_2-1)\gamma'(0)=(a_2a_1-1)\gamma'(0).
\]

Since $dt=ds+\gamma'(r)\,dr$, $\omega=\frac1r\,ds\wedge dr$ on $\{r\neq0\}$. So the computation of Lemma \ref{RADK286UAX}, with $(t,b)$ renamed $(s,b')$, gives the formula for $\Omega=\ttt^*\omega-\sss^*\omega$ on $\{r\neq0\}$, hence on $O_W$ by density.
\end{proof}

\begin{figure}[ht]
\centering
\includegraphics[width=0.53\linewidth]{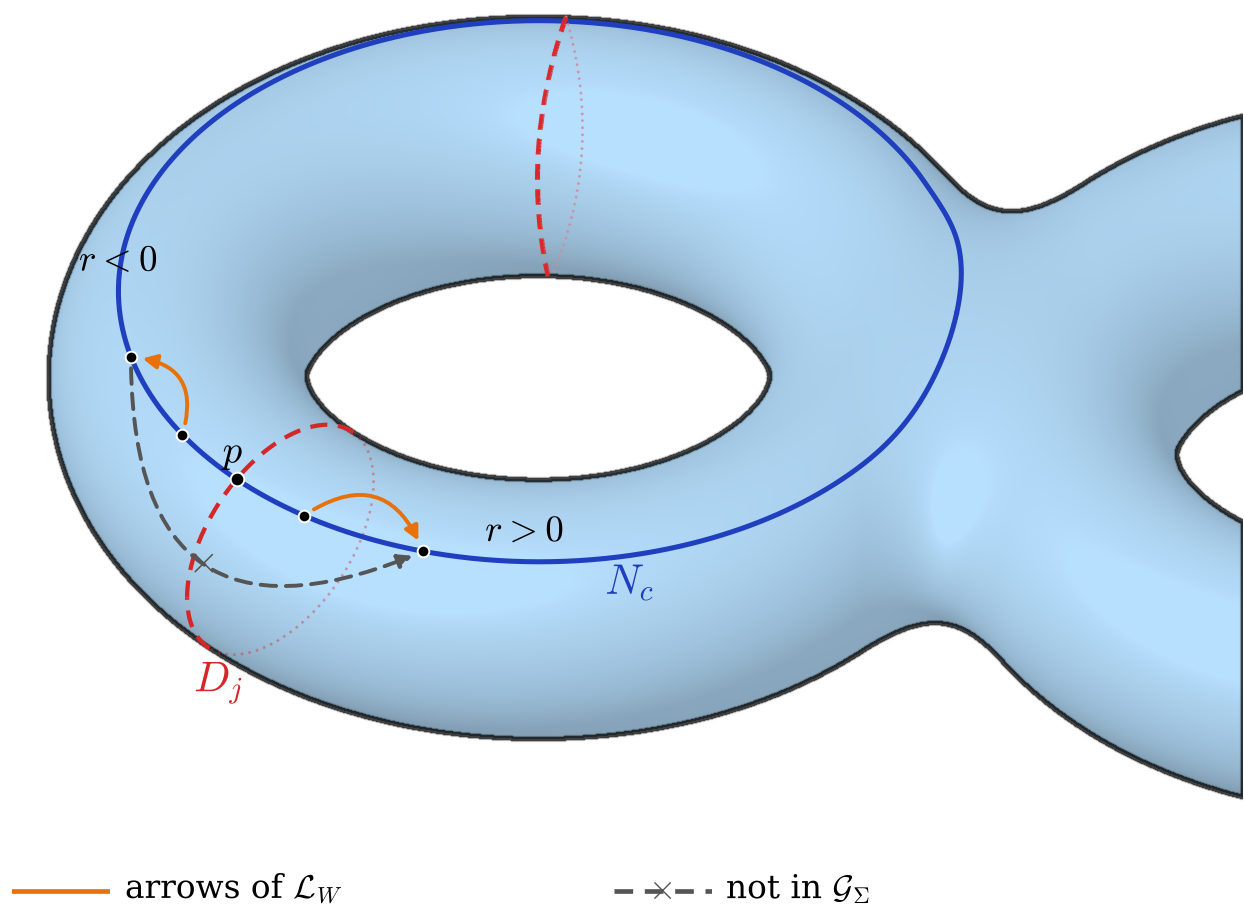}
\caption{A transverse contact.}\label{RADK703FGC}
\end{figure}

\begin{proposition}\label{RADK616RVY}
Let $p\in N_c\cap D_j$ be a transverse contact, with $W=J\times I$, $\gamma$ and the coordinates $(a,b',s,r)$ of Lemmas \ref{RADK227HCV} and \ref{RADK450RVY}. Let $\widehat\Lc_W=\{(a,0,0,r):a>0,\ r\in I,\ ar\in I\}\subset O_W$. Then $\Lc_W=q_j(\widehat\Lc_W)$ is an embedded, source-connected Lagrangian subgroupoid of $\G_\Sigma$ over $N_c\cap W$ integrating $TN^\circ$, with isotropy $\R_{>0}$ at $p$. Moreover, the arrows of $\Lc_W$ with source different from $p$ are
\begin{equation}\label{RADK654KCY}
\Big\{\big((\gamma(r'),r'),(\gamma(r),r)\big):r,r'\in I\setminus\{0\},\ \sgn r=\sgn r'\Big\},
\end{equation}
that is, the pairs of points of $(N_c\cap W)\setminus\{p\}$ on the same side of $p$ (Figure \ref{RADK703FGC}).
\end{proposition}

\begin{proof}
The set $\widehat\Lc_W=\{b'=s=0\}$ is a closed embedded surface in $O_W$, parametrized by $(a,r)$. By Lemma \ref{RADK450RVY},
\[
\sss(a,0,0,r)=(0,r),\quad \ttt(a,0,0,r)=(0,ar),\quad (a_2,0,0,a_1r)\,(a_1,0,0,r)=(a_2a_1,0,0,r).
\] So $\widehat\Lc_W$ is a Lie subgroupoid over $\{s=0\}=N_c\cap W$, with source fibre $\{a>0:ar\in I\}$ over $(0,r)$, equal to $\R_{>0}$ over $p$.

In the coordinates $(a,b,t,r)$, the curve $\alpha\mapsto(e^\alpha,\beta(e^\alpha,r),\gamma(r),r)$ lies in the source fibre of $\widehat\Lc_W$ through $1_{(\gamma(r),r)}$, and its velocity at $\alpha=0$ is $\zeta_1+\gamma'(r)\zeta_2$ since $\frac{\partial\beta}{\partial a}(a,r)=\gamma'(ar)$ for every $r$. Along $N_c\cap W$, we have
\[
\chi(TN_c^\circ)=\chi\big\langle dt-\gamma'(r)\,dr\big\rangle
=\big\langle\zeta_1+\gamma'(r)\,\zeta_2\big\rangle.
\]
Thus $\widehat\Lc_W$ integrates $TN^\circ$ over $N_c\cap W$. Every term of $\Omega$ in Lemma \ref{RADK450RVY} contains $ds$, $db'$ or $b'$, so $\Omega|_{\widehat\Lc_W}=0$, and $\widehat\Lc_W$ is Lagrangian since it has dimension $2$.

Let $\bar r=\sup_I|r|$ and $O=O_W\cap\{|b'|<\ell_j/(2\bar r)\}\supset\widehat\Lc_W$. Two distinct points of $O$ in one fibre of $q_j$ would have the same $a$, $s$ and $r$, and $b'$-coordinates differing by at least $\ell_j/\bar r$, which is impossible. So \eqref{RADK284EQW} holds, $\Lc_W$ is embedded in $\overline\G_j|_{U_j}\subset\G_\Sigma$, and $\Omega_\Sigma|_{\Lc_W}=0$ by \eqref{RADK150JPC}.

For $g=(a,0,0,r)$ with $r\neq0$, $\sss(g)=(\gamma(r),r)$ and $\ttt(g)=(\gamma(ar),ar)$, and $\{ar:a>0,\ ar\in I\}=I\cap r\R_{>0}$. Hence \eqref{RADK654KCY}, by \eqref{RADK618ENP}.
\end{proof}

At a tangential contact, an arrow $(a,b,t,h(t))$ has target $\big(t+b\,h(t),a\,h(t)\big)$. When $t+b\,h(t)$ lies in the chosen lift of $J$, the target belongs to $N_c\cap W$ precisely when $a\,h(t)=h\big(t+b\,h(t)\big)$.
For $t\neq0$, this determines $a$ as a function of $(t,b)$. The next lemma extends that function smoothly across $t=0$ on a domain that keeps source and target on the same side of the contact. We denote the target coordinate by $\tau(t,b)$ and the extended value of $a$ by $\eta(t,b)$.

\begin{lemma}\label{RADK400HUH}
Let $p$, $W=J\times I$ and $h$ be as in Lemma \ref{RADK562YPM}. Write $h(t)=t\,h_1(t)$ with $h_1$ smooth, and set
\[
\tau(t,b)=t+b\,h(t),\quad
S=\big\{(t,b)\in J\times\R:\ \tau(t,b)\in J,\ 1+b\,h_1(t)>0\big\}.
\]
For $t\in J$, write $S_t=\{b\in\R:(t,b)\in S\}$, and for $t\neq0$ set $J_t=\{u\in J:ut>0\}$. Then $S$ is open, $S_0=\R$, and for $t\neq0$ the slice $S_t$ is an interval containing $0$, which $b\mapsto\tau(t,b)$ maps bijectively onto $J_t$. Moreover, there is a smooth function $\eta\colon S\to\R_{>0}$ with $\eta(0,\cdot)\equiv1$ and
\begin{equation}\label{RADK084GSA}
\eta(t,b)\,h(t)=h\big(\tau(t,b)\big),\quad (t,b)\in S .
\end{equation}
\end{lemma}

\begin{proof}
The set $S$ is open, and $S_0=\R$ since $\tau(0,b)=0\in J$ and $h_1(0)=h'(0)=0$. For $t\neq0$, $\tau(t,b)=t\,(1+b\,h_1(t))$ with $h_1(t)\neq0$, so $b\mapsto\tau(t,b)$ is an affine bijection of $\R$, and $1+b\,h_1(t)>0$ exactly when $t\,\tau(t,b)>0$. So it maps $S_t$ onto the interval $J_t$, and $S_t$ is an interval containing $0$ since $\tau(t,0)=t\in J_t$.

To construct $\eta$, work first on the larger open set $S'=\big\{(t,b)\in J\times\R:\tau(t,b)\in J\big\}\supset S$
and set $\eta(t,b)=1+b\int_0^1h'\big(t+u\,b\,h(t)\big)\,du$.
This is smooth, since $t+ubh(t)=(1-u)t+u\,\tau(t,b)\in J$. The fundamental theorem of calculus, applied to $u\mapsto h(t+ubh(t))$, gives \eqref{RADK084GSA}, and $\eta(0,\cdot)\equiv1$ since $h'(0)=0$.

If $t\neq0$ and $\eta(t,b)=0$, then $h(\tau(t,b))=0$ by \eqref{RADK084GSA}, so $\tau(t,b)=0$ and $1+b\,h_1(t)=0$, which is excluded on $S$. So $\eta$ has no zero on the interval $S_t$, and $\eta>0$ there since $\eta(t,0)=1$.
\end{proof}

\begin{figure}[ht]
\centering
\includegraphics[width=0.53\linewidth]{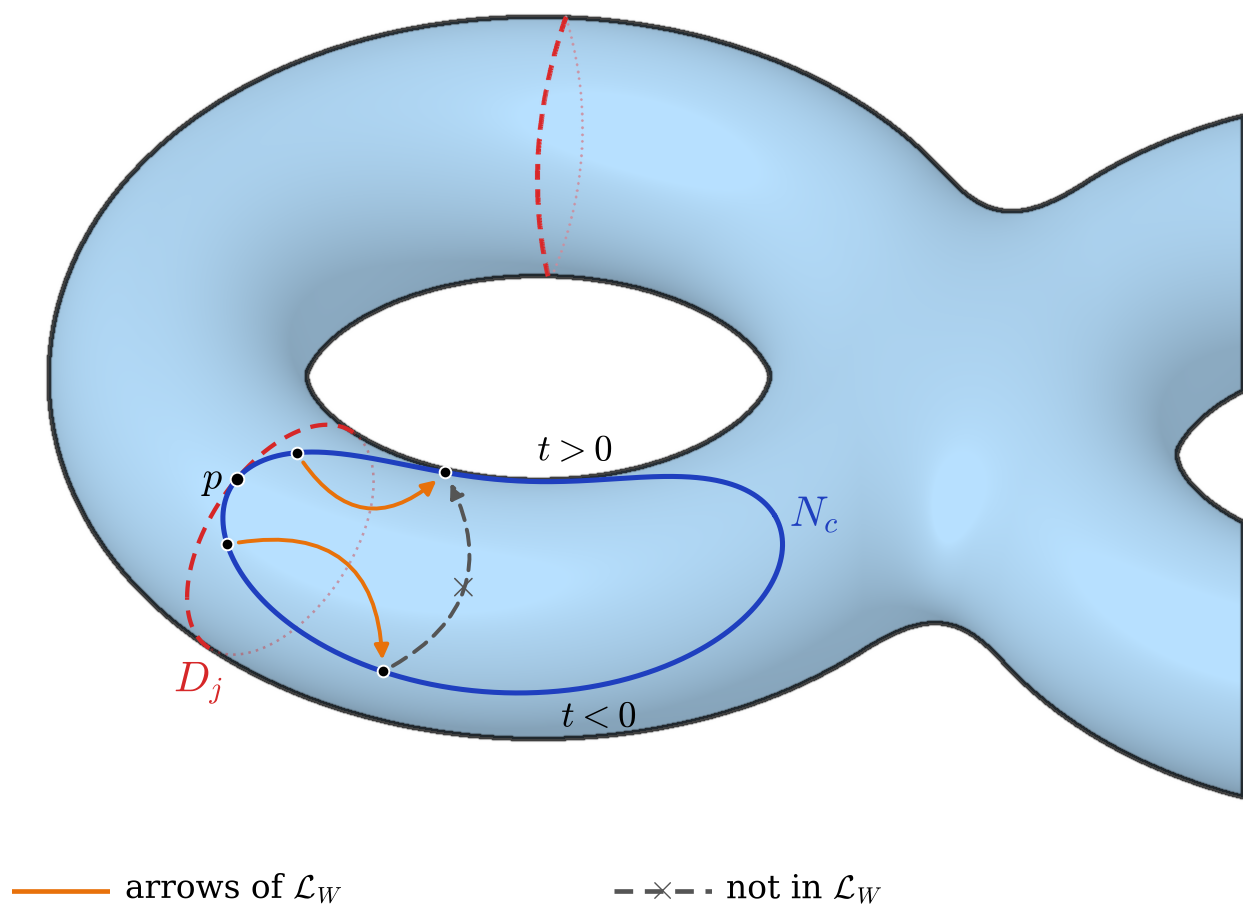}
\caption{A tangential contact.}\label{RADK704FGD}
\end{figure}

\begin{proposition}\label{RADK318JCA}
Let $p$, $W=J\times I$ and $h$ be as in Lemma \ref{RADK562YPM}, and let $\tau$, $\eta$ and $S$ be as in Lemma \ref{RADK400HUH}. Let $\widehat\Lc_W=\{g(t,b)=(\eta(t,b),b,t,h(t)):(t,b)\in S\}\subset\G_j$. Then $\Lc_W=q_j(\widehat\Lc_W)$ is an embedded, source-connected Lagrangian subgroupoid of $\G_\Sigma$ over $N_c\cap W$ integrating $TN^\circ$, with isotropy $\R$ at $p$. Moreover, the arrows of $\Lc_W$ with source different from $p$ are
\begin{equation}\label{RADK821EHP}
\Big\{\big((t',h(t')),(t,h(t))\big):t\in J\setminus\{0\},\ t'\in J_t\Big\},
\end{equation}
that is, the pairs of points of $(N_c\cap W)\setminus\{p\}$ on the same side of $p$ (Figure \ref{RADK704FGD}).
\end{proposition}

\begin{proof}
The set $\widehat\Lc_W$ is an embedded surface in $\G_j$, as the graph over $S$ of $(t,b)\mapsto(\eta(t,b),h(t))$. By \eqref{RADK084GSA}, $g(t,b)$ is an arrow from $(t,h(t))$ to $(\tau(t,b),h(\tau(t,b)))$, both in $N_c\cap W$, and $g(t,0)$ is a unit.

As $J$ has length less than $\ell_j$, composable arrows of $\widehat\Lc_W$ are $g(t_1,b_2)$ and $g(t,b_1)$ with $t_1=\tau(t,b_1)$. With $\eta_1=\eta(t,b_1)$,
\begin{align*}
\tau(t,b_1+\eta_1b_2)&=t_1+b_2\,\eta_1h(t)=t_1+b_2\,h(t_1)=\tau(t_1,b_2),\\
\eta(t,b_1+\eta_1b_2)&=\frac{h(\tau(t_1,b_2))}{h(t_1)}\cdot\frac{h(t_1)}{h(t)}=\eta(t_1,b_2)\,\eta_1
\quad(t\neq0),
\end{align*}
and both sides of the last identity equal $1$ for $t=0$. Moreover, $(t,b_1+\eta_1b_2)\in S$ since $\tau(t_1,b_2)\in J_{t_1}=J_t$ for $t\neq0$. Thus \eqref{RADK795GLQ} gives
\[
g(t_1,b_2)\,g(t,b_1)=g(t,b_1+\eta_1b_2),\quad g(t,b_1)^{-1}=g(t_1,-b_1/\eta_1).
\]
The parameter of the inverse belongs to $S$, since $\tau(t_1,-b_1/\eta_1)=t_1-b_1h(t)=t$.
Hence $\widehat\Lc_W$ is a Lie subgroupoid, with source fibres the intervals $S_t$ and isotropy $\{(1,b,0,0):b\in\R\}\cong\R$ at $p$.

At $x=(t_0,h(t_0))$, the source fibre through $g(t_0,0)$ has tangent vector $\frac{d}{db}\big|_{b=0}g(t_0,b)=\frac{\partial\eta}{\partial b}(t_0,0)\,\zeta_1+\zeta_2$. Differentiating \eqref{RADK084GSA} in $b$ and in $t$ gives
\begin{equation}\label{RADK315DRV}
\begin{aligned}
\frac{\partial\eta}{\partial b}\,h(t)&=h'(\tau)\,h(t),\\
\frac{\partial\eta}{\partial t}\,h+\eta\,h'&=h'(\tau)\,(1+b\,h'),
\end{aligned}
\end{equation}
where $h$ and $h'$ are evaluated at $t$. So $\frac{\partial\eta}{\partial b}=h'(\tau)$ where $t\neq0$, hence everywhere by continuity, and $\frac{\partial\eta}{\partial b}(t_0,0)=h'(t_0)$. On the other hand,
\[
T_xN_c^\circ=\big\langle h'(t_0)\,dt-dr\big\rangle,\quad \chi\big(h'(t_0)\,dt-dr\big)=-h'(t_0)\,\zeta_1-\zeta_2.
\]
Thus the tangent line to the source fibre is $\chi(T_xN_c^\circ)$, and $\widehat\Lc_W$ integrates $TN^\circ$ over $N_c\cap W$.

Along $\widehat\Lc_W$, $dr=h'\,dt$ and $\eta\,d\alpha=\frac{\partial\eta}{\partial t}\,dt+\frac{\partial\eta}{\partial b}\,db$. So, by \eqref{RADK958NSM} and \eqref{RADK315DRV}, the coefficient $C$ of $dt\wedge db$ in $\Omega|_{\widehat\Lc_W}$ satisfies
\begin{align*}
\eta\,C
&=\frac{\partial\eta}{\partial b}\,(1+b\,h')
-\Big(\frac{\partial\eta}{\partial t}\,h+\eta\,h'\Big)\\
&=h'(\tau)\,(1+b\,h')-h'(\tau)\,(1+b\,h')\\
&=0.
\end{align*}
So $\Omega|_{\widehat\Lc_W}=0$, and $\widehat\Lc_W$ is Lagrangian since it has dimension $2$.

The surface $\widehat\Lc_W$ lies in the open set $\G^{\mathrm{loc}}_{j,W}$ of Lemma \ref{RADK788LQC}, since $\tau(t,b)\in\widetilde J$ and $\eta(t,b)\,h(t)=h(\tau(t,b))\in I$, and $q_j$ is injective on $\G^{\mathrm{loc}}_{j,W}$. So \eqref{RADK284EQW} shows that $\Lc_W$ is embedded in $\overline\G_j|_{U_j}\subset\G_\Sigma$, with $\Omega_\Sigma|_{\Lc_W}=0$ by \eqref{RADK150JPC}.

For $t\neq0$, $\tau(t,\cdot)$ maps $S_t$ onto $J_t$ by Lemma \ref{RADK400HUH}, which gives \eqref{RADK821EHP} by \eqref{RADK618ENP}.
\end{proof}

\subsection{The Lagrangian subgroupoid of an admissible curve}\label{RADK374GLB}

We now assemble the local models of \S\ref{RADK352LLS} into a subgroupoid over the whole curve. The only point requiring care is embeddedness near the isotropy groups at the contacts. Enumerate the points of $\mathcal C$ as $p_1,\dots,p_{n_{\mathcal C}}$. For each $i$, let $D_j$ and $N_c$ be the components of $D$ and $N$ containing $p_i$, and let $W_i$ be a box around $p_i$ as in Lemma \ref{RADK227HCV} or Lemma \ref{RADK562YPM}. The components of $N$ other than $N_c$ are closed and miss $p_i$. So, shrinking $I$ and $J$ as in the proof of Lemmas \ref{RADK227HCV} and \ref{RADK562YPM}, we may assume that $N\cap W_i=N_c\cap W_i$. Then $N\cap D\cap W_i=\{p_i\}$.

\begin{figure}[ht]
\centering
\includegraphics[width=0.77\linewidth]{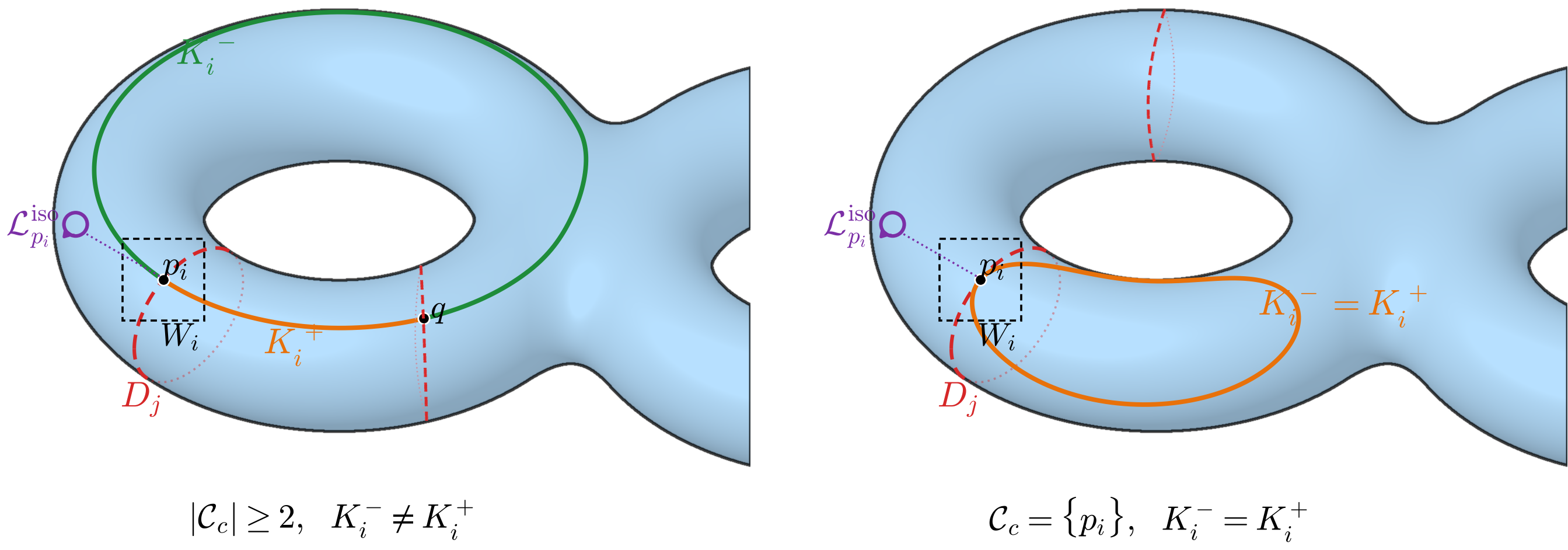}
\caption{The pieces of $\Lc$ near a contact point $p_i$.}\label{RADK705FGE}
\end{figure}

Let $\Lc_{W_i}$ be the local Lagrangian subgroupoid over $N_c\cap W_i$ given by Proposition \ref{RADK616RVY} at a transverse contact and by Proposition \ref{RADK318JCA} at a tangential one. Let $\Lc_{p_i}^{\mathrm{iso}}$ be the isotropy group of $\Lc_{W_i}$ at $p_i$. It is the image under $q_j$ of $\{(a,(a-1)\gamma'(0),0,0):a>0\}\cong\R_{>0}$ at a transverse contact, and of $\{(1,b,0,0):b\in\R\}\cong\R$ at a tangential one. Figure \ref{RADK705FGE} shows the pieces of $\Lc$ near $p_i$.

The following lemma shows that, near each point of $\Lc_{p_i}^{\mathrm{iso}}$, the union of these pieces agrees with the smooth local model $\Lc_{W_i}$. The neighbourhood may depend on the isotropy arrow.

\begin{lemma}\label{RADK719YRS}
Let $K_i^-,K_i^+\in\K$, not necessarily distinct, be the arcs containing the two components of $(N_c\cap W_i)\setminus\{p_i\}$, and set
\[
A_i=(K_i^-\times K_i^-)\cup(K_i^+\times K_i^+)\cup\Lc_{p_i}^{\mathrm{iso}}\subset\G_\Sigma.
\]
\begin{enumerate}
\item If $p_i$ is a transverse contact, or if $K_i^-\neq K_i^+$, then $A_i\cap\overline\G_j|_{W_i}=\Lc_{W_i}$.
\item In general, every $g_0\in\Lc_{p_i}^{\mathrm{iso}}$ has an open neighbourhood $O_{g_0}\subset\overline\G_j|_{W_i}$ with $A_i\cap O_{g_0}=\Lc_{W_i}\cap O_{g_0}$.
\end{enumerate}
\end{lemma}

\begin{proof}
Let $g\in\Lc_{W_i}$. Its endpoints lie in $N_c\cap W_i$, so $g\in\overline\G_j|_{W_i}$. Either $\sss(g)=\ttt(g)=p_i$ and $g\in\Lc_{p_i}^{\mathrm{iso}}$, or $\sss(g)$ and $\ttt(g)$ lie on the same side of $p_i$ by \eqref{RADK654KCY} and \eqref{RADK821EHP}, hence in the same arc $K_i^\pm$. So $\Lc_{W_i}\subset A_i\cap\overline\G_j|_{W_i}$.

(1) Let $g\in(A_i\cap\overline\G_j|_{W_i})\setminus\Lc_{p_i}^{\mathrm{iso}}$. Then $g\in K_i^\pm\times K_i^\pm$ is determined by its endpoints in $(N_c\cap W_i)\setminus\{p_i\}$, and, by \eqref{RADK654KCY} and \eqref{RADK821EHP}, it suffices to show that they lie on the same side of $p_i$. At a transverse contact, a preimage $(a,b',s,r)\in O_{W_i}$ of $g$, in the coordinates of Lemma \ref{RADK450RVY}, has $s=0$ since $N_c\cap W_i=\{s=0\}$, so $\sss(g)=(\gamma(r),r)$ and $\ttt(g)=(\gamma(ar),ar)$, with $r\neq0$ and $\sgn(ar)=\sgn(r)$. At a tangential contact with $K_i^-\neq K_i^+$, the sets $K_i^\pm\cap W_i$ are $\{(t,h(t)):t<0\}$ and $\{(t,h(t)):t>0\}$, so $\sss(g)=(t,h(t))$ and $\ttt(g)=(t',h(t'))$ satisfy $tt'>0$.

(2) By (1), we may assume that $p_i$ is a tangential contact with $K_i^-=K_i^+=K$. By Lemma \ref{RADK788LQC}, $q_j$ identifies $\G^{\mathrm{loc}}_{j,W_i}$ with $\overline\G_j|_{W_i}$. Write $g_0=q_j(\hat g_0)$ with $\hat g_0=(1,\bar b,0,0)\in\G^{\mathrm{loc}}_{j,W_i}$. Since $h(t)/t\to h'(0)=0$ as $t\to0$ and $h(t)\neq0$ for $t\neq0$, there is $\varepsilon_1>0$ with
\begin{equation}\label{RADK814URT}
\frac{|t|}{|h(t)|}>|\bar b|+1\quad\text{for } t\in J,\ 0<|t|<\varepsilon_1 .
\end{equation}
Since $q_j$ is open, $O_{g_0}=q_j(\{(a,b,t,r)\in\G^{\mathrm{loc}}_{j,W_i}:|t|<\varepsilon_1,\ |b-\bar b|<1\})$ is an open neighbourhood of $g_0$ in $\overline\G_j|_{W_i}$.

Let $g\in A_i\cap O_{g_0}$, with preimage $(a,b,t,r)$ in this set. Since $\sss(g)\in N_c\cap W_i$, $r=h(t)$. If $t=0$, then $\sss(g)=\ttt(g)=p_i$ by \eqref{RADK472KHB}, so $g\in\Lc_{p_i}^{\mathrm{iso}}$ since $K$ misses $D$. If $t\neq0$, then $g\in K\times K$ and $\ttt(g)=(t',h(t'))$ with $t'=t+b\,h(t)\in\widetilde J\setminus\{0\}$. If $tt'<0$, then $|t'-t|>|t|$, so \eqref{RADK814URT} gives $|b|=|t'-t|/|h(t)|>|t|/|h(t)|>|\bar b|+1$, hence $|b-\bar b|>1$, a contradiction. Thus $t'\in J_t$, and $g\in\Lc_{W_i}$ by \eqref{RADK821EHP}.
\end{proof}

\begin{remark}\label{RADK577BWY}
If $\mathcal C_c=\{p_i\}$, then the arc $K=N_c\setminus\{p_i\}$ lies in one component of $\Sigma\setminus D$, so $K\cap W_i$ lies on one side of $D_j$ by \eqref{RADK247VWC}. Hence the contact at $p_i$ is tangential, since $N_c\cap W_i$ meets both sides at a transverse contact, and $h(t)h(t')>0$ for $t,t'\in J\setminus\{0\}$. For $tt'<0$, the image of $\big(h(t')/h(t),(t'-t)/h(t),t,h(t)\big)\in\G_j|_{W_i}$, an arrow from $(t,h(t))$ to $(t',h(t'))$, lies in $(K\times K)\cap\overline\G_j|_{W_i}$ but not in $\Lc_{W_i}$ by \eqref{RADK821EHP}, so the equality of Lemma \ref{RADK719YRS}(1) fails. For $h(t)=t^2$ and $t'=-t$, its preimages $\big(1,-2/t+k\ell_j/t^2,t,t^2\big)$, $k\in\Z$, satisfy $|b|\ge2/|t|$ for $0<|t|\le\ell_j/4$, so these arrows have no accumulation point in $\G_\Sigma$ as $t\to0$: a local inverse of $q_j$ near such a point would give a convergent choice of preimages, contradicting the displayed bound. Thus, in general, the comparison in Lemma \ref{RADK719YRS}(2) must be made near the isotropy arrows rather than on the whole restriction to $W_i$.
\end{remark}

\begin{proposition}\label{RADK663CKE}
The subset of $\G_\Sigma$
\[
\begin{aligned}
\Lc={}&\Big(\coprod_{N_c\subset\Sigma\setminus D}N_c\times N_c\Big)
\sqcup\Big(\coprod_{N_c\subset D}\Lc_c\Big)
\sqcup\Big(\coprod_{K\in\K}K\times K\Big)
\sqcup\Big(\coprod_{i=1}^{n_{\mathcal C}}\Lc_{p_i}^{\mathrm{iso}}\Big)
\end{aligned}
\]
is an embedded, source-connected Lagrangian subgroupoid of $\G_\Sigma$ over $N$ integrating $TN^\circ$.
\end{proposition}

\begin{proof}
Each piece of $\Lc$ consists of arrows with both endpoints in one component of $N_0$ or $N_D$, in one arc $K\in\K$, or at one point of $\mathcal C$, and distinct pieces correspond to distinct such subsets of \eqref{RADK376JAG}. So two composable arrows of $\Lc$ lie in the same piece, and $\Lc$ is a subgroupoid over $N$ since each piece is one (Proposition \ref{RADK001UFG} for $\Lc_c$).

Let $g_0\in\Lc$. If an endpoint of $g_0$ is some $p_i\in\mathcal C$, then $g_0\in\Lc_{p_i}^{\mathrm{iso}}$ by \eqref{RADK472KHB}, applied to $g_0$ or $g_0^{-1}$. Otherwise, since $D$ is closed and $\mathcal C$ is finite, each endpoint $x$ has a neighbourhood $U_x$ in $\Sigma$ with $U_x\cap N$ contained in a single component of $N_0$ or $N_D$, or in a single arc of $\K$. On $(\sss,\ttt)^{-1}(U_{\sss(g_0)}\times U_{\ttt(g_0)})$, the subgroupoid $\Lc$ therefore agrees with the restriction of $N_c\times N_c$, $\Lc_c$ or $K\times K$.

If $g_0\in\Lc_{p_i}^{\mathrm{iso}}$, then $\Lc\cap\overline\G_j|_{W_i}=A_i\cap\overline\G_j|_{W_i}$ since $N\cap W_i=N_c\cap W_i\subset K_i^-\cup K_i^+\cup\{p_i\}$, so $\Lc\cap O_{g_0}=\Lc_{W_i}\cap O_{g_0}$ by Lemma \ref{RADK719YRS}(2).

So $\Lc$ coincides near each of its points with one of $N_c\times N_c$, $\Lc_c$, $K\times K$ or $\Lc_{W_i}$, which are embedded Lagrangian Lie subgroupoids integrating $TN^\circ$ over open subsets of $N$ (\S\ref{RADK352LLS} and Propositions \ref{RADK001UFG}, \ref{RADK616RVY} and \ref{RADK318JCA}). Hence $\Lc$ is embedded, $\sss|_\Lc$ is a submersion, and the source fibres $N_c$, $\R$, $K$ and $\Lc_{p_i}^{\mathrm{iso}}$ are connected.
\end{proof}

It remains to identify the isotropy groups at $x\in N\cap D$ inside $\G_{\Sigma,x}\cong\Aff$. By \eqref{RADK795GLQ}, a dilation subgroup has the form
\begin{equation}\label{RADK526KVQ}
\{(a,(a-1)b_0):a>0\}=(1,b_0)\,\{(a,0):a>0\}\,(1,b_0)^{-1}
\end{equation}
and its Lie algebra is $\langle\zeta_1+b_0\,\zeta_2\rangle$. As $b_0$ varies, these are all lines of $\aff$ other than $\langle\zeta_2\rangle$, the Lie algebra of the translation subgroup. Thus every connected one-dimensional Lie subgroup of $\Aff$ is a translation or dilation subgroup. Automorphisms of $\Aff$ preserve these two types, since its translation subgroup is its commutator subgroup.

\begin{proof}[\textnormal{\bfseries Proof of Theorem \ref{RADK953BQH}}]
The subgroupoid $\Lc$ of Proposition \ref{RADK663CKE} is a $1$-shifted Lagrangian by the argument of \cite[Example 3.8]{Maxence2025}, since $\Omega_\Sigma|_\Lc=0$ and $\Lie(\Lc)=TN^\circ$.

The arrows of $\Lc$ between points of $N\setminus D$ are those of the pair groupoids of the components of $N\setminus D$, namely the $N_c\subset N_0$ and the $K\in\K$. Arrows with source in $N\cap D$ are loops by \eqref{RADK472KHB}.

Let $x\in N\cap D_j$. If $D_j\subset N$ or $x=p_i$ is a tangential contact, then $T_xN=T_xD$ and $\Lc_x$ is the translation subgroup, as the image of $\{(1,b,t,0):b\in\R\}$ by Proposition \ref{RADK001UFG}, respectively of $\{(1,b,0,0):b\in\R\}$. At a transverse contact $x=p_i$, $T_xN\neq T_xD$ and $\Lc_x$ is the image of $\{(a,(a-1)\gamma'(0),0,0):a>0\}$, a dilation subgroup by \eqref{RADK526KVQ} with $b_0=\gamma'(0)$.
\end{proof}

\section{Homotopy fibre products of Lagrangian subgroupoids}\label{RADK389HFP}

In this section, we study the shifted intersections of two admissible curves. Let $N_1,N_2$ be admissible curves, with the embedded Lagrangian subgroupoids
$\Lc_1\tto N_1$ and $\Lc_2\tto N_2$ constructed in Proposition \ref{RADK663CKE}.
We determine when their homotopy fibre product is clean, proving Theorem \ref{RADK367ESC},
and describe its induced $0$-shifted Poisson structure up to Morita equivalence.
The cleanness criterion depends only on the contacts of the two curves with $D$.
When both curves contain a component of $D$, the resulting shifted intersection
also recovers its modular period.

\subsection{Cleanness}\label{RADK411CLN}

We first describe the homotopy fibre product and the tangent-space condition defining cleanness. Let
$\Hc=\Lc_1\times^h_{\G_\Sigma}\Lc_2$ be the homotopy fibre product of $\Lc_1$ and $\Lc_2$
\cite[\S4.2]{Maxence2025}. It is a groupoid whose space of objects is
\[
\Hc^0=\big\{(x_1,g,x_2)\in N_1\times\G_\Sigma\times N_2:\sss(g)=x_1,\ \ttt(g)=x_2\big\}.
\]
An arrow with source $(x_1,g,x_2)$ is a pair
$(l,\tilde l)\in\Lc_1\times\Lc_2$ satisfying
$\sss(l)=x_1$ and $\sss(\tilde l)=x_2$, and its target is $\bigl(\ttt(l),\,\tilde l\,g\,l^{-1},\,\ttt(\tilde l)\bigr)$.

We write $\pr_1$, $\pr_0$ and $\pr_2$ for the projections of $\Hc^0$ onto $N_1$, $\G_\Sigma$ and $N_2$. The
\defn{expected tangent space} at $(x_1,g,x_2)$ is
\[
T^{\exp}_{(x_1,g,x_2)}
=\bigl\{(u,w,v)\in T_{x_1}N_1\times T_g\G_\Sigma\times T_{x_2}N_2:
d\sss_g(w)=u,\ d\ttt_g(w)=v\bigr\}.
\]
The homotopy fibre product $\Hc$ is \defn{clean} if $\Hc^0$ is an embedded submanifold of
$N_1\times\G_\Sigma\times N_2$ and its tangent space at every object equals the
expected tangent space \cite[\S\S2.3 and 4.2]{Maxence2025}.
When the object is fixed, we abbreviate this space to $T^{\exp}$.

The map $g\mapsto(\sss(g),g,\ttt(g))$ embeds $\G_\Sigma$ onto a closed submanifold of $\Sigma\times\G_\Sigma\times\Sigma$
and carries $(\sss,\ttt)^{-1}(N_1\times N_2)$ onto $\Hc^0$. We therefore identify $\Hc^0$ with
$(\sss,\ttt)^{-1}(N_1\times N_2)\subset\G_\Sigma$. Under this identification,
\[
T^{\exp}
=\bigl\{w\in T_g\G_\Sigma:
d\sss_g(w)\in T_{x_1}N_1,\ d\ttt_g(w)\in T_{x_2}N_2\bigr\}.
\]

For a $2$-form $\omega$, let $\Gamma_\omega$ be its graph, as in \S\ref{RADK219CNV}. The $1$-shifted Lagrangians $\Lc_1$ and $\Lc_2$ have vanishing $2$-forms, and their $1$-shifted coisotropic structures are
$\mathbb L_{N_1}=\Gamma_0=TN_1$ and $\mathbb L_{N_2}=TN_2$ \cite[\S3.2, Example 3.8 and Remark 4.9]{Maxence2025}.
For graphs, the pullback and sum of Dirac structures satisfy
\[
f^*\Gamma_\omega=\Gamma_{f^*\omega},
\quad
\Gamma_\omega+\Gamma_\eta=\Gamma_{\omega+\eta}.
\]
If $\Hc$ is clean, the Dirac structure defining the induced shifted intersection is
\begin{align*}
\mathbb L
&=-\pr_1^*\mathbb L_{N_1}
  +\pr_2^*\mathbb L_{N_2}
  -\pr_0^*\Gamma_{\Omega_\Sigma}\\
&=\Gamma_{-\pr_0^*\Omega_\Sigma}.
\end{align*}
It is smooth because it is the graph of a smooth $2$-form. Thus it defines a
$0$-shifted Poisson structure on $\Hc$ by \cite[Corollary 4.12]{Maxence2025}.

We separate the objects lying over the symplectic locus from those lying over $D$:
\begin{align*}
\Hc^0_{\mathrm{gen}}
&=\bigl\{(x_1,g,x_2)\in\Hc^0:x_1\notin D\bigr\},\\
\Hc^0_{\mathrm{iso}}
&=\bigl\{(x,g,x)\in\Hc^0:x\in D\cap N_1\cap N_2\bigr\}.
\end{align*}
By \eqref{RADK472KHB}, applied to $g$ and to
$g^{-1}$, we have $x_1\in D$ if and only if $x_2\in D$, and then $x_2=x_1$. Hence
$\Hc^0=\Hc^0_{\mathrm{gen}}\sqcup\Hc^0_{\mathrm{iso}}$. By \eqref{RADK618ENP}, if $x_1\in V_i$, then $x_2\in V_i$.
So $\Hc^0_{\mathrm{gen}}$ is the disjoint union of the open subsets
$\Hc^0_i=\Hc^0_{\mathrm{gen}}\cap(V_i\times\G_\Sigma\times V_i)$.

Since $(\sss,\ttt)$ is a diffeomorphism from $\G_\Sigma|_{V_i}=V_i\times V_i$ onto $V_i\times V_i$, $\Hc^0_i=(\sss,\ttt)^{-1}\big((N_1\cap V_i)\times(N_2\cap V_i)\big)$ is embedded, with the
expected tangent space $d(\sss,\ttt)_g^{-1}\big(T_{x_1}N_1\times T_{x_2}N_2\big)$, and
\begin{equation}\label{RADK917CNZ}
(N_1\cap V_i)\times(N_2\cap V_i)\too\Hc^0_i,\quad (x_1,x_2)\longmapsto\big(x_1,(x_2,x_1),x_2\big),
\end{equation}
is a diffeomorphism. Therefore $\Hc$ is clean on $\Hc^0_{\mathrm{gen}}$.
Moreover, Proposition \ref{RADK181QST} gives
$\Omega_\Sigma=\ttt^*\omega_i-\sss^*\omega_i$ on $V_i\times V_i$.
Since $\sss\circ\pr_0=\pr_1$ and $\ttt\circ\pr_0=\pr_2$, we obtain
\[
\pr_0^*\Omega_\Sigma =\pr_2^*\bigl(\omega_i|_{N_2\cap V_i}\bigr) -\pr_1^*\bigl(\omega_i|_{N_1\cap V_i}\bigr) =0
\]
on $\Hc^0_i$: both restrictions are $2$-forms on curves.

\subsection{The six local configurations}\label{RADK422SLC}

It remains to study $\Hc$ near $\Hc^0_{\mathrm{iso}}$.
Fix $x\in D_j\cap N_1\cap N_2$. After translating the coordinate $t$ in
\eqref{RADK871UDA}, we take $x=(0,0)$ and write an element of
$\G_{\Sigma,x}\cong\Aff$ as $g=(a,b,0,0)$. We retain the coordinate
$\alpha=\log a$ introduced in \S\ref{RADK265ACT}.
We first compute the expected tangent space, before comparing it with the local equations of $\Hc^0$.

\begin{lemma}\label{RADK308NNC}
For $g=(a,b,0,0)\in\G_{\Sigma,x}$, let $M_g=\left(\begin{smallmatrix}1&b\\0&a\end{smallmatrix}\right)$ in the frame
$(\dd t,\dd r)$ of $T_x\Sigma$. Then $d\ttt_g=M_g\circ d\sss_g$, the kernels $\ker d\sss_g=\ker d\ttt_g$ have
dimension $2$, and the expected tangent space at $(x,g,x)$ satisfies
\begin{equation}\label{RADK505NHY}
\dim T^{\exp}=2+\dim\big(T_xN_1\cap M_g^{-1}T_xN_2\big).
\end{equation}
\begin{enumerate}
\item If $T_xN_1=T_xN_2=T_xD_j$, then $\dim T^{\exp}=3$ for every $g$.
\item If exactly one of $T_xN_1$ and $T_xN_2$ equals $T_xD_j$, then $\dim T^{\exp}=2$ for every $g$.
\item If both $N_1$ and $N_2$ are transverse to $D$ at $x$, with the functions $\gamma$ and $\tilde\gamma$ of
Lemma \ref{RADK227HCV}, then $\dim T^{\exp}=2$, except on the curve $\{b=a\,\tilde\gamma'(0)-\gamma'(0)\}$ of
$\G_{\Sigma,x}$, where it equals $3$.
\end{enumerate}
\end{lemma}

\begin{proof}
Since $\ttt(a,b,t,r)=(t+br,ar)$ and $r=0$ at $g$, a tangent vector
$w=(\dot a,\dot b,\dot t,\dot r)$ satisfies
\[
d\sss_g(w)=(\dot t,\dot r),\quad d\ttt_g(w)=(\dot t+b\,\dot r,\ a\,\dot r) =M_g\bigl(\dot t,\dot r\bigr).
\]
Consequently,
\[
\ker d\sss_g=\ker d\ttt_g
=\bigl\langle\dd a,\dd b\bigr\rangle.
\]
The surjection $d\sss_g$ maps $T^{\exp}$ onto
$\{u\in T_xN_1:M_gu\in T_xN_2\}$, with kernel $\ker d\sss_g$. This proves \eqref{RADK505NHY}.

Since $M_g\dd t=\dd t$, the map $M_g$ preserves
$T_xD_j=\langle\dd t\rangle$, proving (1).
In case (2), suppose first that $T_xN_1=T_xD_j\neq T_xN_2$.
Then $M_g^{-1}T_xN_2\neq T_xD_j$, so the intersection in
\eqref{RADK505NHY} is zero. The other order is treated in the same way.
In case (3), the tangent lines are
\begin{align*}
T_xN_1&=\bigl\langle\gamma'(0)\dd t+\dd r\bigr\rangle,\\
T_xN_2&=\bigl\langle\tilde\gamma'(0)\dd t+\dd r\bigr\rangle.
\end{align*}
The vector
\[
M_g\bigl(\gamma'(0)\dd t+\dd r\bigr)
=\bigl(\gamma'(0)+b\bigr)\dd t+a\dd r
\]
lies in $T_xN_2$ exactly when $\gamma'(0)+b=a\,\tilde\gamma'(0)$.
\end{proof}

\begin{figure}[ht]
\centering
\includegraphics[width=0.9\linewidth]{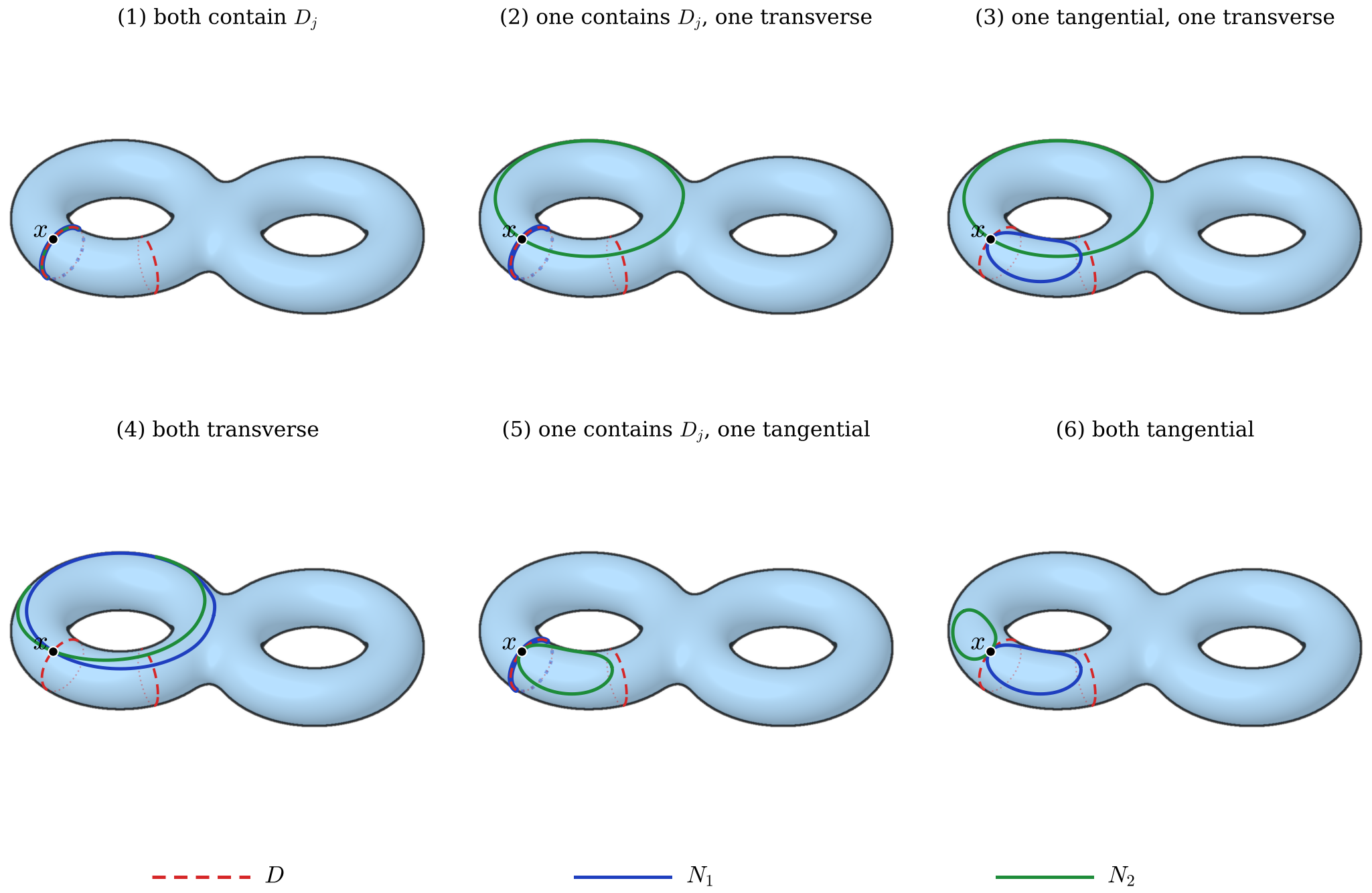}
\caption{The six configurations of two curves at a common point of $D$.}\label{RADK706FGF}
\end{figure}

We now compare these expected tangent spaces with the local equations of $\Hc^0$.
By Lemmas \ref{RADK227HCV} and \ref{RADK562YPM}, there are three possibilities for each curve near $x$:
it contains $D_j$, it is transverse to $D$ and is given by $t=\gamma(r)$,
or it has a tangential contact with $D$ and is given by $r=h(t)$.
For $N_2$, we use the functions $\tilde\gamma$ and $\tilde h$.
Up to interchanging the curves, this gives the six configurations of Figure \ref{RADK706FGF}.

We fix one box $W=J\times I$ around $x$ which only meets the components of $N_1$ and $N_2$ through $x$, and in which
both descriptions hold. For two transverse curves we shrink $I$, and for two tangential curves we shrink $J$. If one
curve contains $D_j$, we shrink $I$ if the other is transverse and $J$ if it is tangential. For a
transverse curve and a tangential one, we take $I=(-\varepsilon,\varepsilon)$ and
$J=(-C\varepsilon,C\varepsilon)$, with $C>\sup|\gamma'|$ fixed and $\varepsilon$ small. This is
possible since $|\gamma(r)|\le\sup|\gamma'|\,|r|$ and $h(t)=O(t^2)$.

We use the open set $\G^{\mathrm{loc}}_{j,W}$ of Lemma \ref{RADK788LQC}. It contains the isotropy group $\G_{j,x}$, and $q_j$ identifies it with $\overline\G_j|_W$, an open subset of $\G_\Sigma$. Under this identification,
\[
\Hc^0\cap\G^{\mathrm{loc}}_{j,W}=\big\{(a,b,t,r)\in\G^{\mathrm{loc}}_{j,W}:(t,r)\in N_1,\ (t+br,ar)\in N_2\big\}.
\]

\begin{proposition}\label{RADK638FNS}
Let $x\in D_j\cap N_1\cap N_2$, $g_0\in\G_{\Sigma,x}$, and let $\hat g_0\in\G^{\mathrm{loc}}_{j,W}$ be its preimage. There is an open neighbourhood $O\subset\G^{\mathrm{loc}}_{j,W}$ of $\hat g_0$ on which $\Hc^0$ is described as follows, in the
six configurations of Figure \ref{RADK706FGF}.
\begin{enumerate}
\item If $D_j\subset N_1$ and $D_j\subset N_2$, then $\Hc^0\cap O=O\cap\{r=0\}$,
an open subset of the embedded $3$-manifold $D_j\times\Aff$, and $\Hc$ is clean there.
\item If one of $N_1,N_2$ contains $D_j$ and the other is transverse, then $\Hc^0\cap O=O\cap
\{t=r=0\}$, and $\Hc$ is clean there.
\item If one of $N_1,N_2$ is tangential and the other is transverse, then $\Hc^0\cap O=O\cap
\{t=r=0\}$, and $\Hc$ is clean there.
\item If both $N_1$ and $N_2$ are transverse, write $\gamma(r)=r\,\gamma_1(r)$ and
$\tilde\gamma(r)=r\,\tilde\gamma_1(r)$. Then
\[
\Hc^0\cap O=O\cap\Big(\{t=r=0\}\ \cup\ \big\{t=\gamma(r),\ b=a\,\tilde\gamma_1(ar)-\gamma_1(r)\big\}\Big).
\]
These two embedded surfaces meet along the curve $\{t=r=0,\ b=a\,\tilde\gamma'(0)-\gamma'(0)\}$, with distinct
tangent planes, and $\Hc^0$ is not a manifold at the points of this curve.
\item If one of $N_1,N_2$ contains $D_j$ and the other is tangential, then $\Hc^0\cap O=O\cap
\{t=r=0\}$ is a surface, but $\dim T^{\exp}=3$. So $\Hc$ is not clean there.
\item If both $N_1$ and $N_2$ are tangential, then $\Hc^0\cap O$ contains $O\cap\{t=r=0\}$, and the expected tangent space has dimension $3$ at every point of this isotropy surface. However, $\Hc^0\cap O$ has empty interior in the $3$-manifold $\{r=h(t)\}$. Thus $\Hc$ is not clean at these isotropy objects. For instance, if $h(t)=\tilde h(t)=t^2$, then $\Hc^0\cap O$ is the union of the sheets $\{t=0\}$ and $\{a=(1+bt)^2\}$ of $\{r=t^2\}$, which meet along $\{t=0,\ a=1\}$.
\end{enumerate}
\end{proposition}

\begin{proof}
The inversion $(x_1,g,x_2)\mapsto(x_2,g^{-1},x_1)$ identifies $\Lc_1\times^h_{\G_\Sigma}\Lc_2$ with $\Lc_2\times^h_{\G_\Sigma}\Lc_1$ and preserves cleanness and $\{t=r=0\}$, so we may fix the order of the curves in cases (2), (3) and (5). Let $O=\G^{\mathrm{loc}}_{j,W}$, except in case (3).

In case (1), $(t,r)\in N_1$ gives $r=0$, and then $(t+br,ar)=(t,0)\in N_2$. Then $T\{r=0\}\subset T^{\exp}$, with equality since $\dim T^{\exp}=3$ by Lemma \ref{RADK308NNC}(1).

In case (2), say $N_1\supset D_j$ and $N_2$ is transverse. Then $r=0$, and $(t,0)\in N_2$ exactly when $t=\tilde\gamma(0)=0$. So $T\{t=r=0\}=T^{\exp}$ by Lemma \ref{RADK308NNC}(2).

In case (3), suppose that $N_1$ is tangential and $N_2$ is transverse. The defining equations are $r=h(t)$, $0=t+b\,h(t)-\tilde\gamma\bigl(a\,h(t)\bigr)$. The second right-hand side vanishes identically at $t=0$ and has derivative $1$ with respect to $t$ there, since $h'(0)=0$. The implicit function theorem therefore gives a neighbourhood $O$ of $\hat g_0$ on which its only zeros have $t=0$. Then $r=0$, and the tangent-space comparison is the same as in case (2).

In case (4), the source condition is $t=\gamma(r)$, and the target condition factors as
\[
0 =\tilde\gamma(ar)-\gamma(r)-br =r\bigl(a\,\tilde\gamma_1(ar)-\gamma_1(r)-b\bigr).
\]
The factor $r$ gives the surface $\{t=r=0\}$, parametrized by $(a,b)$, with tangent plane $\langle\dd a,\dd b\rangle$.
The other factor gives the graph
\[
\bigl\{t=\gamma(r),\ b=a\,\tilde\gamma_1(ar)-\gamma_1(r)\bigr\},
\]
parametrized by $(a,r)$, whose tangent plane contains a vector with $\dd r$-component $1$. Their intersection is the curve $\bigl\{t=r=0,\ b=a\,\tilde\gamma'(0)-\gamma'(0)\bigr\}$, since $\gamma_1(0)=\gamma'(0)$ and
$\tilde\gamma_1(0)=\tilde\gamma'(0)$. If $\Hc^0$ were a submanifold near a point of this curve, its tangent space would contain both tangent planes and have dimension at least $3$. This is impossible for a union of two surfaces.

In case (5), say $N_1\supset D_j$ and $N_2$ is tangential. Then $r=0$ and $\tilde h(t)=0$, that is, $t=0$, while $\dim T^{\exp}=3$ by Lemma \ref{RADK308NNC}(1).

In case (6), the equations are $r=h(t)$, $a\,h(t)=\tilde h\bigl(t+b\,h(t)\bigr)$. They hold on $\{t=r=0\}$, where $\dim T^{\exp}=3$ by Lemma \ref{RADK308NNC}(1).
On the $3$-manifold $\{r=h(t)\}$, parametrized by $(a,b,t)$, set $F(a,b,t)=a\,h(t)-\tilde h\bigl(t+b\,h(t)\bigr)$. Then $\partial F/\partial a=h(t)\neq0$ for $t\neq0$, because the contact is isolated. Thus the zero set has empty interior away from $t=0$, and it cannot contain an open subset supported on $\{t=0\}$ either. Hence $\Hc^0\cap O$ has empty interior in $\{r=h(t)\}$. If $\Hc$ were clean at $(x,g_0,x)$, then $\Hc^0$ would be a $3$-dimensional embedded submanifold of this $3$-manifold near $\hat g_0$, and would therefore contain a neighbourhood of $\hat g_0$, a contradiction. For $h=\tilde h=t^2$, the factorization is $a\,t^2-\bigl(t+b\,t^2\bigr)^2 =t^2\bigl(a-(1+b\,t)^2\bigr)$.
\end{proof}

We can now decide when $\Hc$ is clean.

\begin{proposition}\label{RADK783SZT}
The homotopy fibre product $\Hc=\Lc_1\times^h_{\G_\Sigma}\Lc_2$ is clean if and only if, for every $x\in D\cap N_1\cap N_2$, with $D_j$ the component of $D$ containing $x$,
\begin{enumerate}
\item either $D_j\subset N_1$ and $D_j\subset N_2$,
\item or exactly one of the lines $T_xN_1$ and $T_xN_2$ equals $T_xD_j$.
\end{enumerate}
In that case, $\Gamma_{-\pr_0^*\Omega_\Sigma}$ is a $0$-shifted Poisson structure on $\Hc$, with
\[
\pr_0^*\Omega_\Sigma=
\begin{cases}
0 & \text{on }\Hc^0_{\mathrm{gen}}\text{ and at the points where (2) holds},\\
dt\wedge d\alpha & \text{on the components }D_j\times\Aff\text{ where (1) holds}.
\end{cases}
\]
\end{proposition}

\begin{proof}
Cleanness is local, and $\Hc$ is clean on the open set $\Hc^0_{\mathrm{gen}}$ by \S\ref{RADK411CLN}. At $x\in D\cap N_1\cap N_2$, one of the six cases of Proposition \ref{RADK638FNS} occurs. In cases (1)--(3), $\Hc$ is clean near $\{x\}\times\G_{\Sigma,x}\times\{x\}$, in cases (5) and (6) at no point of this set, and in case (4) not on the nonempty curve $\{b=a\,\tilde\gamma'(0)-\gamma'(0)\}$. Since $T_xN_1=T_xD_j$ exactly when $N_1$ contains $D_j$ or has a tangential contact with $D$ at $x$, and similarly for $N_2$, cases (1)--(3) of Proposition \ref{RADK638FNS} are those allowed by conditions (1) and (2).

The form $\pr_0^*\Omega_\Sigma$ vanishes on $\Hc^0_{\mathrm{gen}}$ by \S\ref{RADK411CLN}. In cases (1)--(3) of Proposition \ref{RADK638FNS}, $\Hc^0$ lies in $\{r=0\}$ near $\{x\}\times\G_{\Sigma,x}\times\{x\}$, and in $\{t=r=0\}$ in cases (2) and (3). By \eqref{RADK958NSM}, $\Omega|_{\{r=0\}}=dt\wedge d\alpha$ and $\Omega|_{\{t=r=0\}}=0$.

If $\Hc$ is clean, then $\Gamma_{-\pr_0^*\Omega_\Sigma}$ is the graph of a smooth $2$-form, hence a $0$-shifted Poisson structure by \cite[Corollary 4.12(1)]{Maxence2025}.
\end{proof}

For the shifted self-intersection of a single curve, the criterion takes a simple form.

\begin{corollary}\label{RADK784SLF}
Let $N$ be an admissible curve with Lagrangian subgroupoid $\Lc$. Then $\Lc\times^h_{\G_\Sigma}\Lc$ is clean if and only if every component of $N$ which meets $D$ is a component of $D$.
\end{corollary}

\begin{proof}
At $x\in D\cap N$, condition (2) of Proposition \ref{RADK783SZT} fails since the two lines coincide. So $\Lc\times^h_{\G_\Sigma}\Lc$ is clean if and only if $D_j\subset N$ for the component $D_j$ of $D$ through every $x\in D\cap N$. Since $D_j$ and the components of $N$ are circles, $D_j\subset N$ means that $D_j$ is the component of $N$ through $x$.
\end{proof}

\subsection{Morita equivalence and modular periods}\label{RADK433MMP}

We now describe the clean homotopy fibre products up to Morita equivalence and identify the function that recovers the modular period. On a component $D_j\times\Aff$ of $\Hc^0$, we use the coordinates $(t,a,b)$ of the arrow $(a,b,t,0)$ and write $\alpha=\log a$. The function $a$ is intrinsic: the map $M_g$ of Lemma \ref{RADK308NNC}, characterized by $d\ttt_g=M_g\circ d\sss_g$, induces multiplication by $a$ on the normal line $T_x\Sigma/T_xD_j$.

\begin{proposition}\label{RADK741QHM}
Suppose that $\Hc$ is clean. Then $\Hc$ is Morita equivalent to the disjoint union of finitely many points and, for each component $D_j$ of $D$ contained in $N_1\cap N_2$, of the trivial bundle of groups $(S^1_{\ell_j}\times\R)\times\R$
over the cylinder $S^1_{\ell_j}\times\R$. Under this Morita equivalence, the $0$-shifted Poisson structure of
Proposition \ref{RADK783SZT} restricts to the graph of the zero form on the points and of the symplectic form $d\alpha\wedge dt$
on each cylinder. Here the last factor $\R$ is the additive isotropy group: the source and target maps are both the projection onto the cylinder.
\end{proposition}

\begin{proof}
We decompose $\Hc^0$ into subsets of three types:
\begin{enumerate}
\item the sets
\[
\Hc^0_{K_1,K_2}
=(\sss,\ttt)^{-1}(K_1\times K_2)\cong K_1\times K_2,
\]
where $K_1$ and $K_2$ are components of $N_1\setminus D$ and $N_2\setminus D$ lying in the same component of $\Sigma\setminus D$;
\item the sets $\{x\}\times\G_{\Sigma,x}\times\{x\}\cong\Aff$, where $x$ is an isolated point of $D\cap N_1\cap N_2$;
\item the sets $D_j\times\Aff$, where $D_j\subset N_1\cap N_2$ is a component of $D$.
\end{enumerate}

An arrow $(l,\tilde l)$ of $\Hc$ moves $x_1$ to $\ttt(l)$ and $x_2$ to $\ttt(\tilde l)$, so, by Proposition \ref{RADK663CKE} and \eqref{RADK472KHB}, it preserves the component of $N_1\setminus D$ or the point of $N_1\cap D$ containing $x_1$, and similarly for $x_2$. If $x_1\in D$, then $x_1=x_2$ by \eqref{RADK472KHB}, and condition (1) or (2) of Proposition \ref{RADK783SZT} holds at $x_1$. Under (2), one curve is transverse to $D$ at $x_1$, so $x_1$ is isolated in $D\cap N_1\cap N_2$. So the sets of types (1)--(3) are finitely many, invariant, and partition $\Hc^0$. Those of type (1) or (3) are open as preimages under $(\pr_1,\pr_2)$ of open subsets of $N_1\times N_2$, and those of type (2) by cases (2) and (3) of Proposition \ref{RADK638FNS}, the only ones allowed at an isolated point by Proposition \ref{RADK783SZT}.

On a set of type (1), $\Lc_1$ and $\Lc_2$ restrict to the pair groupoids of $K_1$ and $K_2$ (Proposition \ref{RADK663CKE}), so, by \eqref{RADK917CNZ}, $\Hc$ restricts to the pair groupoid of $K_1\times K_2$, Morita equivalent to a point, and the form vanishes by Proposition \ref{RADK783SZT}.

On a set of type (2), $\Hc$ restricts to the action groupoid of $\Lc_{1,x}\times\Lc_{2,x}$ on $\Aff$, with action $(l,\tilde l)\cdot g=\tilde l\,g\,l^{-1}$.
By Theorem \ref{RADK953BQH}, one isotropy subgroup is the translation subgroup, the kernel of $(a,b)\mapsto a$, and the other is a dilation subgroup, which this map identifies with $\R_{>0}$. The translation subgroup is normal, and every conjugate of a dilation subgroup is a dilation subgroup. Consequently,
\[
\Lc_{1,x}\cap g^{-1}\Lc_{2,x}g=\{1\},\quad \Lc_{2,x}\bigl(g\Lc_{1,x}g^{-1}\bigr)=\Aff.
\]
The first identity says that the stabilizer of $g$ is trivial, while the second gives
\[
\Lc_{2,x}\,g\,\Lc_{1,x} =\Lc_{2,x}\bigl(g\Lc_{1,x}g^{-1}\bigr)\,g =\Aff.
\]
The action is therefore free and transitive, so its action groupoid is Morita equivalent to a point. The form vanishes by Proposition \ref{RADK783SZT}.

On a set of type (3), $\Lc_1$ and $\Lc_2$ restrict to the bundle of translation groups over $D_j$ (Proposition \ref{RADK001UFG}). For $l=(1,u,t,0)$ and $\tilde l=(1,\tilde u,t,0)$, \eqref{RADK795GLQ} gives
\[
(1,\tilde u)\,(a,b)\,(1,u)^{-1}
=\bigl(a,\,b-u+a\,\tilde u\bigr).
\]
The orbits are the lines $\{(t,a)\}\times\R$, the stabilizers are $\{u=a\tilde u\}\cong\R$, and the slice $\{b=0\}$ meets every orbit once. The restriction to this slice is the trivial bundle of additive groups $\R$ over the cylinder $S^1_{\ell_j}\times\R$, with coordinates $(t,\alpha)$. Its inclusion is a Morita equivalence: the target map on arrows with source in the slice is $(t,a,\tilde u,u)\longmapsto\bigl(t,a,a\,\tilde u-u\bigr)$,
a surjective submersion onto $D_j\times\Aff$. By Proposition \ref{RADK783SZT}, $-\pr_0^*\Omega_\Sigma =-dt\wedge d\alpha =d\alpha\wedge dt$. This form does not involve $b$, and its restriction to the slice is the symplectic form $d\alpha\wedge dt$. Finally, $\Hc^0$ is the disjoint union of the finitely many open invariant sets of types (1)--(3), so $\Hc$ is the disjoint union of the corresponding restrictions. This proves the claimed Morita equivalence and the description of the induced structures.
\end{proof}

The canonical function $\alpha$ recovers the modular period. With the Hamiltonian convention of \S\ref{RADK219CNV},
\[
\iota_{-\dd t}\bigl(d\alpha\wedge dt\bigr)=d\alpha,
\quad X_\alpha=-\dd t.
\]
By \eqref{RADK432CKK}, $X_\alpha$ projects to the modular vector field along $D_j$. Its flow rotates $S^1_{\ell_j}$ at unit speed and has period $\ell_j$. The symplectic cylinder alone does not determine this period: for every $\ell>0$,
the map
\[
(t,\alpha)\longmapsto
\Bigl(\frac{\ell}{\ell_j}\,t,\frac{\ell_j}{\ell}\,\alpha\Bigr)
\]
is a symplectomorphism onto $S^1_\ell\times\R$.
The period is therefore recovered from the symplectic cylinder together with the canonical function $\alpha$.

For a clean self-intersection $\Lc\times^h_{\G_\Sigma}\Lc$ as in Corollary \ref{RADK784SLF}, $D\cap N$ has no isolated points,
the points of Proposition \ref{RADK741QHM} correspond to the ordered pairs of components of $N\setminus D$ in the same component
of $\Sigma\setminus D$, and there is one cylinder for each component of $D$ contained in $N$.

\begin{example}\label{RADK748QSA}
Let $S^2\subset\R^3$ be the unit sphere, with cylindrical coordinates $(\theta,z)$, $\theta\in\R/2\pi\Z$. Its area
form is $d\theta\wedge dz$, so $\dd\theta\wedge\dd z$ extends to a smooth bivector on $S^2$, and
$\sigma=z\,\dd z\wedge\dd\theta$ is a Radko structure, with $D$ the equator and $(t,r)=(\theta,z)$ as in
\eqref{RADK871UDA}. So the modular period of $D$ is $\ell_1=2\pi$. The meridian $N_1$, the great circle consisting of the poles and the points with $\theta\in\{\theta_0,\theta_0+\pi\}$, meets $D$ transversely at two points $p_1,p_2$. If $N_2\neq D$ is a great circle through $p_1$ and $p_2$, for instance $N_2=N_1$, then we are in case (4) of Proposition \ref{RADK638FNS} at
$p_1$, so $\Lc_1\times^h_{\G_\Sigma}\Lc_2$ is not clean, although $N_1$ and $N_2$ meet transversely when
$N_2\neq N_1$. The Lagrangian subgroupoid of $D$, in contrast, has a clean shifted self-intersection.
By Proposition \ref{RADK741QHM}, this homotopy fibre product
is Morita equivalent to $(S^1_{2\pi}\times\R)\times\R$, with the symplectic form $d\alpha\wedge d\theta$ on the
cylinder, and the Hamiltonian flow of $\alpha$ has period $2\pi$. The subgroupoids of the meridian $N_1$ and $D$ also have a clean homotopy fibre product,
by case (2) of Proposition \ref{RADK638FNS}, and it is Morita equivalent to two points.
\end{example}

\begin{proof}[\textnormal{\bfseries Proof of Theorem \ref{RADK367ESC}}]
This is Proposition \ref{RADK783SZT}.
\end{proof}

\section{The homotopy Poisson algebra of a curve}\label{RADK035TQA}

In this section, we study the homotopy Poisson algebras of curves. Let $N$ be an admissible curve in a Radko surface $(\Sigma,\sigma)$, with degeneracy locus $D$, and let $\psi$ be a tubular neighbourhood of $N$. We compute the homotopy Poisson algebra $(\mathfrak{a},\lambda^\psi)$ of \S\ref{RADK059DGU} and determine how its brackets depend on the contacts of $N$ with $D$. The main distinction is between linearization by an arbitrary $\Linf$-isomorphism and linearization by a change of tubular neighbourhood. We first compute the brackets and their cohomology, then construct a formal normal form and identify the obstruction at tangential contacts. Finally, we classify the underlying $\Linf$-algebras and apply Theorem \ref{RADK919BTK} to prove Theorem \ref{RADK420ZBG}.

\subsection{The brackets in adapted coordinates}\label{RADK531BRA}

We first reduce to a single component of $N$ and compute the brackets in coordinates adapted to it.

The normal bundle of $N$ has rank one, so $\wedge^k\nu_N=0$ for $k\ge2$. Thus $\mathfrak{a}$ is concentrated in degrees $0$ and $1$. Since $N$ is compact, it has finitely many components $N_c$, each open and closed in $N$, and $\mathfrak{a}=\bigoplus_c\mathfrak{a}_c$, $\mathfrak{a}_c=C^\infty\bigl(\wedge^\bullet\nu_{N_c}\bigr)=\mathfrak{a}_c^0\oplus\mathfrak{a}_c^1$. Here $C^\infty(\nu_{N_c})$ denotes the smooth sections of the normal bundle, as in \S\ref{RADK088RYM}.
A tubular neighbourhood of $N$ restricts to one of each $N_c$. Conversely, tubular neighbourhoods of the components assemble to one of $N$ after their domains are shrunk so that their images are disjoint. The Taylor coefficients of $Q_\psi$ are multiderivations, so they are local in each argument and vanish on
tuples meeting two different summands $\mathfrak{a}_c$. Hence $(\mathfrak{a},\lambda^\psi)$ is the direct sum of the
$(\mathfrak{a}_c,\lambda^\psi)$. We fix a component $N_c$ and write $\lambda_n=\lambda_n^\psi$ for the brackets
\eqref{RADK840EAG} on $\mathfrak{a}_c$.

Each $N_c$ is an embedded circle in an oriented surface, so its normal bundle $\nu_{N_c}$ is trivial. An \defn{adapted chart} near a point of
$N_c$ is a chart $(x,y)$ of $\Sigma$, with domain $U$, such that $N_c\cap U=\{y=0\}$. It defines a local tubular
neighbourhood $(x,y\,e(x))\mapsto(x,y)$, where $e$ is the local frame of $\nu_{N_c}$ given by the class of $\dd y$.
Conversely, a tubular neighbourhood $\psi$, a coordinate $x$ on an open subset of $N_c$ and a frame $e$ give an
adapted chart, with $y$ the fibre coordinate dual to $e$. The constructions of \S\ref{RADK059DGU} are local over $N$,
and they apply verbatim to local tubular neighbourhoods.

In the adapted chart determined by $\psi$, the formal expansion of the Poisson bivector is
\begin{align*}
\sigma_\psi&=\upsilon(x,y)\,\dd x\wedge\dd y,
&\upsilon(x,y)&=\sum_{k\ge0}\upsilon_k(x)\,y^k,\\
\upsilon_k(x)&=\frac1{k!}\frac{\partial^k\upsilon}{\partial y^k}(x,0).
\end{align*}
The term $\upsilon_k\,y^k\,\dd x\wedge\dd y$ is the component of weight $k-1$ of $\sigma_\psi$
(Definition \ref{RADK189AVW}). We identify functions on $N_c$ with their pullbacks to $\nu_{N_c}$, which do not depend on
$y$, and a section $ge$ of $\nu_{N_c}$ with the vertical vector field $g(x)\,\dd y$.

\begin{proposition}\label{RADK938LYT}
Up to graded antisymmetry, the only possibly nonzero brackets of $(\mathfrak{a}_c,\lambda^\psi)$ are those with exactly one
argument in $\mathfrak{a}_c^0$. For $n\ge1$ and $f,g_1,\dots,g_{n-1}\in C^\infty(N_c)$, they are given by
\begin{equation}\label{RADK001NSV}
\lambda_n\bigl(f,g_1e,\dots,g_{n-1}e\bigr)
=(n-1)!\,\upsilon_{n-1}\,\frac{\partial f}{\partial x}\,g_1\cdots g_{n-1}\,e.
\end{equation}
\end{proposition}

\begin{proof}
We use \eqref{RADK377FMB} in an adapted chart. A bracket $\lambda_n$ with $j$ arguments in $\mathfrak{a}_c^0$ and
$n-j$ in $\mathfrak{a}_c^1$ has degree $(n-j)+(2-n)=2-j$, so it vanishes for $j=0$ and for $j\ge3$. In particular,
$\lambda_0\in \mathfrak{a}_c^2=0$.

Let $j\in\{1,2\}$. By \eqref{RADK049HQH} for $(Q_\psi)_n$ and \eqref{RADK840EAG}, replacing an argument $e$ by $ge$
multiplies the value by $g$, up to a term with $j+1$ arguments in $\mathfrak{a}_c^0$. For $j=2$, this term vanishes by
the case $j=3$, so we may take all the arguments in $\mathfrak{a}_c^1$ equal to $e$. Since $f$ and $g$ do not depend on $y$,
\begin{align*}
\{f,g\}
&=\upsilon\biggl(\frac{\partial f}{\partial x}\frac{\partial g}{\partial y}
-\frac{\partial f}{\partial y}\frac{\partial g}{\partial x}\biggr)\\
&=0.
\end{align*}
The first formula of \eqref{RADK377FMB} therefore gives $\lambda_n(e,\dots,e,f,g)=0$. For $j=1$, the extra term vanishes by the case $j=2$, so again we may take $g_i=1$, and, since
$\{f,y\}=\upsilon\,\frac{\partial f}{\partial x}$, the second formula of \eqref{RADK377FMB} gives
\begin{align*}
\lambda_n(e,\dots,e,f)
&=(-1)^{n-1}\,\frac{\partial^{n-1}}{\partial y^{n-1}}
\biggl(\upsilon\,\frac{\partial f}{\partial x}\biggr)\Big|_{y=0}\,e\\
&=(-1)^{n-1}\,(n-1)!\,\upsilon_{n-1}\,\frac{\partial f}{\partial x}\,e.
\end{align*}
By graded antisymmetry, moving $f$, of degree $0$, across the $n-1$ arguments $e$ gives the sign $(-1)^{n-1}$, hence
\eqref{RADK001NSV}.
\end{proof}

Since $f,g_1,\dots,g_{n-1}$ are arbitrary in \eqref{RADK001NSV}, $\lambda_n=0$ if and only if $\upsilon_{n-1}=0$, for $n\ge1$. In particular, $\lambda_n=0$ for all $n\ge2$ if and only if $\upsilon_k=0$ for all $k\ge1$. The complex
$(\mathfrak{a}_c,\lambda_1)$ is $\lambda_1\colon C^\infty(N_c)\to C^\infty(\nu_{N_c})$, with
\begin{equation}\label{RADK263LMD}
\lambda_1(f)=\upsilon_0\,\frac{\partial f}{\partial x}\,e .
\end{equation}

\begin{corollary}\label{RADK130HNP}
One has $\lambda_2|_{\mathfrak{a}_c^0\otimes \mathfrak{a}_c^0}=0$, and there is a unique vector field $X$ on $N_c$ with $\lambda_2(f,v)=X(f)\,v$ for $f\in \mathfrak{a}_c^0$ and $v\in \mathfrak{a}_c^1$. In an adapted chart, $X=\upsilon_1\,\dd x$.
\end{corollary}

\begin{proof}
The first claim is part of Proposition \ref{RADK938LYT}. By \eqref{RADK001NSV}, $\lambda_2(f,ge)=\upsilon_1\frac{\partial
f}{\partial x}\,ge$ in every adapted chart. These local vector fields agree on overlaps because they are characterized by the intrinsic bracket $\lambda_2(f,v)=X(f)v$. They therefore define $X$ globally. Uniqueness holds since $\nu_{N_c}$ has a nowhere vanishing section.
\end{proof}

\subsection{Contact multiplicities}\label{RADK574NFI}

The unary bracket depends on the restriction of $\sigma$ to $N_c$, while the binary bracket also depends on the tubular neighbourhood. We first record the transformation laws that distinguish these two kinds of dependence. We then define the contact multiplicities and use them to compute the cohomology of $\lambda_1$.

Keep the tubular neighbourhood, replace the coordinate $x$ by $\theta$, with $x=x(\theta)$, and replace $e$ by $\bar e=\mu\,e$, with $\mu$ nowhere zero. Then $y=\mu\,\bar y$, with $\mu$ evaluated at $x(\theta)$, and in the chart
$(\theta,\bar y)$,
\[
\dd\theta=x'(\theta)\,\dd x+\bar y\,\tfrac{d(\mu\circ x)}{d\theta}\,\dd y,\quad
\dd{\bar y}=\mu\,\dd y,\quad
\dd\theta\wedge\dd{\bar y}=x'(\theta)\,\mu\,\dd x\wedge\dd y.
\]
Hence $\bar\upsilon(\theta,\bar y)=\upsilon(x(\theta),\mu\,\bar y)/(x'(\theta)\,\mu)$, whose coefficients of
$\bar y^0$ and $\bar y^1$ are
\begin{align}
\bar\upsilon_0(\theta)&=\frac{\upsilon_0(x(\theta))}{x'(\theta)\,\mu(x(\theta))},\label{RADK743FKC}\\
\bar\upsilon_1(\theta)&=\frac{\upsilon_1(x(\theta))}{x'(\theta)},\quad\text{so that}\quad
\bar\upsilon_1\,\dd\theta=\upsilon_1\,\dd x\quad\text{on }N_c .\label{RADK958WSC}
\end{align}
More intrinsically, there is a canonical identification
\[
TN_c\otimes\nu_{N_c}\longrightarrow\wedge^2T\Sigma|_{N_c},
\quad u\otimes[v]\longmapsto u\wedge v.
\]
It is well defined because $TN_c$ has rank one.
A tubular neighbourhood is the identity on $TN_c$ and induces the identity on $\nu_{N_c}$, so for every tubular
neighbourhood, coordinate and frame,
\begin{equation}\label{RADK716SNT}
\upsilon_0\,\dd x\otimes e=\sigma|_{N_c} .
\end{equation}
So a change of tubular neighbourhood, coordinate or frame multiplies $\upsilon_0$ by a nowhere vanishing function, as
in \eqref{RADK743FKC}, and the zeros of $\upsilon_0$ and their orders do not depend on these choices. Since $N_c$ is
a circle and $\nu_{N_c}$ is trivial, we fix a global angular coordinate $x\in\R/\Z$ on $N_c$ and a global frame $e$, and we
regard the $\upsilon_k$ as functions on $N_c$.

By \eqref{RADK722CVJ}, two tubular neighbourhoods satisfy $\sigma_{\psi'}=\varphi^*\sigma_\psi$ with
$\varphi=\varphi_{\psi\psi'}\in\DiffI(\nu_{N_c})$. By Lemma \ref{RADK768RKF}, $\varphi=e^Z$ with $Z\in\X_{\ge1}(\nu_{N_c})$, so
$\varphi^*=e^{\ad_Z}$ on formal bivectors by \eqref{RADK781EGD}. In an adapted chart, the component of weight $m$ of a
formal vector field is
\[
Z_{(m)}=y^m\,\kappa(x)\,\dd x+y^{m+1}\,\rho(x)\,\dd y,\quad\text{with }\kappa,\rho\in C^\infty(N_c).
\]

\begin{lemma}\label{RADK892KCU}
Let $m\ge1$, let $\pi=\upsilon\,\dd x\wedge\dd y$ be a formal bivector of weight $\ge-1$, let
$Z\in\X_{\ge m}(\nu_{N_c})$ with $Z_{(m)}$ as above, and let $\tilde\upsilon$ be the coefficient of
$e^{\ad_Z}\pi=\tilde\upsilon\,\dd x\wedge\dd y$. Then
\begin{align}
\tilde\upsilon_k&=\upsilon_k\quad(0\le k\le m-1),\notag\\
\tilde\upsilon_m&=\upsilon_m+\kappa\,\upsilon_0'
-\bigl(\kappa'+(m+1)\rho\bigr)\upsilon_0,
\label{RADK012BAV}
\end{align}
where the prime denotes $\frac{d}{dx}$.
\end{lemma}

\begin{proof}
Since $\pi$ has weight $\ge-1$, with component $\upsilon_0\,\dd x\wedge\dd y$ of weight $-1$, and $\ad_Z$ raises the
weight by at least $m$, the components of $e^{\ad_Z}\pi=\pi+[Z,\pi]+\frac12[Z,[Z,\pi]]+\cdots$ of weight $\le m-2$ are
those of $\pi$, and that of weight $m-1$ changes by $[Z_{(m)},\upsilon_0\,\dd x\wedge\dd y]$. By \eqref{RADK808KXY},
\begin{align*}
Z_{(m)}(\upsilon_0)&=y^m\kappa\,\upsilon_0',\\
\bigl[Z_{(m)},\dd x\bigr]&=-y^m\kappa'\,\dd x-y^{m+1}\rho'\,\dd y,\\
\bigl[Z_{(m)},\dd y\bigr]&=-m\,y^{m-1}\kappa\,\dd x-(m+1)y^m\rho\,\dd y,
\end{align*}
so, by the Leibniz rule,
\begin{align*}
\bigl[Z_{(m)},\upsilon_0\,\dd x\wedge\dd y\bigr]
&=Z_{(m)}(\upsilon_0)\,\dd x\wedge\dd y
+\upsilon_0\,\bigl[Z_{(m)},\dd x\bigr]\wedge\dd y
+\upsilon_0\,\dd x\wedge\bigl[Z_{(m)},\dd y\bigr]\\
&=y^m\big(\kappa\,\upsilon_0'-\kappa'\,\upsilon_0-(m+1)\,\rho\,\upsilon_0\big)\,\dd x\wedge\dd y .
\end{align*}
This gives \eqref{RADK012BAV}.
\end{proof}

To interpret $\upsilon_0$ and $\upsilon_1$ geometrically, we compute them in the local models of \S\S\ref{RADK261MGJ} and \ref{RADK793NZP}.

\begin{lemma}\label{RADK595PPB}
\begin{enumerate}
\item At a point $q\in N_c$, one has $\upsilon_0(q)=0$ if and only if $q\in D$. In particular, if
$N_c\subset\Sigma\setminus D$, then $\upsilon_0$ is nowhere zero.
\item If $N_c=D_j$, then in the adapted chart $(x,y)=(t,r)$ of \eqref{RADK871UDA}, $\upsilon=-y$, so that
$\upsilon_0\equiv0$, $\upsilon_1\equiv-1$ and $\upsilon_k=0$ for $k\ge2$.
\item At a transverse contact $p\in N_c\cap D_j$, in the adapted chart $(x,y)=(r,s)$ of Lemma \ref{RADK227HCV},
$\upsilon=x$, so that $\upsilon_0(x)=x$ and $\upsilon_k=0$ for $k\ge1$.
\item At a tangential contact $p\in N_c\cap D_j$, in the adapted chart $x=t$, $y=r-h(t)$ built from
Lemma \ref{RADK562YPM}, $\upsilon=-\big(y+h(x)\big)$, so that $\upsilon_0=-h$, $\upsilon_1\equiv-1$ and $\upsilon_k=0$ for
$k\ge2$.
\end{enumerate}
\end{lemma}

\begin{proof}
(1) For $q\in N_c$, $\sigma_q=\upsilon_0(q)\,\dd x\wedge\dd y$, and $\sigma_q=0$ exactly when $q\in D$. (2)
$\sigma=r\,\dd r\wedge\dd t=-r\,\dd t\wedge\dd r$. (3) $\sigma=r\,\dd r\wedge\dd s$ and $N_c\cap W=\{s=0\}$ by
Lemma \ref{RADK227HCV}. (4) $D_j\cap W=\{r=0\}$ and $N_c\cap W=\{r=h(t)\}$ by Lemma \ref{RADK562YPM}. With $x=t$ and $y=r-h(t)$,
\[ \dd x=\dd t,\quad \dd y=\dd r-h'(t)\,\dd t,\quad \dd x\wedge\dd y=\dd t\wedge\dd r, \quad \sigma=-r\,\dd t\wedge\dd r=-\bigl(y+h(x)\bigr)\,\dd x\wedge\dd y. \]
\end{proof}

\begin{definition}\label{RADK574MUL}
Let $N_c\not\subset D$ be a component of $N$, and let $p\in N_c\cap D$. The \defn{contact multiplicity}
$m_p\in\{1,2,\dots\}\cup\{\infty\}$ of $N_c$ at $p$ is the order of vanishing at $p$ of $r|_{N_c}$, where $r$ is a local defining
function of $D$ near $p$. Here $m_p=\infty$ means that $r|_{N_c}$ is flat at $p$, that is, all its derivatives vanish there.
\end{definition}

The multiplicity $m_p$ does not depend on $r$, since two defining functions differ by a nowhere vanishing
factor.

\begin{proposition}\label{RADK744CHM}
Let $N_c$ be a component of $N$ with $N_c\not\subset D$.
\begin{enumerate}
\item For every tubular neighbourhood, coordinate and frame, $\upsilon_0$ vanishes at $p\in N_c\cap D$ to order
exactly $m_p$. Moreover, $m_p=1$ if and only if the contact at $p$ is transverse.
\item $H^0(\mathfrak{a}_c,\lambda_1)=\R$, and
\[
\im\lambda_1=\Big\{\upsilon_0\,u\,e:\ u\in C^\infty(N_c),\ \oint_{N_c}u\,dx=0\Big\}.
\]
\item If every $m_p$ is finite, then $\sum_pm_p$ is even and
\[
\dim H^1(\mathfrak{a}_c,\lambda_1)=1+\sum_{p\in N_c\cap D}m_p .
\]
\item If some $m_p$ is infinite, then $H^1(\mathfrak{a}_c,\lambda_1)$ is infinite-dimensional.
\end{enumerate}
\end{proposition}

\begin{proof}
(1) Near $p\in D_j$, $\sigma=r\,\dd r\wedge\dd t$, so, by \eqref{RADK716SNT}, $\upsilon_0$ is $r|_{N_c}$ times a
nowhere vanishing function. Moreover, $r|_{N_c}$ has a simple zero at $p$ if and only if $dr|_{T_pN_c}\neq0$, that is,
$T_pN_c\neq T_pD_j$.

(2) By \eqref{RADK263LMD}, $\lambda_1(f)=\upsilon_0f'e$, and $\upsilon_0$ vanishes only on the finite set $N_c\cap D$
by Lemma \ref{RADK595PPB}(1). So $\lambda_1(f)=0$ implies $f'=0$ on a dense set, and $f$ is constant. The image follows,
since $u=f'$ for some $f\in C^\infty(N_c)$ if and only if $\oint u\,dx=0$.

(3) As in (2), $u\mapsto\upsilon_0u$ is injective on $C^\infty(N_c)$. Let $I_{\mathbf m}\subset C^\infty(N_c)$, with $\mathbf m=(m_p)_p$, be the
space of functions whose jet of order $m_p-1$ vanishes at every $p\in N_c\cap D$. Clearly
$\upsilon_0\,C^\infty(N_c)\subset I_{\mathbf m}$. Conversely, for $g\in I_{\mathbf m}$, Taylor's formula gives
$\upsilon_0=(x-x(p))^{m_p}\,w_p$ and $g=(x-x(p))^{m_p}\,g_p$ near $p$, with $w_p$ and $g_p$ smooth and $w_p(p)\neq0$.
So $g/\upsilon_0$, which is smooth away from $N_c\cap D$, equals $g_p/w_p$ near $p$, so $g\in\upsilon_0\,C^\infty(N_c)$. Polynomials in $x-x(p)$ times
bump functions show that $C^\infty(N_c)$ maps onto the jets of order $m_p-1$ at the points $p$, so $I_{\mathbf m}$ has
codimension $\sum_pm_p$. By (2), $\im\lambda_1$ is the image of the hyperplane $\{\oint u\,dx=0\}$ under
$u\mapsto\upsilon_0u\,e$, so it has codimension $1$ in $I_{\mathbf m}e$, and codimension $1+\sum_pm_p$ in
$\mathfrak{a}_c^1$. Finally, $\upsilon_0$ changes sign at $p$ if and only if $m_p$ is odd, and an even number of times
around $N_c$, so $\sum_pm_p$ is even.

(4) If $m_p=\infty$, every element of $\im\lambda_1$ is flat at $p$, so, as in (3), $H^1(\mathfrak{a}_c,\lambda_1)$
maps onto the jets of every finite order at $p$.
\end{proof}

For transverse contacts, the dimension of $H^1(\mathfrak{a}_c,\lambda_1)$ agrees with the Mazzeo--Melrose formula for $b$-manifolds \cite[Theorem 27]{gmp2014} (see also \cite[Proposition 3.19]{kirchhofflukat2022}). For contacts of equal order, it agrees with the cohomology of $b^k$-manifolds \cite[Proposition 4.3]{scott2016}.

\begin{example}\label{RADK546QHZ}
On the sphere of Example \ref{RADK748QSA}, let $g\colon\R/2\pi\Z\to[0,1)$ be smooth, with $g(\theta)=\exp(-1/\theta^2)$ for
$0<|\theta|\le\frac12$, $g(0)=0$ and $g>0$ elsewhere. Then $N=\{z=g(\theta)\}$ is a connected admissible curve which
meets $D$ only at $p=(0,0)$, with $m_p=\infty$, so $H^1(\mathfrak{a}_c,\lambda_1)$ is infinite-dimensional by Proposition \ref{RADK744CHM}(4).
\end{example}

\begin{remark}\label{RADK745VEP}
By \cite[Theorem 4]{crainic2003vanest}, for a Hausdorff source-connected Lie groupoid whose source fibres are
homologically $n$-connected, the van Est map from its differentiable cohomology to the cohomology of its Lie
algebroid is an isomorphism in degrees $\le n$ and injective in degree $n+1$. We apply this to the Lagrangian
subgroupoid $\Lc$ of Proposition \ref{RADK663CKE}. Every arrow of $\Lc$ has its source and its target in the same component of
$N$, so $\Lc|_{N_c}=\sss^{-1}(N_c)$ is open and closed in $\Lc$. It is Hausdorff, since $\Lc$ is embedded in the
Hausdorff groupoid $\G_\Sigma$, and source-connected, with Lie algebroid $TN_c^\circ$, whose Chevalley--Eilenberg
complex is isomorphic to $(\mathfrak{a}_c,\lambda_1)$ by Remark \ref{RADK135FUR}. If $N_c$ meets $D$, its source fibres are
contractible, since by Proposition \ref{RADK663CKE} they are arcs of $N_c\setminus D$ and isotropy groups isomorphic to
$\R_{>0}$ or $\R$ if $N_c\not\subset D$, and copies of $\R$ if $N_c=D_j$. So the van Est map is an isomorphism
$H^\bullet_{\mathrm{diff}}(\Lc|_{N_c})\cong H^\bullet(\mathfrak{a}_c,\lambda_1)$ in all degrees. If moreover
$N_c\not\subset D$, then, by Proposition \ref{RADK744CHM}, $\dim H^1_{\mathrm{diff}}(\Lc|_{N_c})=1+\sum_pm_p$ if every $m_p$ is
finite, and $H^1_{\mathrm{diff}}(\Lc|_{N_c})$ is infinite-dimensional otherwise. If $N_c\subset\Sigma\setminus D$,
then $\Lc|_{N_c}=N_c\times N_c$ is proper, so $H^1_{\mathrm{diff}}(\Lc|_{N_c})=0$ \cite[Proposition 1]{crainic2003vanest},
while $H^1(\mathfrak{a}_c,\lambda_1)=\R$ by Proposition \ref{RADK744CHM}(3). This is consistent with
\cite[Theorem 4]{crainic2003vanest}, which only gives injectivity in degree $1$ here, since the source fibres $N_c$
are not simply connected.
\end{remark}

\subsection{A normal form}\label{RADK536TRC}

We now construct a tubular neighbourhood in which the formal expansion $\sigma_\psi$ is affine in the fibre coordinate. By \eqref{RADK001NSV}, this eliminates all brackets of arity at least three. When all contacts are transverse, the same construction also eliminates the binary bracket. The normal form concerns the formal expansion along $N_c$. It does not assert that the smooth bivector is affine on an open neighbourhood.

Throughout this subsection, we fix a tubular neighbourhood $\psi_0$ of $N_c$, the global angular coordinate $x\in\R/\Z$ and
the global frame $e$. The frame identifies $\nu_{N_c}$ with $N_c\times\R$, and, after shrinking the domain of $\psi_0$, we work on $N_c\times(-\varepsilon,\varepsilon)$ with coordinates $(x,y)$. Formal functions along $N_c$ are then series $\sum_k f_k(x)y^k$, with $f_k\in C^\infty(N_c)$. As in \S\ref{RADK574NFI}, a change of coordinate $x=x(\theta)$ and of frame
$\bar e=\mu\,e$ replaces the coefficient of a formal bivector by $\upsilon(x(\theta),\mu\,\bar y)/(x'(\theta)\,\mu)$.
So the property that $\upsilon$ is a polynomial of degree at most one in $y$ does not depend on the coordinate and the
frame.

\begin{lemma}\label{RADK512RLZ}
\begin{enumerate}
\item For every $\varphi\in\DiffI(\nu_{N_c})$, there is a tubular neighbourhood $\psi$ of $N_c$ with $\varphi_{\psi_0\psi}=\varphi$.
\item Let $p_1,\dots,p_s\in N_c$ be distinct, and let $\psi_i$ be a tubular neighbourhood of an open arc $J_i\ni p_i$
of $N_c$. Then there is a tubular neighbourhood $\psi$ of $N_c$ which agrees with $\psi_i$ near $p_i$ for every $i$.
\end{enumerate}
\end{lemma}

\begin{proof}
Let $a,b$ be smooth functions on $N_c\times(-\varepsilon,\varepsilon)$ with $a(x,0)=b(x,0)=\frac{\partial b}{\partial
y}(x,0)=0$, and $\gamma(x,y)=\big(x+a(x,y),\,y+b(x,y)\big)$, with $x+a$ modulo $\Z$. Then $\gamma$ is the identity on
the zero section, and its differential there is $(\dot x,\dot y)\mapsto\big(\dot x+\frac{\partial a}{\partial
y}(x,0)\,\dot y,\ \dot y\big)$, so $\gamma$ induces the identity on $\nu_{N_c}$ and is a local diffeomorphism near the
zero section. It is injective near the zero section, since otherwise pairs $u_k\neq u'_k$ converging to the zero
section with $\gamma(u_k)=\gamma(u'_k)$ would have subsequences converging to $u,u'\in N_c$ with
$u=\gamma(u)=\gamma(u')=u'$, contradicting local injectivity at $u$. So $\psi_0\circ\gamma$ is a tubular
neighbourhood, and $\varphi_{\psi_0,\psi_0\circ\gamma}$ is the expansion of $\gamma$.

(1) Since $\varphi\in\DiffI(\nu_{N_c})$, $\varphi^*(x)=x+\sum_{k\ge1}a_k(x)\,y^k$ and $\varphi^*(y)=y+\sum_{k\ge2}b_k(x)\,y^k$, where
$\varphi^*(x)$ is a formal map to $\R/\Z$. By Borel's lemma, uniformly on the compact circle $N_c$, there are $a$ and
$b$ as above with these Taylor series in $y$. Let $\varphi_\gamma\in\DiffI(\nu_{N_c})$ be the expansion of $\gamma$, and set $\varphi'=\varphi_\gamma\circ\varphi^{-1}$. Then $\varphi'^*$ fixes $x$ and $y$, so every power of $\varphi'^*-\id$ vanishes on these coordinates. Hence $\log\varphi'=0$. By Lemma \ref{RADK768RKF}(2), $\varphi_\gamma=\varphi$, and we may take $\psi=\psi_0\circ\gamma$.

(2) Since $\psi_i$ and $\psi_0$ are the identity on $N_c$ and on $\nu_{N_c}$, for pairwise disjoint open arcs $J'_i\ni
p_i$ with compact closure in $J_i$, $\psi_0^{-1}\circ\psi_i(x,y)=\big(x+a_i(x,y),\,y+b_i(x,y)\big)$ near
$J'_i\times\{0\}$, with $a_i(x,0)=b_i(x,0)=\frac{\partial b_i}{\partial y}(x,0)=0$. Let $\beta_i\in C^\infty(N_c)$ be
supported in $J'_i$ and equal to $1$ near $p_i$. With $a=\sum_i\beta_ia_i$ and $b=\sum_i\beta_ib_i$, after shrinking
$\varepsilon$, $\psi=\psi_0\circ\gamma$ agrees with $\psi_i$ near $p_i$.
\end{proof}

The construction has two steps. First, we choose a tubular neighbourhood that is already affine near the tangential contacts. We then remove the remaining coefficients one at a time without changing the neighbourhood formally near those contacts. The following lemma provides the second step.

\begin{lemma}\label{RADK523SLV}
Let $N_c\not\subset D$, let $A\subset N_c$ be an open set which contains the tangential contacts of $N_c$ and whose
closure contains no transverse contact, and let $\pi=\upsilon\,\dd x\wedge\dd y$ be a formal bivector whose
coefficient of $y^0$ is the function $\upsilon_0$ of $\sigma_{\psi_0}$. Let $m\ge1$, and let $g\in C^\infty(N_c)$
vanish on $A$. Then there is $Z\in\X_{\ge m}(\nu_{N_c})$ of weight $m$, whose coefficients vanish on $A$, such that the
coefficient $\tilde\upsilon$ of $e^{\ad_Z}\pi$ satisfies $\tilde\upsilon_k=\upsilon_k$ for $k<m$,
$\tilde\upsilon_m=\upsilon_m-g$, and $\tilde\upsilon=\upsilon$ over $A$.
\end{lemma}

\begin{proof}
By Proposition \ref{RADK744CHM}(1), $\upsilon_0'(q)\neq0$ at every transverse contact $q$. Let $B$ be a union of open arcs around
the transverse contacts on which $\upsilon_0'\neq0$, with closure disjoint from that of $A$, and $O=N_c\setminus D$.
Then $A$, $B$ and $O$ cover $N_c$. Let $\chi_A+\chi_B+\chi_O=1$ be a subordinate partition of unity, and
\[
Z=y^m\kappa\,\dd x+y^{m+1}\rho\,\dd y,\quad
\kappa=-\chi_B\,g/\upsilon_0',\quad
\rho=\tfrac1{m+1}\bigl(\chi_O\,g/\upsilon_0-\kappa'\bigr).
\]
These functions are smooth, since $\upsilon_0'\neq0$ on $B$ and $\upsilon_0\neq0$ on $O$ by Lemma \ref{RADK595PPB}(1). On $A$, $\chi_B=g=0$, so $\kappa$,
$\kappa'$ and $\rho$ vanish. Since $\chi_Ag=0$,
\[
\kappa\,\upsilon_0'-\bigl(\kappa'+(m+1)\rho\bigr)\upsilon_0 =-\chi_Bg-\chi_Og =-g,
\]
and Lemma \ref{RADK892KCU} gives $\tilde\upsilon_k=\upsilon_k$ for $k<m$ and $\tilde\upsilon_m=\upsilon_m-g$. Finally,
$\ad_Z$ is local and $Z$ vanishes over $A$, so $e^{\ad_Z}\pi=\pi$ over $A$.
\end{proof}

\begin{proposition}\label{RADK357CDF}
Let $x$ be a global coordinate on $N_c$ and $e$ a global frame of $\nu_{N_c}$.
\begin{enumerate}
\item There is a tubular neighbourhood $\psi$ of $N_c$ with $\sigma_\psi=(\upsilon_0+\upsilon_1\,y)\,\dd x\wedge\dd y$.
For such a $\psi$, one has $\lambda^\psi_n=0$ for $n\ge3$, so that $(\mathfrak{a}_c,\lambda^\psi)$ is a differential graded
Poisson algebra.
\item If $N_c\not\subset D$ and $U$ is a neighbourhood of the tangential contacts of $N_c$, then $\psi$ can be chosen
with $\upsilon_1=0$ outside $U$. In particular, if all contacts of $N_c$ with $D$ are transverse, then there is a
tubular neighbourhood with $\sigma_\psi=\upsilon_0\,\dd x\wedge\dd y$, for which $\lambda^\psi_n=0$ for $n\ge2$.
\end{enumerate}
\end{proposition}

\begin{proof}
It suffices to construct the stated formal expansion, since the claims about the brackets then follow from \eqref{RADK001NSV}. If $N_c=D_j$, the chart $(t,r)$ of \eqref{RADK871UDA} gives a tubular neighbourhood with $\upsilon=-y$ in the coordinate $t$ and the frame $\dd r$, by Lemma \ref{RADK595PPB}(2). The coefficient remains affine after changing to $x$ and $e$.

Suppose henceforth that $N_c\not\subset D$, with tangential contacts $p_1,\dots,p_s$.

We first choose an affine model near the tangential contacts. By the beginning of this subsection, we may use $x$
and $e$. By Lemma \ref{RADK595PPB}(4) and Lemma \ref{RADK512RLZ}(2), we may assume that $\psi_0$ agrees near each $p_i$ with a
local tubular neighbourhood $\psi_i$ for which $\upsilon$ is affine in $y$. Let $U$ be a neighbourhood of the tangential contacts, arbitrary for (1). Choose an open set $A\subset U$, a union of small arcs around the $p_i$ on which $\psi_0=\psi_i$, whose closure contains no transverse contact, and $\beta\in
C^\infty(N_c)$ supported in $U$ with $\beta|_A=1$. Then $\upsilon_k=0$ on $A$ for $k\ge2$, since $\psi_0=\psi_i$ there. If $s=0$, then $A=\emptyset$ and $\beta=0$.

We then remove the remaining coefficients successively. Let $\pi^{(0)}=\sigma_{\psi_0}$. By Lemma \ref{RADK523SLV} with $m=1$ and $g=(1-\beta)\,\upsilon_1$, there is $Z_1$ of
weight $1$, vanishing over $A$, such that $\pi^{(1)}=e^{\ad_{Z_1}}\pi^{(0)}$ has coefficients $\upsilon_0$ and $\beta\upsilon_1$ of $y^0$ and $y^1$ and agrees with $\pi^{(0)}$ over $A$. Inductively, for $m\ge2$, Lemma \ref{RADK523SLV}
with $g$ the coefficient of $y^m$ in $\pi^{(m-1)}$, which vanishes on $A$, gives $Z_m$ of weight $m$, vanishing over
$A$, such that $\pi^{(m)}=e^{\ad_{Z_m}}\pi^{(m-1)}$ has the coefficients of $\pi^{(m-1)}$ in $y^0,\dots,y^{m-1}$, no
term in $y^m$, and agrees with $\pi^{(m-1)}$ over $A$.

Finally, we take the formal limit and realize it smoothly. Set $\varphi_M=e^{Z_1}\circ\cdots\circ e^{Z_M}$.
Then $\varphi_M^*\sigma_{\psi_0}=\pi^{(M)}$ by \eqref{RADK781EGD}. Since $Z_{M+1}$ maps $\Of$ into $\Ic^{M+1}$, so does $\varphi_{M+1}^*-\varphi_M^*=(e^{Z_{M+1}}-\id)\circ\varphi_M^*$, and the $\varphi_M^*$ converge to a continuous algebra endomorphism $\varphi^*$ of $\Of$, and $\varphi^*-\id$, like each $\varphi_M^*-\id$, maps $\Ic^w$ into $\Ic^{w+1}$. Hence $\varphi^*$ is invertible, with inverse
$\sum_{k\ge0}(\id-\varphi^*)^k$, and $\varphi\in\DiffI(\nu_{N_c})$. Similarly, $e^{Z_{M+1}}\circ\cdots\circ e^{Z_{M'}}$ converges, as
$M'\to\infty$, to $\varphi'_M\in\DiffI(\nu_{N_c})$ with $\varphi=\varphi_M\circ\varphi'_M$, and $\varphi'^*_M-\id$, hence
$\log\varphi'_M$, maps $\Ic^w$ into $\Ic^{w+M+1}$. So $\log\varphi'_M\in\X_{\ge M+1}(\nu_{N_c})$, since $f\,y^k\dd x$
and $f\,y^{k+1}\dd y$ have weight $k$. By \eqref{RADK781EGD},
$\varphi^*\sigma_{\psi_0}=e^{\ad_{\log\varphi'_M}}\pi^{(M)}$ differs from $\pi^{(M)}$ only in weights $\ge M$, that
is, in the coefficients of $y^k$ with $k\ge M+1$. Since $M$ is arbitrary,
$\varphi^*\sigma_{\psi_0}=(\upsilon_0+\beta\upsilon_1\,y)\,\dd x\wedge\dd y$. By Lemma \ref{RADK512RLZ}(1),
$\varphi=\varphi_{\psi_0\psi}$ for a tubular neighbourhood $\psi$, and $\sigma_\psi=\varphi^*\sigma_{\psi_0}$ by
\eqref{RADK722CVJ}. This proves (1) and (2), with $U=\emptyset$ for the last claim.
\end{proof}

For transverse contacts, the last claim also follows from the $b$-symplectic Gotay theorem of Geudens and Zambon
\cite[Theorem 2.13 and Proposition 2.15]{GeudensZambon2021}. Indeed, $N_c$ is then a Lagrangian $b$-coisotropic
submanifold \cite[Definition 2.2]{GeudensZambon2021}, whose model is the $b$-cotangent bundle of $(N_c,N_c\cap D)$
\cite[Example 1.17]{GeudensZambon2021}, up to a fibrewise linear change of coordinates which makes the neighbourhood induce the identity on $\nu_{N_c}$.

\begin{proposition}\label{RADK327GZZ}
Let $N_c\not\subset D$ be a component whose contacts with $D$, if any, are all transverse. Then $|N_c\cap D|=2k$ for
some $k\ge0$, and there are a tubular neighbourhood of $N_c$, a coordinate $\theta\in\R/2\pi\Z$ on $N_c$ and a frame
$e$ of $\nu_{N_c}$ for which $\lambda_n=0$ for $n\ge2$ and
\begin{equation}\label{RADK984ZHD}
\lambda_1(f)=
\begin{cases}
\dfrac{\partial f}{\partial\theta}\,e,&k=0,\\[2pt]
\sin(k\theta)\,\dfrac{\partial f}{\partial\theta}\,e,&k\ge1.
\end{cases}
\end{equation}
\end{proposition}

\begin{proof}
By Proposition \ref{RADK744CHM}(1), $m_p=1$ for every $p\in N_c\cap D$, so $|N_c\cap D|=\sum_pm_p$ is even by
Proposition \ref{RADK744CHM}(3). By \eqref{RADK716SNT}, $\upsilon_0$ depends only on the coordinate and the frame. If $k=0$, then
$\upsilon_0$ is nowhere zero by Lemma \ref{RADK595PPB}(1). Choose a diffeomorphism $\theta\mapsto x(\theta)$ from
$\R/2\pi\Z$ onto $N_c$, and let $\mu(x(\theta))=\upsilon_0(x(\theta))/x'(\theta)$. Then $\bar\upsilon_0\equiv1$ by
\eqref{RADK743FKC}.

If $k\ge1$, with $N_c\cap D=\{p_1,\dots,p_{2k}\}$ in cyclic order, choose such a diffeomorphism with
$x(i\pi/k)=x(p_i)$. By Lemma \ref{RADK595PPB}(1) and Proposition \ref{RADK744CHM}(1), $\upsilon_0\circ x$ and $\sin(k\theta)$ have
simple zeros exactly at the points $i\pi/k$. So, by Hadamard's lemma,
$\mu(x(\theta))=\upsilon_0(x(\theta))/\big(x'(\theta)\sin(k\theta)\big)$ is smooth and nowhere zero, and
$\bar\upsilon_0(\theta)=\sin(k\theta)$ by \eqref{RADK743FKC}.

By Proposition \ref{RADK357CDF}(2) for $\theta$ and $\bar e$, $\sigma_\psi=\bar\upsilon_0\,\dd\theta\wedge\dd{\bar y}$ for some
tubular neighbourhood $\psi$, and \eqref{RADK001NSV} gives the brackets.
\end{proof}

\begin{example}\label{RADK615PRL}
On the sphere of Example \ref{RADK748QSA}, let $N=\{z=z_0\}$ with $0<|z_0|<1$. Consider the adapted charts $(x,y)=(\theta,z-z_0)$ and $(x,\tilde y)$, where $\tilde y=z_0\log\bigl(1+y/z_0\bigr)$, $|y|<|z_0|$.
Since $\partial\tilde y/\partial y=z_0/(z_0+y)$, which equals $1$ on $N$, the two charts induce the same frame of the normal bundle. On a neighbourhood of $N$,
\[
\sigma=-(z_0+y)\,\dd x\wedge\dd y =-z_0\,\dd x\wedge\dd{\tilde y}.
\] In the first chart $\upsilon_1=-1$ and
$\lambda_2\neq0$, while in the second all the brackets $\lambda_n$, $n\ge2$, vanish. So, away from the tangential
contacts, the binary bracket depends on the tubular neighbourhood, which is why part (2) of
Theorem \ref{RADK420ZBG} only describes it at the tangential contacts.
\end{example}

\subsection{Tangential contacts and the modular vector field}\label{RADK547TGC}

The preceding construction removes the binary bracket when all contacts are transverse. We now identify the obstruction at a tangential contact. The vector field $X$ defined by the binary bracket in Corollary \ref{RADK130HNP} agrees there with the modular vector field. Higher contact multiplicity forces agreement of higher jets, and along a component of $D$ the two vector fields agree everywhere.

\begin{proposition}\label{RADK890FAH}
Suppose that $N_c$ has a tangential contact at $p\in D_j$, or that $N_c=D_j$ and $p\in D_j$. Then the vector
$\upsilon_1(p)\,\dd x|_p$ does not depend on the tubular neighbourhood, the coordinate or the frame, and, with
$X_{\mathrm{vol}}$ the modular vector field of \eqref{RADK432CKK},
\begin{equation}\label{RADK432CKM}
\upsilon_1(p)\,\dd x\big|_p=X_{\mathrm{vol}}(p)\neq0 .
\end{equation}
In particular, if $N_c=D_j$, then the vector field $X$ of Corollary \ref{RADK130HNP} is $X_{\mathrm{vol}}|_{D_j}$ for every
tubular neighbourhood. In both cases, $\lambda_2\neq0$ for every tubular neighbourhood of $N_c$.
\end{proposition}

\begin{proof}
By Lemma \ref{RADK595PPB}(2) and (4) and \eqref{RADK716SNT}, $\upsilon_0(p)=\upsilon_0'(p)=0$ for every tubular
neighbourhood, coordinate and frame. A change of coordinate and frame leaves $\upsilon_1\dd x$ unchanged by
\eqref{RADK958WSC}. Under a change of tubular neighbourhood, possibly only defined near $p$, Lemma \ref{RADK892KCU} with
$m=1$, whose proof is local over $N_c$, gives
\begin{align*}
\tilde\upsilon_1(p)
&=\upsilon_1(p)+\kappa(p)\upsilon_0'(p)
-\bigl(\kappa'(p)+2\rho(p)\bigr)\upsilon_0(p)\\
&=\upsilon_1(p).
\end{align*}
Since every change of data is a composite of the two, $\upsilon_1(p)\,\dd x|_p$ is invariant. In the charts of
Lemma \ref{RADK595PPB}(2) and (4), $\upsilon_1(p)=-1$ and $\dd x|_p=\dd t|_p$, since $\dd x=\dd t$ if $N_c=D_j$, and $\dd
x=\dd t+h'(t)\,\dd r$ with $h'(0)=0$ otherwise. So $\upsilon_1(p)\,\dd x|_p=-\dd t|_p=X_{\mathrm{vol}}(p)\neq0$ by
\eqref{RADK432CKK}. If $N_c=D_j$, this holds at every $p\in D_j$, so $X=\upsilon_1\dd x=X_{\mathrm{vol}}|_{D_j}$.
Finally, $\lambda_2(f,v)=X(f)\,v$ with $X=\upsilon_1\dd x$ by Corollary \ref{RADK130HNP}, and $X(p)\neq0$.
\end{proof}

A section of $T\Sigma|_{N_c}$ \emph{vanishes to order at least $k$} at $p$ if, in a local coordinate $x$ on $N_c$ centred at $p$, its components in a local frame of $T\Sigma$ are $O(x^k)$. This condition is independent of the coordinate and frame. It vanishes to infinite order if the condition holds for every $k$.

\begin{proposition}\label{RADK904JET}
Let $p$ be a tangential contact of $N_c$ with $D_j$, of multiplicity $m_p\in\{2,3,\dots\}\cup\{\infty\}$, let
$\mathrm{vol}$ be a volume form near $p$, and let $X$ be the vector field of Corollary \ref{RADK130HNP}.
\begin{enumerate}
\item If $m_p<\infty$, then, for every tubular neighbourhood of $N_c$, the section $X-X_{\mathrm{vol}}|_{N_c}$ of $T\Sigma|_{N_c}$ vanishes
to order $m_p-1$ at $p$. In particular, the jet of order $m_p-2$ of $X$ at $p$ does not depend on the tubular
neighbourhood. If $m_p=\infty$, this section vanishes to infinite order at $p$, and every finite-order jet of $X$ at $p$ is independent of the tubular neighbourhood.
\item If $m_p<\infty$, the jet of order $m_p-1$ of $X$ at $p$ depends on the tubular neighbourhood.
\end{enumerate}
\end{proposition}

\begin{proof}
(1) In the chart $(t,r)$ of \eqref{RADK871UDA} near $p=(0,0)$, Lemma \ref{RADK562YPM} gives $N_c=\{r=h(t)\}$. If $m_p$ is finite, $h$ vanishes at $0$ to order exactly $m_p$, and if $m_p=\infty$, then $h$ is flat at $0$. For $\mathrm{vol}_0=dr\wedge dt$, $X_{\mathrm{vol}_0}=-\dd t$ near $p$, as in
\eqref{RADK432CKK}. Changing the volume form adds to the modular vector field a Hamiltonian vector field
$\sigma(du)=r\big(\frac{\partial u}{\partial r}\dd t-\frac{\partial u}{\partial t}\dd r\big)$, so
$X_{\mathrm{vol}}=-\dd t+r\,Y$ for a smooth vector field $Y$.

Let $\psi_p$ be the local tubular neighbourhood of Lemma \ref{RADK595PPB}(4), with $x=t$, $y=r-h(t)$ and frame $\dd y$.
There, $\upsilon_0=-h$ and $\upsilon_1\equiv-1$, and $\dd x=\dd t+h'(t)\,\dd r$ along $N_c$. By \eqref{RADK958WSC}, we may
compute $X$ for every tubular neighbourhood $\psi$ in this coordinate and frame. Near $p$, the transition from
$\psi_p$ to $\psi$ is $e^Z$ with $Z\in\X_{\ge1}$, by the local version of Lemma \ref{RADK768RKF}, and Lemma \ref{RADK892KCU}
with $m=1$ gives
\begin{align*}
\upsilon_1^\psi
&=-1+\kappa\,\upsilon_0'-\bigl(\kappa'+2\rho\bigr)\upsilon_0\\
&=-1-\kappa h'+\bigl(\kappa'+2\rho\bigr)h.
\end{align*}
Suppose first that $m_p<\infty$. Then $\upsilon_1^\psi=-1+O(t^{m_p-1})$, and, along $N_c$,
\[
X-X_{\mathrm{vol}} =\upsilon_1^\psi\,\dd x+\dd t-hY =-h'(t)\,\dd r+O(t^{m_p-1})-h(t)Y =O(t^{m_p-1}).
\]
Here we used $h'(t)=O(t^{m_p-1})$ and $h(t)=O(t^{m_p})$. The claim on jets follows because $X_{\mathrm{vol}}|_{N_c}$ does not depend on $\psi$. In fact, this order is attained: $X$ is tangent to $N_c$, whereas the component of $X_{\mathrm{vol}}$ along $\dd r$ in the frame $(\dd x,\dd r)$ is $h'+O(t^{m_p})$, which has order exactly $m_p-1$.

If $m_p=\infty$, both $h$ and $h'$ are flat. The displayed formula for $\upsilon_1^\psi$ shows that $\upsilon_1^\psi+1$ is flat as well. Since $\dd x=\dd t+h'\dd r$ along $N_c$ and $X_{\mathrm{vol}}=-\dd t+hY$ there, the section $X-X_{\mathrm{vol}}|_{N_c}$ is flat. This proves the assertion for every finite-order jet.

(2) Fix $\psi$, a global coordinate $x$ on $N_c$ with $x(p)=0$ and a global frame, identifying a neighbourhood of the
zero section of $\nu_{N_c}$ with $N_c\times(-\varepsilon,\varepsilon)$. Let $\kappa\in C^\infty(N_c)$ with
$\kappa(p)\neq0$, and $\gamma(x,y)=(x+y\,\kappa(x),\,y)$. As in the proof of Lemma \ref{RADK512RLZ},
$\psi'=\psi\circ\gamma$ is a tubular neighbourhood, and $\varphi_{\psi\psi'}$ is the expansion of $\gamma$. Let
$Z=\log\varphi_{\psi\psi'}$. Since $\varphi_{\psi\psi'}^*(y)=y$ and $(\varphi^*_{\psi\psi'}-\id)^k$ maps $\Of$ into
$\Ic^k$, $Z(y)=0$ and $Z(x)=y\,\kappa$ modulo $\Ic^2$. So $Z_{(1)}=y\,\kappa\,\dd x$, and, by Lemma \ref{RADK892KCU}, the
coefficients of $\sigma_\psi$ and $\sigma_{\psi'}$ satisfy
\[
\tilde\upsilon_1-\upsilon_1=\kappa\,\upsilon_0'-\kappa'\,\upsilon_0.
\]
By Proposition \ref{RADK744CHM}(1), $\upsilon_0=x^{m_p}w$ near $p$, with $w(0)\neq0$. Therefore
\[
\tilde\upsilon_1-\upsilon_1 =m_p\kappa(0)w(0)\,x^{m_p-1}+O(x^{m_p}),\quad m_p\kappa(0)w(0)\neq0.
\]
Thus the jets of order $m_p-1$ of $\upsilon_1\dd x$ and $\tilde\upsilon_1\dd x$ at $p$ differ.
\end{proof}

\subsection{\texorpdfstring{$\Linf$}{L-infinity}-classification and induced isomorphisms}\label{RADK558LCL}

We now classify the underlying $\Linf$-algebras of the homotopy Poisson algebras $(\mathfrak{a}_c,\lambda^\psi)$. For components not contained in $D$, the classification depends only on the total contact multiplicity
\[
M(N_c)=\sum_{p\in N_c\cap D}m_p\in\{0,1,2,\dots\}\cup\{\infty\},
\]
where the empty sum is $0$ and the sum is $\infty$ if some $m_p$ is infinite. We then apply the criterion of Theorem \ref{RADK919BTK} to determine which isomorphisms are induced by changes of tubular neighbourhood.

Throughout this classification, $\Linf$-isomorphisms are algebraic: they need not be continuous, local or compatible with the commutative products. The classification therefore concerns the underlying $\Linf$-algebras, not isomorphisms preserving the full homotopy Poisson structure. A flat $\Linf$-structure $Q$ on a graded vector space $W$ is \defn{minimal} if $Q_1=0$. A complex $(C,d_C)$ is
viewed as the flat $\Linf$-structure on $C$ whose only nonzero Taylor coefficient is $d_C$. In particular, we
write $(\mathfrak{a}_c,\lambda_1)$ for the flat $\Pinf$-algebra whose only nonzero bracket is $\lambda_1$. Its coderivation,
which we denote by $Q^{\mathrm{lin}}$, is the flat $\Linf$-structure on $\mathfrak{a}_c[1]$ whose only nonzero Taylor
coefficient is $(Q_\psi)_1$.

By the decomposition theorem of Kontsevich \cite[\S4.5.1, Lemma]{kontsevich2003}, every flat $\Linf$-algebra $(V,Q)$
is $\Linf$-isomorphic to a direct sum $(W,Q^{\min})\oplus(C,d_C)$, where $(W,Q^{\min})$ is minimal and $(C,d_C)$ is a
contractible complex.

\begin{proposition}\label{RADK667KEG}
Let $\psi$ be a tubular neighbourhood of $N_c$, and write $Q=Q_\psi$ and $V=\mathfrak{a}_c[1]$.
\begin{enumerate}
\item If $N_c=D_j$, then $Q'_2\neq0$ for every flat $\Linf$-algebra $(V',Q')$ which is $\Linf$-isomorphic to
$(V,Q)$.
\item If $N_c\not\subset D$, there is an $\Linf$-isomorphism $F\colon(\mathfrak{a}_c,\lambda^\psi)\to(\mathfrak{a}_c,\lambda_1)$ with
identity linear part.
\item Suppose that $N_c\not\subset D$, and let $N'_{c'}\not\subset D'$ be a component of an admissible curve in a Radko surface $(\Sigma',\sigma')$,
with tubular neighbourhood $\psi'$ and $\mathfrak{a}'_{c'}=C^\infty(\wedge^\bullet\nu_{N'_{c'}})$. Then
$(\mathfrak{a}_c,\lambda^\psi)$ and $(\mathfrak{a}'_{c'},\lambda^{\psi'})$ are $\Linf$-isomorphic if and only if
$\sum_pm_p=\sum_{p'}m_{p'}$ in $\{0,1,2,\dots\}\cup\{\infty\}$.
\item If $N_c=D_j$ and $N'_{c'}=D'_{j'}$ are components of the degeneracy loci of two Radko surfaces, with tubular
neighbourhoods $\psi$ and $\psi'$, then $(\mathfrak{a}_c,\lambda^\psi)$ and $(\mathfrak{a}'_{c'},\lambda^{\psi'})$ are
$\Linf$-isomorphic.
\end{enumerate}
\end{proposition}

\begin{proof}
(1) By Lemma \ref{RADK595PPB}(2) and \eqref{RADK716SNT}, $\upsilon_0\equiv0$, so $Q_1=0$, and $Q_2\neq0$ by
Proposition \ref{RADK890FAH}. For an $\Linf$-isomorphism $F\colon(V,Q)\to(V',Q')$, \eqref{RADK170LCB} for the coalgebra morphisms
$F$ and $F^{-1}$ gives $(F^{-1})_1=F_1^{-1}$. By \eqref{RADK170LCB} and \eqref{RADK450GWZ}, the components of length
one of $F\circ Q=Q'\circ F$ on $V$ and on $S^2V$ are
\[
F_1\circ Q_1=Q'_1\circ F_1,\quad F_1\circ Q_2+F_2\circ Q_1^{(2)} =Q'_1\circ F_2+Q'_2\circ(F_1\cdot F_1),
\]
where $Q_1^{(2)}(a_1a_2)=Q_1(a_1)\,a_2+(-1)^{|a_1|}a_1\,Q_1(a_2)$ and $(F_1\cdot F_1)(a_1a_2)=F_1(a_1)\,F_1(a_2)$.
Hence $Q'_1=0$ and $Q'_2=F_1\circ Q_2\circ(F_1\cdot F_1)^{-1}\neq0$.

(2) We first show that the minimal summand in the decomposition theorem has zero brackets. We then compose with a linear isomorphism to obtain a linearization with identity linear part.

By Proposition \ref{RADK744CHM}(2), $H^0(\mathfrak{a}_c,\lambda_1)=\R\cdot1$. Choose a decomposition $V=V_H\oplus V_B\oplus
V'$ of graded vector spaces with $V_B=\im Q_1$, $\ker Q_1=V_H\oplus V_B$ and $1\in V_H$, so that
$V_H\subset\R\cdot1\oplus \mathfrak{a}_c^1$. By Proposition \ref{RADK938LYT}, $Q_n$ with $n\ge2$ vanishes unless exactly one
argument is a nonconstant function, and $1\cdot1=0$ in $S^2V$ since $1$ is odd. So $Q$ vanishes on $S(V_H)$, and the
inclusion $i\colon(V_H,0)\to(V,Q)$ is a strict quasi-isomorphism. Let $G\colon(V,Q)\to(W,Q^{\min})\oplus(C,d_C)$ be as
in the decomposition theorem, and $\mathsf q$ the projection onto $W$. Then $G'=\mathsf q\circ G\circ
i\colon(V_H,0)\to(W,Q^{\min})$ is a quasi-isomorphism between algebras with zero differential, so $G'_1$ is bijective,
$G'$ is invertible by induction on the length, and $Q^{\min}=G'\circ0\circ G'^{-1}=0$. It follows that $G_1\colon(V,Q_1)\to(W\oplus C,0\oplus d_C)$ is an isomorphism of complexes. Hence
\[
F=S(G_1^{-1})\circ G\colon(V,Q)\longrightarrow(V,Q^{\mathrm{lin}})
\]
is an $\Linf$-isomorphism with $F_1=\id$, where $S(G_1^{-1})$ is the coalgebra isomorphism induced by the linear map $G_1^{-1}$.

(3) As in (1), the linear part of an $\Linf$-isomorphism is an isomorphism of complexes, so the total multiplicities
agree by Proposition \ref{RADK744CHM}(3) and (4). Conversely, if they agree, the spaces $H^1$ have the same dimension, by
Proposition \ref{RADK744CHM}(3) in the finite case. Otherwise both dimensions are $2^{\aleph_0}$. Indeed, $\mathfrak{a}_c^1$ has
cardinality $2^{\aleph_0}$, since a continuous function is determined by its values on a countable dense set. If
$m_p=\infty$, every element of $\im\lambda_1$ is flat at $p$ (proof of Proposition \ref{RADK744CHM}(4)), so, by Borel's lemma,
$H^1(\mathfrak{a}_c,\lambda_1)$ maps onto $\R[[x]]$ for a coordinate $x$ centred at $p$, and the $\sum_ks^kx^k$,
$s\in\R$, are linearly independent by the Vandermonde determinant. By (2), it suffices to find an isomorphism of
complexes $(\mathfrak{a}_c,\lambda_1)\cong (\mathfrak{a}'_{c'},\lambda'_1)$. Choose complements $E^0$ of $\R$ in
$\mathfrak{a}_c^0$ and $E^1$ of $\im\lambda_1$ in $\mathfrak{a}_c^1$, and similarly $E'^0$, $E'^1$. The spaces $E^0$ and $E'^0$ are complements of a line in the spaces $C^\infty(N_c)$ and $C^\infty(N'_{c'})$, which are isomorphic through a diffeomorphism $N_c\to N'_{c'}$, and $\dim E^1=\dim E'^1$. So there are linear isomorphisms $L^0\colon E^0\to E'^0$ and $L^1\colon E^1\to
E'^1$. Since $\lambda_1$ and $\lambda'_1$ are injective on $E^0$ and $E'^0$, the maps $1\mapsto1$, $f\mapsto L^0(f)$,
$\lambda_1(f)\mapsto\lambda'_1(L^0(f))$ and $v\mapsto L^1(v)$, for $f\in E^0$ and $v\in E^1$, define an isomorphism of
complexes.

(4) Let $\psi_j$ be the tubular neighbourhood given by the chart $(t,r)$ of \eqref{RADK871UDA}, with $t\in
S^1_{\ell_j}$. By Corollary \ref{RADK250MSL},
$\Phi(\varphi_{\psi\psi_j})\colon(\mathfrak{a}_c,\lambda^\psi)\to(\mathfrak{a}_c,\lambda^{\psi_j})$ is an
$\Linf$-isomorphism, and, by Proposition \ref{RADK938LYT} and Lemma \ref{RADK595PPB}(2), the only nonzero brackets of
$\lambda^{\psi_j}$ in the frame $e$ given by $\dd r$ are, up to graded antisymmetry, $\lambda_2(f,ge)=-\frac{\partial
f}{\partial t}\,g\,e$. Let $u_j\colon\R/\Z\to D_j$, $\theta\mapsto\ell_j\theta$, and let
$\mathfrak{b}=C^\infty(\R/\Z)\oplus C^\infty(\R/\Z)\,e_0$ with the single nonzero bracket
$\lambda_2(f,ge_0)=-\frac{\partial f}{\partial\theta}\,g\,e_0$. The linear map $\Theta(f)=\ell_j^{-1}\,f\circ u_j$,
$\Theta(ge)=(g\circ u_j)\,e_0$ is bijective, preserves degrees, and satisfies
\begin{align*}
\lambda_2\bigl(\Theta(f),\Theta(ge)\bigr)
&=-\ell_j^{-1}\,\frac{\partial(f\circ u_j)}{\partial\theta}\,(g\circ u_j)\,e_0\\
&=-\biggl(\frac{\partial f}{\partial t}\,g\biggr)\circ u_j\,e_0\\
&=\Theta\bigl(\lambda_2(f,ge)\bigr).
\end{align*}
So $\Theta$ is a strict $\Linf$-isomorphism onto $\mathfrak{b}$, which does not depend on $\ell_j$, and similarly for
$D'_{j'}$.
\end{proof}

By Proposition \ref{RADK667KEG}(4), components of $D$ with different modular periods still have $\Linf$-isomorphic algebras. The next proposition shows that an isomorphism whose linear part respects multiplication of functions by normal sections must preserve the modular period.

\begin{proposition}\label{RADK669PER}
Let $N_c=D_j$ and $N'_{c'}=D'_{j'}$ be as in Proposition \ref{RADK667KEG}(4), and let $F\colon(\mathfrak{a}_c,\lambda^\psi)\to
(\mathfrak{a}'_{c'},\lambda^{\psi'})$ be an $\Linf$-isomorphism with $F_1(f\,v)=F_1(f)\,F_1(v)$ for $f\in\mathfrak{a}_c^0$ and
$v\in\mathfrak{a}_c^1$. Then $\ell_j=\ell'_{j'}$.
\end{proposition}

\begin{proof}
By the proof of Proposition \ref{RADK667KEG}(1), $F_1$ is bijective, preserves degrees and satisfies $Q'_2\circ(F_1\cdot
F_1)=F_1\circ Q_2$, so $\lambda'_2(F_1a,F_1b)=F_1\big(\lambda_2(a,b)\big)$, since the signs of \eqref{RADK840EAG} only
depend on the degrees. By Proposition \ref{RADK890FAH}, $\lambda_2(f,v)=X(f)\,v$ and $\lambda'_2(f',v')=X'(f')\,v'$, with
$X=X_{\mathrm{vol}}|_{D_j}$ and $X'=X'_{\mathrm{vol}}|_{D'_{j'}}$. For $f\in C^\infty(D_j)$ and
$v\in\mathfrak{a}_c^1$,
\[
X'\bigl(F_1(f)\bigr)\,F_1(v) =\lambda'_2\bigl(F_1f,F_1v\bigr) =F_1\bigl(X(f)v\bigr) =F_1\bigl(X(f)\bigr)F_1(v).
\]
Choosing $v$ with $F_1(v)$ nowhere zero gives $X'\circ F_1=F_1\circ X$ on $C^\infty(D_j)$, so the complexified
operators $X$ and $X'$ have the same eigenvalues. In the coordinate $t\in S^1_{\ell_j}$ of \eqref{RADK871UDA}, $X=-\dd
t$ by \eqref{RADK432CKK}, and $-\frac{df}{dt}=\mu f$ has a nonzero solution $f=e^{-\mu t}$ on $S^1_{\ell_j}$ if and
only if $\mu\in\frac{2\pi i}{\ell_j}\Z$. Consequently, $\frac{2\pi i}{\ell_j}\Z=\frac{2\pi i}{\ell'_{j'}}\Z$.
Since both periods are positive, equality of these subgroups implies $\ell_j=\ell'_{j'}$.
\end{proof}

It remains to identify the isomorphisms induced by changes of tubular neighbourhood. By Theorem \ref{RADK919BTK}, an isomorphism $F$ with $F_1=\id$ is induced in this way precisely when the Taylor coefficients of $\log F$ are multiderivations. We therefore describe those multiderivations explicitly in rank one. Fix a global frame $e$ of
$\nu_{N_c}$, with fibre coordinate $y$. A vector field $Y$ on $N_c$ lifts to $\nu_{N_c}\cong N_c\times\R$ as the
vector field which is horizontal for this product structure, and we write the lift again $Y$. In an adapted chart,
$Y=Y^x(x)\,\dd x$.

\begin{proposition}\label{RADK100GZF}
Let $n\ge2$. A graded symmetric map $l\colon S^n(\mathfrak{a}_c[1])\to \mathfrak{a}_c[1]$ of degree $0$ is a multiderivation of $\mathfrak{a}_c$ if
and only if there are a vector field $X_n$ on $N_c$ and a function $\rho_n\in C^\infty(N_c)$ such that, for all
$f,g_1,\dots,g_n\in C^\infty(N_c)$,
\begin{align}
l(f,g_1e,\dots,g_{n-1}e)&=g_1\cdots g_{n-1}\,X_n(f),\label{RADK667YLQ}\\
l(g_1e,\dots,g_ne)&=\Big(\sum_{i=1}^n\Big(\prod_{k\neq i}g_k\Big)X_n(g_i)+\rho_n\,g_1\cdots g_n\Big)\,e.
\label{RADK947FVY}
\end{align}
The pair $(X_n,\rho_n)$ is then unique, and $l=(\D_{Z_{(n-1)}})_n$ for the vector field of weight $n-1$
\begin{equation}\label{RADK006TFE}
Z_{(n-1)}=\frac{(-1)^{n-1}}{(n-1)!}\,y^{n-1}X_n+\frac{(-1)^n}{n!}\,\rho_n\,y^n\,\dd y .
\end{equation}
\end{proposition}

\begin{proof}
In $\mathfrak{a}_c[1]$, functions have degree $-1$ and sections of $\nu_{N_c}$ degree $0$, so $l$ vanishes on tuples
with at least two arguments in $\mathfrak{a}_c^0$, and permuting arguments gives no sign when at most one is a
function.

Let $l$ be a multiderivation. For $g,h\in C^\infty(N_c)$ and $a_2,\dots,a_n\in \mathfrak{a}_c[1]$, \eqref{RADK692LBY}
with $d=0$ gives
\begin{align*}
l(gh,a_2,\dots,a_n)&=l(g,a_2,\dots,a_n)\,h+g\,l(h,a_2,\dots,a_n),\\
l(ge,a_2,\dots,a_n)&=\pm\,l(g,a_2,\dots,a_n)\,e+g\,l(e,a_2,\dots,a_n),
\end{align*}
with the sign $+$ when all the $a_k$ are sections, while the first term of the second identity vanishes if some $a_k$
is a function. By the first identity, the map $f\mapsto l(f,e,\dots,e)$ is a derivation of $C^\infty(N_c)$, and hence defines a vector field $X_n$. Define $\rho_n$ by $l(e,\dots,e)=\rho_n e$. The second identity in each argument gives \eqref{RADK667YLQ} and \eqref{RADK947FVY}.

Conversely, since $\mathfrak{a}_c^1=C^\infty(N_c)\,e$, a graded symmetric map of degree $0$ satisfying
\eqref{RADK667YLQ} and \eqref{RADK947FVY} is determined by $(X_n,\rho_n)$, and taking $g_i=1$ recovers $(X_n,\rho_n)$
from it. For $Z_{(n-1)}$ given by \eqref{RADK006TFE}, \eqref{RADK862PHS} and \eqref{RADK914PAN} give
\begin{align*}
(\D_{Z_{(n-1)}})_n(e,\dots,e,x)&=(-1)^{n-1}\,\frac{\partial^{n-1}}{\partial y^{n-1}}
\Big(\frac{(-1)^{n-1}}{(n-1)!}\,y^{n-1}X_n^x\Big)=X_n^x,\\
(\D_{Z_{(n-1)}})_n(e,\dots,e)&=(-1)^n\,\frac{\partial^n}{\partial y^n}\Big(\frac{(-1)^n}{n!}\,\rho_n\,y^n\Big)e
=\rho_n\,e.
\end{align*}
By \S\ref{RADK059DGU}, $(\D_{Z_{(n-1)}})_n$ is a multiderivation. So, by the first part and these two identities, it satisfies \eqref{RADK667YLQ} and \eqref{RADK947FVY} for the pair $(X_n,\rho_n)$. Hence every map satisfying \eqref{RADK667YLQ} and \eqref{RADK947FVY} equals $(\D_{Z_{(n-1)}})_n$ and is a multiderivation.
\end{proof}

Let $\psi,\psi'$ be tubular neighbourhoods of $N_c$, and let $F\colon(\mathfrak{a}_c,\lambda^\psi)\to(\mathfrak{a}_c,\lambda^{\psi'})$ be an $\Linf$-isomorphism. If $F_1=\id$, then $\log F\in\CoderII(\mathfrak{a}_c)$. So, by Theorem \ref{RADK919BTK} and Proposition \ref{RADK100GZF}, $F$ is induced by a change of tubular neighbourhood if and only if $F_1=\id$ and every $(\log F)_n$, $n\ge2$, is given by \eqref{RADK667YLQ}--\eqref{RADK947FVY} for some $(X_n,\rho_n)$. In that case, $\log F=\D_Z$ by Lemma \ref{RADK332ENG}(1), and $F=e^{\D_Z}=\Phi(e^Z)$, with $Z=\sum_{n\ge2}Z_{(n-1)}$ given by \eqref{RADK006TFE}. For $n=2$,
$Z_{(1)}=-y\,X_2+\frac12\rho_2\,y^2\,\dd y$, that is, $\kappa=-X_2^x$ and $\rho=\frac12\rho_2$ in Lemma \ref{RADK892KCU},
so, by \eqref{RADK012BAV}, the coefficients $\upsilon_1$ of $\sigma_\psi$ and $\tilde\upsilon_1$ of $\sigma_{\psi'}$
satisfy
\begin{equation}\label{RADK284CLH}
\tilde\upsilon_1=\upsilon_1-X_2(\upsilon_0)+\Big(\frac{\partial X_2^x}{\partial x}-\rho_2\Big)\,\upsilon_0 .
\end{equation}
To apply the criterion to linearization, we also need to specify when the target complex is realized geometrically. We say that an $\Linf$-isomorphism
\[
F\colon(\mathfrak{a}_c,\lambda^\psi)\longrightarrow(\mathfrak{a}_c,\lambda_1)
\]
is \emph{induced by a change of tubular neighbourhood} if there is a tubular neighbourhood $\psi''$ with $\lambda^{\psi''}_n=0$ for all $n\ge2$, and if $F=\Phi(\varphi)$ for some $\varphi\in\DiffI(\psi,\psi'')$. The unary bracket is then automatically $\lambda_1$, by Remark \ref{RADK135FUR}.

\begin{corollary}\label{RADK652QDW}
Let $N_c\not\subset D$, and let $\psi$ be a tubular neighbourhood of $N_c$.
\begin{enumerate}
\item If $N_c$ has a tangential contact with $D$, then no $\Linf$-isomorphism $(\mathfrak{a}_c,\lambda^\psi)\to(\mathfrak{a}_c,\lambda_1)$, in particular not the one of Proposition \ref{RADK667KEG}(2), is induced by a change of tubular neighbourhood. More generally, $\Phi(\varphi)\circ Q_\psi\circ\Phi(\varphi)^{-1}\neq Q^{\mathrm{lin}}$ for every $\varphi\in\DiffI(\nu_{N_c})$.
\item If every contact of $N_c$ with $D$ is transverse, then $Q^{\mathrm{lin}}=Q_{\psi''}$ for some tubular neighbourhood
$\psi''$, and $\Phi(\varphi_{\psi\psi''})$ is an induced $\Linf$-isomorphism $(\mathfrak{a}_c,\lambda^\psi)\to(\mathfrak{a}_c,\lambda_1)$
with identity linear part.
\end{enumerate}
\end{corollary}

\begin{proof}
(1) By Proposition \ref{RADK890FAH}, the binary bracket of $Q_{\psi'}$ is nonzero for every tubular neighbourhood $\psi'$, so $Q^{\mathrm{lin}}\neq Q_{\psi'}$. For $\varphi\in\DiffI(\nu_{N_c})$, Proposition \ref{RADK286QVV}(1) gives $\Phi(\varphi)\circ Q_\psi\circ\Phi(\varphi)^{-1}=Q_{\varphi^*\sigma_\psi}$. With $\varphi^*\sigma_\psi=\tilde\upsilon\,\dd x\wedge\dd y$ in an adapted chart near a tangential contact $p$, Lemma \ref{RADK892KCU} with $m=1$ and $Z=\log\varphi$ gives $\tilde\upsilon_1(p)=\upsilon_1(p)\neq0$, as in the proof of Proposition \ref{RADK890FAH}. By the computation of
Proposition \ref{RADK938LYT} for the formal bivector $\varphi^*\sigma_\psi$, the binary bracket of $Q_{\varphi^*\sigma_\psi}$ is $(f,ge)\mapsto\tilde\upsilon_1\frac{\partial f}{\partial x}\,ge$ near $p$, which is nonzero.

(2) By Proposition \ref{RADK327GZZ}, $\lambda^{\psi''}_n=0$ for $n\ge2$ for some tubular neighbourhood $\psi''$, and
$(Q_{\psi''})_1=(Q_\psi)_1$ by Remark \ref{RADK135FUR}, so $Q_{\psi''}=Q^{\mathrm{lin}}$. The last claim is
Corollary \ref{RADK250MSL}.
\end{proof}

\begin{example}\label{RADK008DGT}
Let $M=(\R/\Z)\times\R$ with coordinates $(x,y)$ and the symplectic Poisson structure $\sigma=\dd x\wedge\dd y$, and let $N=\{y=0\}$. This is the local model of a component contained in $\Sigma\setminus D$. Then $\nu_N=N\times\R$, with frame $e=\dd y$. The maps $\psi(x,y)=(x,y)$ and $\psi'(x,y)=(x,y+y^2)$, for $|y|<\frac12$, are tubular
neighbourhoods, and the transition $\varphi=\varphi_{\psi\psi'}$ is the expansion of
$(x,y)\mapsto(x,y+y^2)$. With $\tilde y=y+y^2$, we have $\dd{\tilde y}=(1+2y)^{-1}\dd y$, so
\begin{align*}
\sigma_\psi&=\dd x\wedge\dd y,\\
\sigma_{\psi'}&=\frac{1}{1+2y}\,\dd x\wedge\dd y
=\sum_{k\ge0}(-2)^k y^k\,\dd x\wedge\dd y.
\end{align*}
By \eqref{RADK001NSV}, $\lambda^\psi_n=0$ for $n\ge2$, while
\[
\lambda^{\psi'}_n(f,g_1e,\dots,g_{n-1}e)=(n-1)!\,(-2)^{n-1}\,\frac{\partial f}{\partial x}\,g_1\cdots g_{n-1}\,e.
\]
So the two $\Pinf$-structures differ, although $\lambda_1^\psi=\lambda_1^{\psi'}$.

Let $Z=\log\varphi=z(y)\,\dd y$. The condition $e^Z(y)=y+y^2$ reads $y+z+\frac12zz'+\cdots=y+y^2$.
It determines the coefficients of $z=\sum_{k\ge2}z_ky^k$ successively and gives the iterative logarithm $z(y)=y^2-y^3+\frac32y^4-\cdots$. In Proposition \ref{RADK100GZF}, $X_n=0$ and $\rho_n=(-1)^nn!\,z_n$, so $(\log\Phi(\varphi))_n(e,\dots,e)=\rho_n\,e$, with $\rho_2=2$, $\rho_3=6$ and $\rho_4=36$.

By Proposition \ref{RADK100GZF} with $X_n=0$, the coderivation $\D_Z$ vanishes on every tuple containing an element of $C^\infty(N)$, and it is $C^\infty(N)$-multilinear on sections of $\nu_N$. By \eqref{RADK450GWZ}, the same holds for the component of length one of every power of $\D_Z$. Hence the coefficients $\Phi(\varphi)_n$, $n\ge2$, are multiderivations, given by \eqref{RADK947FVY} with $X_n=0$. In rank at least two, the coefficients of an induced isomorphism need not be multiderivations (Remark \ref{RADK881TQT}).
\end{example}

\begin{example}\label{RADK207TNC}
On the sphere of Example \ref{RADK748QSA}, let $N'$ be the circle of angular radius $\pi/4$ centred at the point of latitude $\pi/4$ on the meridian $\theta=0$. It passes through the north pole, it meets $D$ only at $p=(0,0)$, and away from the north pole
\begin{align*}
N'&=\biggl\{z=h(\theta),\ |\theta|<\frac\pi2\biggr\},\\
h(\theta)&=\frac{\sin^2\theta}{1+\cos^2\theta}
=\frac12\theta^2+O(\theta^4),
\end{align*}
so $m_p=2$. The meridian $N_1$ of Example \ref{RADK748QSA} has two transverse contacts with $D$, so its total multiplicity is also $2$. By Proposition \ref{RADK667KEG}(3), the homotopy Poisson algebras of $N'$ and $N_1$ are $\Linf$-isomorphic. By Proposition \ref{RADK327GZZ}, $N_1$ has a tubular neighbourhood in which $\lambda_n=0$ for $n\ge2$. By Proposition \ref{RADK890FAH} and \eqref{RADK432CKK}, $N'$ has none, and $\lambda_2(f,v)=X(f)\,v$ with $X(p)=-\dd\theta$ for every tubular neighbourhood.
\end{example}

\begin{proof}[\textnormal{\bfseries Proof of Theorem \ref{RADK420ZBG}}]
Part (1) is Proposition \ref{RADK667KEG}(2) and (3). In part (2), an induced $\Linf$-isomorphism $(\mathfrak{a}_c,\lambda^\psi)\to(\mathfrak{a}_c,\lambda_1)$ exists if all contacts are transverse, by
Corollary \ref{RADK652QDW}(2), and none exists otherwise, by Corollary \ref{RADK652QDW}(1). At a tangential contact $p$,
$\lambda^\psi_2(f,v)=X(f)\,v$ for the vector field $X$ of Corollary \ref{RADK130HNP}, and $X(p)=X_{\mathrm{vol}}(p)$ by
\eqref{RADK432CKM}.
\end{proof}

\bibliographystyle{amsplain}
\bibliography{radko}

\end{document}